\input ./style/arxiv-general.cfg
\documentclass[MSNbibl,number,citesort,seceqn,dvips]{arxbj}
\makeatletter
   \@ifpackageloaded{graphicx}{}{\usepackage{graphicx}}
\makeatother
\usepackage{accents}
\renewcommand{\mathring}[1]{\accentset{\circ}{#1}}

\volume{22}
\issue{3}
\pubyear{2016}
\firstpage{1894}
\lastpage{1936}
\doi{10.3150/15-BEJ714}
\docsubty{FLA}

\makeatletter
\def\sfrac#1#2{#1/#2}

\def\sklfrac#1#2{(#1/#2)}

\newcommand{\ee}{\mathrm{e}}
\renewcommand{\mathring}[1]{\accentset{\circ}{#1}}
\newcommand{\lleft}{\left}
\newcommand{\rrvert}{\vert}
\newcommand{\rright}{\right}
\newcommand{\rrVert}{\Vert}
\newcommand{\llvert}{\vert}
\newcommand{\llVert}{\Vert}
\renewcommand{\mid}{|}
\newcommand{\xrightarrow}[1]{\stackrel{\mathrm{#1}}{\to}}
\newcommand{\eqref}[1]{(\ref{#1})}
\newtheorem{theorem}{Theorem}[section]
\newtheorem{lem}{Lemma}[section]
\newtheorem{prop}{Proposition}[section]
\newremark{rmk}{Remark}[section]
\newremark{example}{Example}[section]
\newcommand{\MRC}{\operatorname{MRC}}
\newcommand{\vectorize}{\operatorname{vec}}
\newcommand{\Cov}{\operatorname{Cov}}
\makeatother

\begin{document}
\begin{frontmatter}

\title{Quadratic covariation estimation of an irregularly observed
semimartingale with~jumps and noise}
\runtitle{Quadratic covariation estimation}

\begin{aug}
\author[A]{\inits{Y.}\fnms{Yuta}~\snm{Koike}\corref{}\ead[label=e1]{kyuta@ism.ac.jp}}
\address[A]{Risk Analysis Research Center, The Institute of Statistical
Mathematics and CREST JST, 10-3 Midori-cho, Tachikawa, Tokyo~190-8562, Japan.
\printead{e1}}
\end{aug}

%
\received{\smonth{8} \syear{2014}}
%
\revised{\smonth{1} \syear{2015}}

%
\begin{abstract}
This paper presents a central limit theorem for a pre-averaged version
of the realized covariance estimator for the quadratic covariation of a
discretely observed semimartingale with noise. The semimartingale
possibly has jumps, while the observation times show irregularity,
non-synchronicity, and some dependence on the observed process. It is
shown that the observation times' effect on the asymptotic distribution
of the estimator is only through two characteristics: the observation
frequency and the covariance structure of the noise. This is completely
different from the case of the realized covariance in a pure
semimartingale setting.
\end{abstract}

%
\begin{keyword}
\kwd{jumps}
\kwd{microstructure noise}
\kwd{non-synchronous observations}
\kwd{quadratic covariation}
\kwd{stable limit theorem}
\kwd{time endogeneity}
\end{keyword}
\end{frontmatter}


\section{Introduction}

The quadratic covariation matrix of a semimartingale is one of the
fundamental quantities in statistics of semimartingales.
In the context of the estimation of the diffusion coefficient of an It\^
o process observed discretely in a fixed interval, limit theorems
associated with the discretized quadratic covariation play a key role,
and such research has a long history (cf.~\cite{Dohnal1987,GCJ1993}).
Furthermore, in recent years such an asymptotic theory has been applied
to measuring the covariance structure of financial assets from
high-frequency data. This was pioneered by \cite{AB1998,BNS2002}, and
has become one of the most active areas in financial econometrics. In
such a context, the discretized quadratic covariation is also called
the \textit{realized covariance}.

However, raw high frequency data typically deviates from the ideal
situation where we observe a continuous semimartingale at equidistant
times, and this motivates statisticians to develop the theory in more
complicated settings.
One topic is the treatment of measurement errors in the data. For
financial high-frequency data, such errors originate from market
microstructure noise and have attracted vast attention in the past
decade; among various studies, see, for example, \mbox{\cite
{ZMA2005,BNHLS2008,PV2009,Z2006,Reiss2011}}.
In the univariate context, central limit theorems under irregular
sampling settings have also been studied by many authors, especially
assuming the independence between the observed process and the
observation times; see, for example,~\cite{MZ2009,HJY2011,GCJ1993}.
In the multivariate case, the irregularity of the observation times
causes the non-synchronicity which makes the analysis more complicated.
The prominent works on this topic are the Fourier analysis approach of
\cite{MM2002}, the sampling design kernel method of \cite{HY2005} and
the quasi-likelihood analysis of \cite{OY2014}.
In addition, recently various approaches to deal with these issues
simultaneously have been proposed by many authors; see, for
example,~\cite{AFX2010,BNHLS2011,Bibinger2012,BHMR2014,SX2014,CPV2013}.

Another important issue is incorporating jumps into the model. In such
a situation interest is often paid to estimating the integrated
volatility and the integrated covariance matrix, that is,~the
integrated diffusion coefficient, and there are many studies on this
issue in various settings. Regarding the central limit theories, see,
for example,~Chapters~11 and 13 of \cite{JP2012} for the basic setting,
\cite{PV2009SPA} for the noise setting, \cite{Koike2014thy} for the
non-synchronous observation setting and \cite{BW2015} for the noisy and
non-synchronous observation setting.

In contrast, turning to the \textit{entire} quadratic covariation
estimation in the presence of jumps, there are fewer works.
A central limit theorem for the realized covariance of an equidistantly
observed L\'evy process has been proved in \citet{JP1998} in the
context of the analysis of the Euler scheme.
This result has been extended to general It\^o semimartingales in \citet
{J2008} as a special case of the asymptotic results on various
functionals of semimartingale increments.
The situation where measurement errors are present has been treated by
Jacod, Podolskij and Vetter \citet{JPV2010} who focus on the ``pre-averaging'' counterparts of the
functionals discussed in \cite{J2008}, which were introduced in \citet
{PV2009} to extend classical power variation based methods to a noisy
observation setting.
The theory requires a different treatment in the absence of the
diffusion coefficient, and this case has been studied in \citet{DJT2013}.

When we further focus on the situation where the observation times are
irregular, at least to the best of the
author's knowledge, there is no comprehensive study on the central
limit theory for the quadratic covariation estimation,
except for the recent work of Bibinger and Vetter \cite{BibV2013} and
Bibinger and Winkelmann \cite{BW2015}; the
former have derived central limit theorems for the realized covariance
and the Hayashi--Yoshida estimator of Hayashi and Yoshida
\cite{HY2005} for a general It\^o
semimartingale observed irregularly and non-synchronously, while the
latter have established a central limit theorem for an adjusted version
of the spectral covariance estimator of Bibinger and Rei{\ss} \cite{BR2014} in a
non-synchronous and noisy observation setup, focusing on asymptotically
regular observation times in the sense that they satisfy conditions in
Proposition 2.54 of \cite{MZ2012}.
The aim of this study is to develop such a theory in the situation
where the observation data is contaminated by noise and the observation
times are as general as possible. More precisely, we derive a central
limit theorem for the pre-averaged version of the realized covariance
proposed in \citet{CKP2010} (called the \textit{modulated realized
covariance}) under an irregular sampling setting in the presence of
jumps. The main finding of this paper is that in the synchronous case
the observation times' effect on the asymptotic distribution of the
estimator is \textit{only} through their conditional expected
durations, provided that the limit of such quantities are well-defined.
In other words, the irregularity of the observation times has \textit
{no} impact on the asymptotic distribution of the estimator because the
conditional expected durations of the observation times naturally link
with the magnitude of the observation frequency, and thus their effect
is not due to the irregularity. This is completely different from the
pure semimartingale setting of \cite{BibV2013} where the distribution
of the durations around the jump times of the semimartingale directly
affects the asymptotic distribution of the realized covariance.

To deal with non-synchronous observations we rely on a data
synchronization method proposed in \citet{AFX2010}, which also matches
the proposal of Section~3.6 of \cite{CKP2010}. The non-synchronicity
naturally links with the covariance structure of the noise, hence it
affects the asymptotic distribution through that relation. 
On the other hand, the interpolations to the synchronized sampling
times do not matter asymptotically.
This can be seen as a counterpart of the finding of \citet
{Bibinger2012} in the continuous case.

Another issue we attempt to solve is how the dependence between the
observed process and the observation times (called the \textit{time
endogeneity}) affects the asymptotic theory in our setting. This issue
has recently been highlighted by several authors such as Fukasawa \cite{Fu2010b},
Li \textit{et al.} \cite{LMRZZ2014},
Li, Zhang and Zheng \cite{LZZ2012},
Rosenbaum and Tankov \cite{RT2011} in various settings, and it is
indeed known that such dependence possibly causes a non-standard limit
theorem even in the continuous semimartingale setting. In this paper,
this issue is partly solved in the sense that we do not rule out the
dependence between the continuous component of the process and the
observation times, but partly rule out the dependence between the jump
component and the observation times. The result shows that the time
endogeneity is also immaterial in our setting.


This paper is organized as follows.
Section~\ref{setting} presents the mathematical model and the
construction of the estimator we are focusing on.
Section~\ref{main} is devoted to the main result of this paper.
Section~\ref{examples} provides some illustrative examples of the
observation times,
while Section~\ref{simulation} provides a simulation study.
All proofs are given in Section~\ref{proofs}.


\section{The set up}\label{setting}

Given a stochastic basis $\mathcal{B}^{(0)}=(\Omega^{(0)},\mathcal
{F}^{(0)},(\mathcal{F}^{(0)}_t)_{t\geq0} ,P^{(0)})$, we consider a
$d$-dimensional semimartingale $X=(X_t)_{t\in\mathbb{R}_+}$ of the form
\[
X_t=X_0+\int_0^t
b_s\,\mathrm{d}s+\int_0^t
\sigma_s\,\mathrm{d}W_s+(\delta1_{\{\llVert \delta\rrVert \leq1\}})\star(
\mu-\nu)_t+(\delta1_{\{\llVert \delta\rrVert >1\}})\star\mu_t,
\]
where $W$ is a $d'$-dimensional $(\mathcal{F}^{(0)}_t)$-standard
Brownian motion, $\mu$ is an $(\mathcal{F}^{(0)}_t)$-Poisson random
measure on $\mathbb{R}_+\times E$ with $E$ being a Polish space, $\nu$
is the intensity measure of $\mu$ of the form $\nu(\mathrm{d}t,\mathrm
{d}z)=\mathrm{d}t\otimes\lambda(\mathrm{d}z)$ with $\lambda$ being a
$\sigma$-finite measure on $E$, $b$ is an $(\mathcal
{F}^{(0)}_t)$-progressively measurable $\mathbb{R}^d$-valued process,
$\sigma$ is an $(\mathcal{F}^{(0)}_t)$-progressively measurable $\mathbb
{R}^d\otimes\mathbb{R}^{d'}$-valued process, and $\delta$ is an
$(\mathcal{F}^{(0)}_t)$-predictable $\mathbb{R}^d$-valued function on
$\Omega^{(0)}\times\mathbb{R}_+\times E$. Also, $\star$ denotes the
integral (either stochastic or ordinary) with respect to some
(integer-valued) random measure. Here and below, we use standard
concepts and notation in stochastic calculus, which are described in
detail in, for example,~Chapter~2 of \cite{JP2012}. Our aim is to
estimate the quadratic covariation matrix process
$[X,X]=([X^k,X^l])_{1\leq k,l\leq d}$ of $X$ from noisy and discrete
observation data of $X$.

The observed process $Y$ is subject to additional measurement errors as follows:
\[
Y_t=X_t+\epsilon_t.
\]
The mathematical construction of the noise process $\epsilon$ is
explained later. We observe the components of the $d$-dimensional
process $Y=(Y^1,\dots,Y^d)$ discretely and non-synchronously. For each
$k=1,\dots,d$ the observation times for $Y^k$ are denoted by
$t^k_0,t^k_1,\dots,$ that is, the observation data $(Y^k_{t^k_i})_{i\in
\mathbb{Z}_+}$ is available.
We assume\vspace*{1pt} that $(t^k_i)_{i=0}^\infty$ is a sequence of $(\mathcal
{F}^{(0)}_t)$-stopping times which implicitly depend on a parameter
$n\in\mathbb{N}$ representing the observation frequency and satisfy
that $t^k_i\uparrow\infty$ as $i\to\infty$ and $\sup_{i\geq
0}(t^k_i\wedge t-t^k_{i-1}\wedge t)\to^p0$ as $n\to\infty$ for any $t\in
\mathbb{R}_+$, with setting $t^k_{-1}=0$ for notational convenience
(hereafter we will refer to such a sequence as a \textit{sampling
scheme} for short).

Now we introduce the precise definition of the noise process $\epsilon
$. It is basically the same as the one from Chapter~16 of \cite
{JP2012}, but we need a slight modification to ensure the (joint)\vspace*{1pt}
measurability of the process $\epsilon$, which is necessary for us to
consider variables\vspace*{-1pt} such as $\epsilon^k_{t^k_i}$. For any $t\in\mathbb
{R}_+$ there is a transition probability $Q_t(\omega^{(0)},\mathrm
{d}u)$ from $(\Omega^{(0)},\mathcal{F}^{(0)}_t)$ into $\mathbb{R}^d$
satisfying $\int u Q_t(\omega^{(0)},\mathrm{d}u)=0$ (this will
correspond to the conditional distribution of the noise at the time $t$
given $\mathcal{F}^{(0)}_t$). Then, at each frequency $n\in\mathbb{N}$,
the stochastic basis $\mathcal{B}=(\Omega,\mathcal{F},(\mathcal
{F}_t)_{t\in\mathbb{R}_+} ,P)$ supporting the observed process $Y$ is
constructed in the following manner (for notational simplicity we
subtract the index $n$ from $\mathcal{B}$): we endow the space $\Omega
^{(1)}=(\mathbb{R}^d)^{\mathbb{N}}$ with the product Borel $\sigma
$-field $\mathcal{F}^{(1)}$ and with the probability measure $Q(\omega
^{(0)},\mathrm{d}\omega^{(1)})$ which is the product
$\bigotimes_{i\in\mathbb{N}}Q_{\mathcal{T}^n_i(\omega^{(0)})}(\omega
^{(0)},\cdot)$.
Here, $(\mathcal{T}^n_i)_{i\geq0}$ is the increasing reordering of
total observation times $\{t^k_i\dvt k=1,\dots,d$ and $i\in\mathbb
{Z}_+\}$. More formally, it is defined sequentially by $\mathcal
{T}^n_0=\min_{k=1,\dots,d}t^k_0$ and $\mathcal{T}^n_i=\min_{k=1,\dots
,d}\min\{t^k_j\dvt t^k_j>\mathcal{T}^n_{i-1}\}$ for $i=1,2\dots.$ Note
that $\mathcal{T}^n_i$ is an $(\mathcal{F}^{(0)}_t)$-stopping time since
$\mathcal{T}^n_i=\min_{k=1,\dots,d}\inf_{j\geq1} (t^k_i )_{\{
t^k_j>\mathcal{T}^n_{i-1}\}}$,
where for an $(\mathcal{F}^{(0)}_t)$-stopping time $\tau$ and a set
$A\in\mathcal{F}^{(0)}_\tau$, we define $\tau_A$ by $\tau_A(\omega
^{(0)})=\tau(\omega^{(0)})$ if $\omega^{(0)}\in A$; $\tau_A(\omega
^{(0)})=\infty$ otherwise (see Claim~1.15 from Chapter~1 of \cite{JS}).
Then we define the probability space $(\Omega,\mathcal{F},P)$ by
%
\begin{eqnarray}\label{defspace}
\Omega &=&\Omega^{(0)}\times\Omega^{(1)},\qquad
\mathcal{F}=\mathcal{F}^{(0)}\otimes\mathcal{F}^{(1)},
\nonumber\\[-8pt]\\[-8pt]\nonumber
P \bigl(\mathrm{d}\omega^{(0)},\mathrm{d}\omega^{(1)}
\bigr)&=&P^{(0)}\bigl(\mathrm{d}\omega^{(0)}\bigr)Q\bigl(
\omega^{(0)},\mathrm{d}\omega^{(1)}\bigr).
\end{eqnarray}
After that, the noise process $\epsilon=(\epsilon_t)_{t\geq0}$ is
defined on this probability space by $\epsilon_t=\epsilon^0_{\mathsf
{N}_n(t)}$, where $(\epsilon^0_i)_{i\in\mathbb{N}}$ denotes the
canonical process on $(\Omega^{(1)},\mathcal{F}^{(1)})$ and $\mathsf
{N}_n(t)=\sum_{i=0}^\infty1_{\{\mathcal{T}^n_i\leq t\}}$.
By construction, given $\mathcal{F}^{(0)}$, $(\epsilon_{\mathcal
{T}^n_i})_{i\in\mathbb{Z}_+}$ is (serially) independent and $\epsilon
_{\mathcal{T}^n_i}$ obeys the law $Q_{\mathcal{T}^n_i(\omega
^{(0)})}(\omega^{(0)},\cdot)$ for every $i$.
Finally, the filtration $(\mathcal{F}_t)_{t\geq0}$ is defined as the
one generated by $(\mathcal{F}^{(0)}_t)_{t\geq0}$ and $(\epsilon
_t)_{t\geq0}$.

Any variable or process defined on either $\Omega^{(0)}$ or $\Omega
^{(1)}$ is considered in the usual way as a variable or a process on
$\Omega$. Specifically,\vspace*{1pt} our noisy process $Y=(Y_t)_{t\geq0}$ is the
process defined as the sum of the latent process $X$ on $\Omega^{(0)}$
and the noise process $\epsilon$ on $\Omega$.

%
\begin{rmk}\label{remarkmeasure}
To ensure that the probability measure $P$ in \eqref{defspace} is
well-defined, we further need the measurability of the map $\omega
^{(0)}\mapsto Q(\omega^{(0)},A)$ for any Borel subset $A$ of $\mathbb
{R}^d$. This is ensured by the progressive measurability of the process
$(Q_t(\cdot,A))_{t\geq0}$ which we will assume later (see assumption~\textup{[A4]}). This assumption also ensures that $\mathcal{B}$ is the
\textit{very good} filtered extension of~$\mathcal{B}^{(0)}$, which is
necessary to apply the version of Jacod's stable limit theorem
described by Theorem~2.2.15 of \cite{JP2012}.
\end{rmk}


To deal with the non-synchronicity of the observation times we rely on
a data synchronization method, which is commonly used in the
literature; see, for example,~\cite
{AFX2010,BNHLS2011,CKP2010,Zhang2011}. 
Let $(T_p)_{p=0}^\infty$ and $(\tau^k_p)_{p=0}^\infty$ $(k=1,\dots,d)$
be sampling schemes such that
%
\begin{equation}
\label{H1} \tau^k_0\leq T_0\quad
\mbox{and}\quad T_{p-1}<\tau^k_p\leq
T_p\qquad\mbox{for any $p\geq1$ and any $k=1,\dots,d$.}
\end{equation}
We\vspace*{-1pt} assume that the observation data $(Y^k_{\tau^k_p})_{p\in\mathbb
{Z}_+}$ is available for every $k=1,\dots,d$, that is, $\{\tau^k_p\dvt
p\geq0\}\subset\{t^k_i\dvt i\geq0\}$, and construct statistics based on
this synchronized\vspace*{1pt} data set $(Y^k_{\tau^k_p})_{p\in\mathbb{Z}_+}$,
$k=1,\dots,d$. In \citet{AFX2010} this type of synchronization method
is called the \textit{Generalized Synchronization method} and
$(T_p)_{p=0}^\infty$ is called the \textit{Generalized Sampling Time}.
One way to implement such synchronization is the so-called \textit
{refresh time sampling method} introduced by Barndorff-Nielsen \textit{et al.} \citet{BNHLS2011} to this
area. Namely, we first define the \textit{refresh times} $T_0,T_1,\dots
$ of the\vspace*{1pt} sampling schemes $\{(t^k_i)\}_{k=1}^d$ sequentially by
$T_0=\max\{t^1_0,\dots,t^d_0\}$ and
$T_p=\max_{k=1,\dots,d}\min\{t^k_i\dvt t^k_i>T_{p-1}\}$ for $p=1,2,\dots
.$ After that, for each $k$, $(\tau^k_p)$ is defined by interpolating
the next-ticks into $(T_p)$ as follows:
\begin{eqnarray*}
\tau^k_0=t^k_0\quad\mbox{and}\quad\tau^k_p=\min\bigl\{t^k_i\dvt t^k_i>T_{p-1}
\bigr\},\qquad p=1,2,\dots.
\end{eqnarray*}
Note that $\tau^k_p$ is an $(\mathcal{F}^{(0)}_t)$-stopping time due to
an analogous reason to that for $\mathcal{T}^n_i$.

Now the modulated realized covariance (henceforth MRC) estimator we
focus on is constructed in the following way. First, we choose a
sequence $k_n$ of positive integers and a number $\theta\in(0,\infty)$
such that
$k_n=\theta\sqrt{n} +\mathrm{o}(n^{1/4})$
as $n\to\infty$. We also choose a continuous \mbox{function} $g\dvtx
[0,1]\rightarrow\mathbb{R}$ which is piecewise $C^1$ with a piecewise
Lipschitz derivative $g'$ and satisfies $g(0)=g(1)=0$ and $\int_0^1
g(x)^2\,\mathrm{d}x>0$.
After that, for any $d$-dimensional stochastic process $V=(V^1,\dots
,V^d)$ we define the quantity
%
\begin{equation}
\label{defPA} \overline{V}^k_{i}=\sum
_{p=1}^{k_n-1}g \biggl(\frac{p}{k_n} \biggr)
\bigl(V^k_{\tau^k_{i+p}}-V^k_{\tau^k_{i+p-1}} \bigr),
\end{equation}
and set $\overline{V}_i=(\overline{V}^1_i,\dots,\overline{V}^d_i)^*$
(hereafter an asterisk denotes the transpose of a matrix). The MRC
estimator is defined by
\begin{eqnarray*}
\MRC[Y]^n_t=\frac{1}{\psi_2 k_n}\sum
_{i=0}^{N^n_t-k_n+1}\overline{Y}_{i} (
\overline{Y}_{i} )^* -\frac{\psi_1}{2\psi_2k_n^2}[Y,Y]^n_t,
\end{eqnarray*}
where $N^n_t=\max\{p\dvt T_p\leq t\}$, $\psi_1=\int_0^1g'(x)^2\,\mathrm{d}x$, $\psi_2=\int_0^1g(x)^2\,\mathrm{d}x$ and
\begin{eqnarray*}
[Y,Y]^n_t=\sum_{p=1}^{N^n_t}
\Delta_{p} Y (\Delta_{p} Y )^*,\qquad
\Delta_{p} Y= \bigl(Y^1_{\tau^1_p}-Y^1_{\tau^1_{p-1}},
\dots,Y^d_{\tau^d_p}-Y^d_{\tau^d_{p-1}} \bigr)^*
\end{eqnarray*}
for each $t\in\mathbb{R}_+$. Here, we set $\sum_{i=p}^q\equiv0$ if
$p>q$ by convention.
In the synchronous and equidistant sampling setting, the asymptotic
distribution of the MRC estimator has been derived in \citet{JPV2010}
(see Theorem 4.6 of that paper, and see also Section~4 of \citet{HP2012}).
Our purpose is to develop the asymptotic distribution of the MRC
estimator in the situation where the observation times are possibly
irregular, non-synchronous and endogenous.


\section{Main result}\label{main}

\subsection{Notation}

In this subsection, some notation is introduced in order to state our
main result. 
{First} we introduce notation appearing { in} the assumptions stated in
the next subsection. We write
$X^n\xrightarrow{ucp}X$ for processes $X^n$ and $X$ to express\vspace*{1pt} shortly
that $\sup_{0\leq t\leq T}\llvert X^n_t-X_t\rrvert \rightarrow^p0$
{for any $T>0$}. $\varpi$ denotes some (fixed) positive constant.
{We denote by $(\mathcal{G}^{(0)}_t)$ (resp.,~$(\mathcal{G}_t)$) the
smallest filtration containing $(\mathcal{F}^{(0)}_t)$
(resp.,~$(\mathcal{F}_t)$) such that $\mathcal{G}^{(0)}_0$
(resp.,~$\mathcal{G}_0$) contains the \mbox{$\sigma$-}field generated by $\mu
$, that is,~the \mbox{$\sigma$-}field generated by all the variables $\mu(A)$,
where $A$ ranges all measurable subsets of $\mathbb{R}_+\times E$.}

{Next,} we introduce some quantities appearing in the representation of
the asymptotic variance of the estimator. We set $\Sigma_s=\sigma
_s\sigma_s^*$ for each $s\in\mathbb{R}_+$, that is,~$\Sigma$ denotes
the diffusion coefficient matrix process. We denote by $\Upsilon_t$ the
covariance matrix of $\epsilon_t$, that is,~$\Upsilon_t(\cdot)=\int
uu^*Q_t(\cdot,\mathrm{d}u)$ (we will assume the existence of the second
moment of the noise later, so this matrix always exists). For any
real-valued bounded measurable functions $u,v$ on $[0,1]$, we define
the function $\phi_{u,v}$ on $[0,1]$ by $\phi_{u,v}(y)=\int_{y}^1
u(x-y)v(x)\,\mathrm{d}x$. Then, we put
\begin{eqnarray*}
\Phi_{22}=\int_{0}^{1}
\phi_{g,g}(y)^2\,\mathrm{d}y,\qquad\Phi_{12}=\int
_{0}^{1}\phi_{g,g}(y)
\phi_{g',g'}(y)\,\mathrm{d}y,\qquad\Phi_{11}=\int
_{0}^{1}\phi_{g',g'}(y)^2
\,\mathrm{d}y.
\end{eqnarray*}
On the other hand, for any $k,l=1,\dots,d$ we define the process
$\mathfrak{J}^{kl}$ by
\begin{eqnarray*}
\mathfrak{J}^{kl}_s=\Delta X^k_s
\Delta X^{l}_s \biggl\{\Phi_{22}\theta\bigl(
\Sigma^{kl}_{s-}G_{s-}+\Sigma^{kl}_{s}G_{s}
\bigr)+\frac{\Phi_{12}}{\theta} \bigl(\Upsilon^{kl}_{s-}
\chi^{kl}_{s-}+\Upsilon^{kl}_{s}
\chi^{kl}_{s} \bigr) \biggr\}.
\end{eqnarray*}

%
\begin{rmk}[(Properties of $\phi_{u,v}$)]\label{rmkphi}
We will use the following properties of $\phi_{u,v}$: first, for any
real-valued bounded measurable function $u$ on $[0,1]$, $\phi_{u,u}$ is
non-negative. In fact, setting $u(x)=0$ for $x\notin[0,1]$, we have
$\phi_{u,u}(y)=\int_{-\infty}^\infty u(-(y-x))u(x)\,\mathrm{d}x$ for all
$y\in[0,1]$, and we can extend the domain of $\phi_{u,u}$ to the whole
real line using this expression. Then, denoting by $\hat{f}$ the
Fourier transform of a function $f$ on $\mathbb{R}$, we have $\hat{\phi
}_{u,u}=\llvert \hat{u}\rrvert ^2\geq0$. Hence, $\phi_{u,u}$ is a
positive definite function and, in particular, $\phi_{u,u}(y)\geq0$ for
all $y\in\mathbb{R}$. Next, we can easily check that $\phi_{g,g}'=\phi
_{g,g'}=-\phi_{g',g}$ and $\phi_{g,g}''=-\phi_{g',g'}$. In particular,
$\Phi_{12}=\int_0^1\phi_{g',g}(y)^2\,\mathrm{d}y$ due to integration by parts.
\end{rmk}

\subsection{Assumptions}

We impose the following condition on the sampling schemes $(T_p)_{p\geq
0}$ and $(\tau^k_p)_{p\geq0}$ ($k=1,\dots,d$):
%
\begin{longlist}[{[A1]}]
\item[{[A1]}]  $(T_p)_{p\geq0}$\vspace*{1pt} and $(\tau^k_p)_{p\geq0}$
($k=1,\dots,d$) are sequences of {$(\mathcal{F}^{(0)}_t)$}-stopping
times and satisfy~$(\ref{H1})$. It also holds that
%
\begin{equation}
\label{A4} r_n(t):=\sup_{p\geq0}(T_p
\wedge t-T_{p-1}\wedge t)=\mathrm{o}_p\bigl(n^{-\xi}\bigr)
\end{equation}
as $n\to\infty$ (note that $T_{-1}=0$ by convention) for every $t>0$
and every $\xi\in(0,1)$. Moreover, for each $n$ we have {a $(\mathcal
{G}^{{(0)}}_t)$}-progressively measurable positive-valued process
$G^n_t$, {a $(\mathcal{G}^{(0)}_t)$}-progressively measurable
$[0,1]^d\otimes[0,1]^d$-valued process $\chi^n_t=(\chi^{n,kl}_t)_{1\leq
k,l\leq d}$ and a random subset $\mathcal{N}^n$ of $\mathbb{Z}_+$
satisfying the following conditions:

%
\begin{enumerate}[(iii)]
\item[(i)] $\{(\omega,p)\in\Omega\times\mathbb{Z}_+:p\in\mathcal{N}^n(\omega
)\}$ is a measurable set of $\Omega\times\mathbb{Z}_+$. Moreover, there
is a constant $\kappa\in(0,\frac{1}{2})$ such that $\#(\mathcal{N}^n\cap
\{p\dvt T_p\leq t\})=\mathrm{O}_p(n^\kappa)$ as $n\to\infty$ for every $t>0$.

\item[(ii)] $E[n(T_{p+1}-T_p)\mid{\mathcal{G}^{(0)}_{T_p}}]=G^n_{T_p}$ and
$E[1_{\{\tau^k_{p+1}=\tau^l_{p+1}\}}\mid {\mathcal
{G}^{(0)}_{T_p}}]=\chi^{n,kl}_{T_p}$ for every $n$, every $p\in\mathbb
{Z}_+-\mathcal{N}^n$ and any $k,l=1,\dots,d$.

\item[(iii)]  There is a cadlag $(\mathcal
{F}^{(0)}_t)$-adapted positive valued process $G$ such that:
%
\begin{longlist}[(iii-a)]

\item[(iii-a)]$n^\varpi(G^n-G)\xrightarrow{ucp}0$,

\item[(iii-b)]$G_{t-}>0$ for every $t>0$,

\item[(iii-c)]$G$ is an It\^o semimartingale of the form
\begin{eqnarray*}
G_t=G_0+\int_0^t
\widehat{b}_s\,\mathrm{d}s+\int_0^t
\widehat{\sigma}_s\,\mathrm{d}W_s +(\widehat{
\delta}1_{\{\llvert \widehat{\delta}\rrvert \leq1\}})\star(\mu-\nu
)_t+(\widehat{
\delta}1_{\{\llvert \widehat{\delta}\rrvert >1\}})\star\mu_t,
\end{eqnarray*}
where $\widehat{b}_s$ is a locally bounded and $(\mathcal
{F}^{(0)}_t)$-progressively measurable real-valued process, $\widehat
{\sigma}_s$ is a cadlag $(\mathcal{F}^{(0)}_t)$-adapted $\mathbb
{R}\otimes\mathbb{R}^{d'}$-valued process, and $\widehat{\delta}$ is an
$(\mathcal{F}^{(0)}_t)$-predictable real-valued function on $\Omega
^{(0)}\times\mathbb{R}_+\times E$ such that there is a sequence
$(\widehat{\rho}_j)$ of $(\mathcal{F}^{(0)}_t)$-stopping times
increasing to infinity and, for each $j$, a deterministic non-negative
function $\widehat{\gamma}_j$ on $E$ satisfying $\int\widehat{\gamma
}_j(z)^2\wedge1\lambda(\mathrm{d}z)<\infty$ and $\llvert \widehat
{\delta}(\omega^{(0)},t,z)\rrvert \leq\widehat{\gamma}_j(z)$ for all
$(\omega^{(0)},t,z)$ with $t\leq\widehat{\rho}_j(\omega^{(0)})$.
\end{longlist}

\item[(iv)] There is a cadlag $(\mathcal{F}^{(0)}_t)$-adapted $[0,1]^d\otimes
[0,1]^d$-valued process $\chi$ such that $n^\varpi(\chi^n-\chi
)\xrightarrow{ucp}0$ as $n\to\infty$. Furthermore,\vspace*{1pt} for each $j\in\mathbb
{N}$ we have a cadlag $(\mathcal{F}^{(0)}_t)$-adapted $[0,1]^d\otimes
[0,1]^d$-valued process\vspace*{1pt} $\chi(j)$, an $(\mathcal{F}^{(0)}_t)$-stopping
time $\check{\rho}_j$, and a constant $\check{\Lambda}_j$ such that
$\check{\rho}_j\uparrow\infty$ as $j\to\infty$ and $\chi(\omega
^{(0)})_t=\chi(j)(\omega^{(0)})_t$ if $t<\check{\rho}_j(\omega^{(0)})$ and
\begin{eqnarray*}
E \bigl[\bigl\llVert\chi(j)_{t_1}-\chi(j)_{t_2}\bigr\rrVert
^2\mid\mathcal{F}_{t_1\wedge t_2} \bigr] \leq\check{
\Lambda}_j E \bigl[\llvert t_1-t_2\rrvert
^\varpi\mid\mathcal{F}_{t_1\wedge t_2} \bigr]
\end{eqnarray*}
for every $j$ and any $(\mathcal{F}^{(0)}_t)$-stopping times $t_1$ and
$t_2$ bounded by $j$.
\end{enumerate}
\end{longlist}

%
\begin{rmk}\label{rmkA1}
(i) The assumptions on $(T_p)$ are motivated by the concept of the
\textit{restricted discretization scheme} discussed in detail in
Chapter~14 of \cite{JP2012}. In fact, suppose that $T_p$'s are of the form
\[
T_p=T_{p-1}+\theta^n_{T_{p-1}}
\varepsilon(n,p),\qquad p=1,2,\dots,
\]
where $\theta^n$ is a cadlag $(\mathcal{F}^{(0)}_t)$-adapted process,
$(\varepsilon(n,p))_{p\geq1}$ is a sequence of i.i.d.~positive
variables independent of $b$, $\sigma$, $\delta$, $W$, $\mu$, and such
that $E[\varepsilon(n,p)]=1$ and $E[\varepsilon(n,p)^r]<\infty$ for
every $r>0$, and $T_0=0$.
By constructing the filtration $(\mathcal{F}^{(0)}_t)$ suitably, we may
assume that $\epsilon(n,p)$ is independent of $\mathcal
{F}^{(0)}_{T_{p-1}}$ for all $n,p$.
Then we have~[A1](i)--(ii) regarding $G^n$ while we set
$\mathcal{N}^n=\varnothing$ and $G^n=n\theta^n$. In this case~[A1](iii) corresponds to (a weaker version of) assumption (E) from
\cite{JP2012}, and $(\ref{A4})$ follows from Lemma 14.1.5 of \cite
{JP2012}. Unlike their setting, however, our assumption does not rule
out the dependence between $\varepsilon(n,p)$'s and $X$ (see,
e.g.,~Example~\ref{exhit} in the next section). The importance of such
dependence has recently been emphasized in econometric literature; see,
for example, \citet{RW2011}.

(ii) The assumptions on the quantities $1_{\{\tau^k_p=\tau
^l_p\}}$ are necessary for the treatment of the ($\mathcal
{F}^{(0)}$-con\-ditional) covariance between $\epsilon^k_{\tau^k_p}$ and
$\epsilon^l_{\tau^l_p}$, which is given by $\Upsilon^{kl}_{\tau
^k_p}1_{\{\tau^k_p=\tau^l_p\}}$ (a\vspace*{1pt} similar kind of assumption also
appears in \citet{BM2013} due to the same reason as ours). Therefore,
those assumptions can be dropped when $\Upsilon^{kl}\equiv0$ if $k\neq
l$; this is often assumed in the literature on the covariance
estimation of non-synchronously observed semimartingales with noise.
The quantity $\chi^n$ measures the degree of the non-synchronicity, and
$\chi^n_s$ is a matrix all of whose components are equal to 1 in the
synchronous case while it is an identity matrix in the completely
non-synchronous case. Hence,~[A1](iv) is satisfied in these two
extreme cases.

(iii) The possibility of the set $\mathcal{N}^n$ being
non-empty excludes the following trivial exception of~[A1] with
$\mathcal{N}^n$ being empty: if $T_0=\log n/n$ and $T_p=T_{p-1}+1/n$
for $p\geq1$,~[A1] with $\mathcal{N}^n=\varnothing$ is not
satisfied because $G^n_{T_0}\to\infty$ as $n\to\infty$. This assumption
is also useful to ensure the stability under the localization used in
the proof; see Lemma~\ref{HJYlem41}.

(iv) 
The fact that we consider the conditional expected durations given
$\mathcal{G}^{(0)}_{T_p}$'s instead of $\mathcal{F}^{(0)}_{T_p}$ rules
out some dependence between the sampling schemes and the jumps of the
observed process.
For example, if $\mu$ is a jump measure of a one-dimensional L\'evy
process (i.e.,~$E=\mathbb{R}$) and $T_p$'s are of the form $T_p=\inf\{
t>T_{p-1}\dvt\llvert \int_{T_{p-1}}^t\int_{\llvert z\rrvert \leq
1}z(\mu-\nu)(\mathrm{d}s,\mathrm{d}z)\rrvert >\eta_n\}$ for
$p=1,2,\dots$ and for some appropriate sequence $(\eta_n)_{n\geq1}$ of
positive numbers, then~[A1] obviously fails because $T_p$'s are
$\mathcal{G}^{(0)}_0$-measurable (this type of sampling scheme is well
studied in \citet{RT2011}). On the other hand, it still allows the
presence of the instantaneous causality between the sampling schemes
and the jumps: see Example~\ref{exhitJ}.

 (v) Under~[A1], it holds that
%
\begin{equation}
\label{HJYlem22} \frac{1}{n}N^n_t\to^p
\int_0^t\frac{1}{G_s}\,\mathrm{d}s
\end{equation}
as $n\to\infty$ for every $t\in\mathbb{R}_+$ ({see Section~6.1 of \cite
{Koike2014time} for the proof}). In particular,~[A1] ensures
that the parameter $n$ controls the magnitude of the number of observations.
\end{rmk}

We impose the following structural assumption on the latent process $X$:
%
\begin{longlist}[{[A2]}]
\item[{[A2]}]  The volatility process $\sigma$ is an It\^o
semimartingale of the form
\begin{eqnarray*}
\sigma_t=\sigma_0+\int_0^t
\widetilde{b}_s\,\mathrm{d}s+\int_0^t
\widetilde{\sigma}_s\,\mathrm{d}W_s +(\widetilde{
\delta}1_{\{\llvert \widetilde{\delta}\rrvert \leq1\}})\star(\mu-\nu
)_t +(\widetilde{
\delta}1_{\{\llvert \widetilde{\delta}\rrvert >1\}})\star\mu_t,
\end{eqnarray*}
where $\widetilde{b}_s$ is a locally bounded and $(\mathcal
{F}^{(0)}_t)$-progressively measurable $\mathbb{R}^d\otimes\mathbb
{R}^{d'}$-valued process, $\widetilde{\sigma}_s$ is a cadlag {$(\mathcal
{F}^{(0)}_t)$}-adapted $\mathbb{R}^d\otimes\mathbb{R}^{d'}\otimes\mathbb
{R}^{d'}$-valued process, and $\widetilde{\delta}$ is an {$(\mathcal
{F}^{(0)}_t)$}-predictable $\mathbb{R}^d\otimes\mathbb{R}^{d'}$-valued
function on $\Omega^{(0)}\times\mathbb{R}_+\times E$.

Moreover, for each $j$ there is an {$(\mathcal{F}^{(0)}_t)$}-stopping
time $\rho_j$, a bounded {$(\mathcal{F}^{(0)}_t)$}-progressively
measurable $\mathbb{R}^d$-valued process $b(j)_s$, a deterministic
non-negative function $\gamma_j$ on $E$, and a constant $\Lambda_j$
such that $\rho_j\uparrow\infty$ as $j\to\infty$ and, for each $j$,
%
\begin{enumerate}[(iii)]
\item[(i)] $b(\omega^{(0)})_s=b(j)(\omega^{(0)})_s$ if $s<\rho_j(\omega^{(0)})$,

\item[(ii)] $E [\llVert b(j)_{t_1}-b(j)_{t_2}\rrVert ^2\mid\mathcal
{F}_{t_1\wedge t_2} ]
\leq\Lambda_j E [\llvert t_1-t_2\rrvert ^\varpi\mid\mathcal
{F}_{t_1\wedge t_2} ]$ for any $(\mathcal{F}^{(0)}_t)$-stopping times
$t_1$ and $t_2$ bounded by $j$,

\item[(iii)] {$\int\{\gamma_j(z)^2\wedge1 \}\lambda(\mathrm{d}z)<\infty$} and
$\llVert \delta(\omega^{(0)},t,z)\rrVert \vee\llVert \widetilde
{\delta}(\omega^{(0)},t,z)\rrVert \leq\gamma_j(z)$ for all $(\omega
^{(0)},t,z)$ with $t\leq\rho_j(\omega^{(0)})$,

\item[(iv)] $E [\llVert \delta(t_1\wedge\rho_j,z)-\delta(t_2\wedge\rho
_j,z)\rrVert ^2\mid\mathcal{F}_{t_1\wedge t_2} ]
\leq\Lambda_j\gamma_j(z)^2E [\llvert t_1-t_2\rrvert ^\varpi\mid
\mathcal{F}_{t_1\wedge t_2} ]$ for any $(\mathcal{F}^{(0)}_t)$-stopping
times $t_1$ and $t_2$ bounded by $j$.
\end{enumerate}
\end{longlist}

%
\begin{rmk}
An~\textup{[A2]} type assumption is commonly used in the literature of
power variations (see, e.g.,~\cite{JP2012}), except for assumptions
(ii) and (iv), that is,~continuity assumptions on the drift and the
jump coefficient. Such assumptions are necessary for the treatment of
the irregularity and the non-synchronicity of the observation times as
in \cite{HY2011}.
\end{rmk}

We also impose the following regularity condition on the noise process:
%
\begin{longlist}[{[A3]}]
\item[{[A3]}] There is a constant $\Gamma>4$ and a sequence
$(\rho'_{j})_{j\geq1}$ of {$(\mathcal{F}^{(0)}_t)$}-stopping times
increasing to infinity such that
\[
\sup_{\omega^{(0)}\in\Omega^{(0)},t<\rho'_{j}(\omega^{(0)})}\int
\llVert z\rrVert^\Gamma
Q_t\bigl(\omega^{(0)},\mathrm{d}z\bigr)<\infty.
\]
Moreover, for each $j$ there is a bounded cadlag {$(\mathcal
{F}^{(0)}_t)$}-adapted $\mathbb{R}^d\otimes\mathbb{R}^d$-valued process
$\Upsilon(j)_t$ and a constant $\Lambda_j'$ such that:
%
\begin{enumerate}[(ii)]
\item[(i)] $\Upsilon(j)(\omega^{(0)})_t=\Upsilon(\omega^{(0)})_t$ if $t<\rho
'_j(\omega^{(0)})$,

\item[(ii)] $E [\llVert \Upsilon(j)_{t_1}-\Upsilon(j)_{t_2}\rrVert ^2\mid
\mathcal{F}_{t_1\wedge t_2} ]
\leq\Lambda'_j E [\llvert t_1-t_2\rrvert ^\varpi\mid\mathcal
{F}_{t_1\wedge t_2} ]$ for any $(\mathcal{F}^{(0)}_t)$-stopping times
$t_1$ and $t_2$ bounded by $j$.
\end{enumerate}
\end{longlist}

%
\begin{rmk}
The locally boundedness of the moment process of the noise is used for
verifying a Lyapunov type condition for central limit theorems and
proving the negligibility of the edge effect. The continuity assumption
of the covariance matrix process of the noise is necessary due to the
same reason as for~\textup{[A2]}. If the noise is assumed to be
i.i.d.~and independent of $\mathcal{F}^{(0)}$,~\textup{[A3]} simply
means the $\Gamma$th moment of the noise is finite for some $\Gamma>4$.
\end{rmk}

Finally, we introduce the following technical condition to avoid some
measure-theoretic problems:
%
\begin{enumerate}[{[A4]}]
\item[{[A4]}] (i) A regular conditional probability of $P^{(0)}$
given $\mathcal{H}$ exists for any sub-$\sigma$-field $\mathcal{H}$ of~$\mathcal{F}^{(0)}$. 

\noindent (ii) The process $(Q_t(\cdot,A))_{t\geq0}$ is $(\mathcal
{F}^{(0)}_t)$-progressively measurable for any Borel set $A$ of~$\mathbb{R}^d$.
\end{enumerate}

%
\begin{rmk}
(i)~\textup{[A4]}(i) is satisfied, for example, when $(\Omega
^{(0)},\mathcal{F}^{(0)})$ is a standard measurable space, that is,~it
is Borel isomorphic to some Polish space (see, e.g.,~Theorem I-3.1 of
\cite{IW1989}). In fact, this assumption is not restrictive for applications.

 (ii)~\textup{[A4]}(ii) is satisfied, for example, when
$Q_t\equiv Q$ for some probability measure $Q$ on $\mathbb{R}^d$, that
is,~the noise is modeled by an i.i.d.~sequence. Another example is the
case where $Q_t(\omega^{(0)},\cdot)$ has a density of the form $f(\cdot
,X_t(\omega^{(0)}))$, where $f$ is a measurable function on $\mathbb
{R}^d\times\mathbb{R}^d$ into $[0,1]$ such that $\int_{\mathbb
{R}^d}f(x,\theta)\,\mathrm{d}x=1$ for every $\theta\in\mathbb{R}^d$. 
{Example 16.1.5 of \cite{JP2012} is encompassed with this type of
model.} Thus, this assumption also seems to be unrestrictive for applications.
\end{rmk}

\subsection{Result}

To state the main result, we need the notion of \textit{stable
convergence} as common in this area. For each $n\geq1$, let $X^n$ be a
random variable which is defined on $\mathcal{B}$ and takes values in a
Polish space~$S$. The variables $X^n$ are said to \textit{converge
stably in law} to an $S$-valued random variable $X$ defined on an
extension of $\mathcal{B}^{(0)}$ if $E[Uf(X^n)]\rightarrow\widetilde
{E}[Uf(X)]$ for any $\mathcal{F}^{(0)}$-measurable bounded random
variable $U$ and any bounded continuous function $f$ on $S$, where
$\widetilde{E}$ denotes the expectation with respect to the probability
measure of the extension. We then write $X^n\to^{d_s}X$. Note that we
need a slightly generalized definition of stable convergence described
at the end of Section~2.2.1 of~\cite{JP2012} because $\mathcal{B}$
changes as $n$ varies. The most important property of stable
convergence is the following: if the real-valued variables $V_n$
defined on $\mathcal{B}$ converge in probability to a variable $V$
defined on $\mathcal{B}^{(0)}$, then $X^n\to^{d_s}X$ implies that
$(X^n,V_n)\to^{d_s}(X,V)$ for the product topology on the space $S\times
\mathbb{R}$.

%
\begin{theorem}\label{mainthm}
Suppose that~\textup{[A1]}--\textup{[A4]} are satisfied. Then
\[
n^{1/4} \bigl(\MRC[Y]^n_t-[X,X]_t
\bigr)\to^{d_s}\mathcal{W}_t+\mathcal{Z}_t
\]
as $n\to\infty$ for any $t>0$, where $\mathcal{W}$ and $\mathcal{Z}$
are $\mathbb{R}^d\otimes\mathbb{R}^d$-valued processes defined on an
extension of {$\mathcal{B}^{(0)}$}, which conditionally on {$\mathcal
{F}^{(0)}$} are mutually independent, centered Gaussian with
independent increments, the first one being continuous and the second
one being purely discontinuous, and with (conditional) covariances
%
\begin{eqnarray}\label{avarW}
&& \widetilde{E} \bigl[\mathcal{W}^{kl}_t
\mathcal{W}^{k'l'}_t\mid{\mathcal{F}^{(0)}} \bigr]\nonumber
\\
&&\quad =\frac{2}{\psi_2^{2}}\int_0^t \biggl[
\Phi_{22}\theta\bigl\{\Sigma^{kk'}_s
\Sigma^{ll'}_s+\Sigma^{kl'}_s
\Sigma^{lk'}_s \bigr\}G_s
\nonumber\\[-8pt]\\[-8pt]\nonumber
&&\hspace*{53pt}{} +\frac{\Phi_{11}}{\theta^3}
\bigl\{\Upsilon^{kk'}_s\chi^{kk'}_{s}
\Upsilon^{ll'}_s\chi^{ll'}_{s}+
\Upsilon^{kl'}_s\chi^{kl'}_{s}
\Upsilon^{lk'}_s\chi^{lk'}_{s} \bigr\}
\frac{1}{G_s}
\\
&&\hspace*{53pt}{}+\frac{\Phi_{12}}{\theta} \bigl\{\Sigma^{kk'}_s
\Upsilon^{ll'}_s\chi^{ll'}_{s}+
\Sigma^{lk'}_s\Upsilon^{kl'}_s
\chi^{kl'}_{s}+\Sigma^{ll'}_s
\Upsilon^{kk'}_s\chi^{kk'}_{s}+
\Sigma^{kl'}_s\Upsilon^{lk'}_s
\chi^{lk'}_{s} \bigr\} \biggr]\,\mathrm{d}s\nonumber
\end{eqnarray}
and
%
\begin{equation}
\label{avarZ} \widetilde{E} \bigl[\mathcal{Z}^{kl}_t
\mathcal{Z}^{k'l'}_t\mid{\mathcal{F}^{(0)}} \bigr]=
\frac{1}{\psi_2^{2}}\sum_{s\leq t} \bigl(
\mathfrak{J}^{kk'}_s+\mathfrak{J}^{kl'}_s+
\mathfrak{J}^{lk'}_s+\mathfrak{J}^{ll'}_s
\bigr).
\end{equation}
Here, $\widetilde{E}$ denotes the expectation with respect to the
probability measure of the extension.

When further $X$ is continuous, the processes $n^{1/4} (\MRC[Y]^n-[X,X]
)$ converge stably in law to the process $\mathcal{W}$ for the
Skorokhod topology.
\end{theorem}


%
\begin{rmk}
The above theorem shows that the observation times' effect on the
asymptotic distribution of the MRC estimator is only through the
asymptotic conditional expected duration process $G$ and the limiting
process $\chi$ measuring the degree of the non-synchronicity. As was
indicated in Remark~\ref{rmkA1}(ii), $\chi$ simply reflects the
covariance structure of the noise process, while $G$ naturally affects
the asymptotic distribution of the estimator because it links with the
(spot) sampling frequency, as seen from $(\ref{HJYlem22})$.
Consequently, the irregularity and the endogeneity of the observation
times have no impact on the asymptotic distribution of the estimator.
\end{rmk}


%
\begin{rmk}\label{rmkavg}
In the proof of the theorem, it plays a key role to replace the
duration $(T_{p+1}-T_{p})$ with its conditional expectation
$G^n_{T_p}$. Such replacement is possible because our estimator
contains a local averaging procedure $(\ref{defPA})$. More formally,
this procedure makes it possible to apply a standard martingale
argument described in {Lemma 2.3 of \cite{Fu2010b}} to the durations.
The benefits of this fact appear in the treatments of the irregularity
and the endogeneity of the observation times in {Lemmas~\ref
{appstoptime} and~\ref{JPlem1637}}.
Also, this is why the higher (conditional) moments of the durations do
not affect the asymptotic distribution of the estimator.
\end{rmk}


%
\begin{rmk}[(Covariance structure of $\mathcal{W}_t$)]\label{rmkW}
It is convenient to observe that the covariance structure of $\mathcal
{W}_t$ is analogous to the asymptotic covariance of the realized
covariance in a standard setting. For this purpose, in the following we
use some concepts from matrix algebra found in, for example,~\citet{HJ1991}.
For each $s\in\mathbb{R}_+$, we denote by $\widetilde{\Upsilon}_s$ the
Hadamard product of $\Upsilon_s$ and $\chi_s$, that is, $\widetilde
{\Upsilon}^{kl}_s=\Upsilon^{kl}_s\chi^{kl}_s$ for $k,l=1,\dots,d$, and
set $\overline{\Sigma}_s(y)=\frac{\sqrt{2}}{\psi_2} \{\phi
_{g,g}(y)\theta^{\sfrac{1}{2}}\Sigma_s\sqrt{G_s}+\phi_{g',g'}(y)\theta
^{-\sfrac{3}{2}}\widetilde{\Upsilon}_s/\sqrt{G_s} \}$. Since both
$\Upsilon_s$ and $\chi_s$ is positive semidefinite, so is $\widetilde
{\Upsilon}_s$ due to the Schur product theorem, that is,~Theorem 5.2.1
of \cite{HJ1991} (note that the positive semidefiniteness of $\chi_s$
can be checked directly using the fact that $\chi^{kk}_s=1$ and $0\leq
\chi^{kl}_s\leq1$ for any $k,l$). Therefore, $\overline{\Sigma}_s(y)$~is positive semidefinite as well because both $\phi_{g,g}$ and $\phi
_{g',g'}$ are non-negative (see Remark~\ref{rmkphi}). Then the
left-hand side of \eqref{avarW} can be rewritten as
\[
\int_0^t \biggl[\int_0^1
\bigl\{\overline{\Sigma}_s(y)^{kk'}\overline{
\Sigma}_s(y)^{ll'}+\overline{\Sigma}_s(y)^{kl'}
\overline{\Sigma}_s(y)^{lk'} \bigr\}\,\mathrm{d}y \biggr]
\,\mathrm{d}s.
\]
The integrand of the above expression is nothing but the $\mathcal
{F}^{(0)}$-conditional covariance between the $(k,l)$th and $(k',l')$th
entries of the variable $\overline{\Sigma}_s(y)^{1/2}\zeta(\overline
{\Sigma}_s(y)^{1/2}\zeta)^*$, where $\zeta$ is a $d$-dimensional
standard normal variable independent of $\mathcal{F}^{(0)}$. In other
words, $\vectorize(\mathcal{W}_t)$ is centered Gaussian with covariance
matrix $\int_0^1\mathfrak{S}_s\,\mathrm{d}s$, where $\mathfrak{S}_s=\int
_0^1(\overline{\Sigma}_s(y)\otimes\overline{\Sigma}_s(y))\Cov[\vectorize
(\zeta\zeta^*)]\,\mathrm{d}y$ and, $\vectorize$ and $\otimes$ denote the
vec-operator and the Kronecker product of matrices, respectively
(cf.~Section~2.2 of~\cite{BHMR2014}). In particular, the process
$\mathfrak{S}_s$ is cadlag, $(\mathcal{F}^{(0)}_t)$-adapted and takes
values in the set of $d\times d$ positive semidefinite matrices, hence
we can construct the process $\mathcal{W}$ stated as in the theorem by
Proposition 4.1.2 of \cite{JP2012}. More precisely, $\mathcal{W}$ can
be realized as $\vectorize(\mathcal{W}_t)=\int_0^t\mathfrak
{S}_s^{1/2}\,\mathrm{d}W'_s$, where $W'$ is a $d^2$-dimensional standard
Brownian motion defined on an extension of $\mathcal{B}^{(0)}$ and
independent of $\mathcal{F}^{(0)}$.

Note that the Fisher information matrix for covariance matrix
estimation of a multivariate diffusion process from non-synchronous and
noisy observations is \textit{not} analogous to that for a pure
diffusion setting; see Section~2.2 of \cite{BHMR2014} for details.
\end{rmk}


%
\begin{rmk}[(Covariance structure of $\mathcal{Z}_t$)]\label{rmkZ}
$\mathcal{Z}_t$ apparently has an analogous covariance structure to the
asymptotic covariance of the realized covariance due to jumps in the
regular sampling case (cf.~equation (5.4.4) of \cite{JP2012}), and it
can be realized as follows. Set $A_m=\{z\dvt\gamma(z)>1/m\}$ for each
$m\in\mathbb{N}$, and denote by $(S(m,j))_{j\geq1}$ the successive jump
times of the Poisson process $1_{A_m\setminus A_{m-1}}\star\mu$. Let
$(S_r)_{r\geq1}$ be a reordering of the double sequence $(S(m,j))$.
Suppose that sequences $(\Psi_{r-})_{r\geq1}$ and $(\Psi_{r+})_{r\geq
1}$ of i.i.d.~standard $d'$-dimensional normal variables and sequences
$(\Psi'_{r-})_{r\geq1}$ and $(\Psi'_{r+})_{r\geq1}$ of i.i.d.~standard
$d$-dimensional normal variables are defined on an extension of
$\mathcal{B}^{(0)}$ and that all of them are mutually independent and
independent of $\mathcal{F}^{(0)}$. Now, the variable $\widetilde
{\Upsilon}_s$ defined in Remark~\ref{rmkW} is positive semidefinite, it
admits the (positive semi-definite) square root $\widetilde{\upsilon
}_s:=\Upsilon_s^{1/2}$. Since the process $\widetilde{\Upsilon}$ is
cadlag and $(\mathcal{F}^{(0)}_t)$-adapted, so is $\widetilde{\upsilon
}_s$. Then $\mathcal{Z}$ is realized as $\mathcal{Z}_t=\sum_{r\dvt
S_r\leq t}(\mathfrak{Z}_r+\mathfrak{Z}_r^*)$, where
\begin{eqnarray*}
\mathfrak{Z}_r=\frac{1}{\psi_2}\Delta X_{S_r} \biggl\{
\sqrt{\Phi_{22}\theta} (\sigma_{S_r-}\sqrt{G_{S_r-}}
\Psi_{r-}+\sigma_{S_r}\sqrt{G_{S_r}}
\Psi_{r+} )+\sqrt{\frac{\Phi_{12}}{\theta}} \bigl(\widetilde{
\upsilon}_{S_r-}\Psi'_{r-}+\widetilde{
\upsilon}_{S_r}\Psi'_{r+} \bigr) \biggr\}^*.
\end{eqnarray*}
This is indeed the desired one; see Proposition 4.1.4 of \cite{JP2012}.
\end{rmk}


%
\begin{rmk}[(Comparison with a pure semimartingale setting)]
It would be interesting to observe how our result is different from
\citet{BibV2013}'s one in a pure semimartingale setting. For
simplicity, we focus on the univariate case, that is, we assume that
$d=d'=1$, and assume that $T_p=t^1_p$ for every $p$ for notational
simplicity. Now let us recall their result briefly. Suppose that $b$,
$\sigma$ and $\delta$ are continuous. Suppose also that the sequence
$(T_p)$ is independent of $b,\sigma,\delta,W,\mu$ and satisfies~$(\ref
{A4})$. Then, according to Theorem 2 of~\cite{BibV2013}, for any
$t>0$ we have the following convergence:
%
\begin{equation}
\label{BVresult} \sqrt{n}\bigl([X,X]^n_t-[X,X]_t
\bigr)\to^{d_s}\sqrt{2}\int_0^t
\sigma_s^2\sqrt{H'(s)}
\,\mathrm{d}W'_s+2\sum_{r\dvt S_r\leq t}
\Delta X_{S_r}\sigma_{S_r}\sqrt{\eta(S_r)}
\Psi_r,
\end{equation}
where $W'$ is a standard Brownian motion, $H$ is a (possibly random)
$C^1$ function such that $n\sum_{p\dvt T_p\leq t}(T_p-T_{p-1})^2\to
^pH(t)$ for every $t\in\mathbb{R}_+$ (the existence is assumed),
$(S_r)_{r\geq1}$ is a sequence of stopping times exhausting the jumps
of $X$, $(\Psi_r)_{r\geq1}$ is a sequence of i.i.d.~standard normal
variables, and $(\eta(t))_{t\in\mathbb{R}_+}$ is a family of
independent random variables with uniformly bounded first moments, and
such that the processes $(n(T_+(t)-T_-(t)))_{t\in\mathbb{R}_+}$
converge finite-dimensionally in law to $(\eta(t))_{t\in\mathbb{R}_+}$
(the existence is assumed, and this condition can be weakened; see
Assumption 2 of \cite{BibV2013} for details). Here, $T_+(t)=\min\{
T_p\dvt T_p\geq t\}$ and $T_-(t)=\max\{T_p\dvt T_p\leq t\}$ for any
$t\in\mathbb{R}_+$ and $W'$, $(\Psi_r)$ and $(\eta(t))$ are defined on
an extension of $\mathcal{B}$ and mutually independent as well as
independent of $\mathcal{F}$. On the other hand, provided that $\Upsilon
\equiv0$ (so the noise is absent), the corresponding result to our
estimator can be written as follows:
%
\begin{eqnarray}\label{univresult}
&& n^{1/4}\bigl(\MRC[Y]^n_t-[X,X]_t
\bigr)
\nonumber\\[-8pt]\\[-8pt]\nonumber
&&\quad \to^{d_s}\frac{\sqrt{2\Phi_{22}\theta}}{\psi_2} \biggl(\sqrt
{2}\int
_0^t\sigma_s^2
\sqrt{G_s}\,\mathrm{d}W'_s+2\sum
_{r\dvt S_r\leq t}\Delta X_{S_r}\sigma_{S_r}
\sqrt{G_{S_r}}\Psi_r \biggr),
\end{eqnarray}
where we also assume that $G$ is continuous for simplicity. Compared
with the above equation with $(\ref{BVresult})$, the quantities $H'$
and $\eta$ coming from the irregularity of the observation times in the
left-hand side of $(\ref{BVresult})$ are replaced with $G$ in $(\ref
{univresult})$. Since the quantity $H$ contains the information of the
second moments of the durations and $\eta$ contains that of all the
moments of the durations around the jump times, the distributional
future of the durations strongly affects the asymptotic distribution in
$(\ref{BVresult})$. In contrast, the first moments of the durations
only affect the asymptotic distribution in $(\ref{univresult})$.
\end{rmk}


%
\begin{rmk}[(Comparison with the continuous case)]
The result of the theorem is not new if $X$ is continuous. In fact, in
the case that $X$ is continuous, a central limit theorem for the MRC
estimator can be derived with a somewhat weaker assumption on the
limiting process $G$; see Theorem 3.1 and assumption~\textup{[A4]} of
Koike \cite{Koike2014time} for details. In the discontinuous case, we need
some regularity of the path of the left limit process $G_-$ to verify
the approximation given in Proposition~\ref{JPlem1218}, so the
structural assumption~\textup{[A4]}(iii-c) is necessary.

It is worth mentioning that the structural assumption on $G$ is
necessary to deal with the irregularity of observation times in the
discontinuous case. In contrast, such a condition is only required to
handle the time endogeneity in the continuous case. In fact, if the
observation times have a kind of pre-determination property (the
so-called \textit{strong predictability}), convergence in probability
of $G^n$ to $G$ for the Skorokhod topology is sufficient to derive a
central limit theorem; see \citet{Koike2014phy} for details.
\end{rmk}


%
\begin{rmk}[(Feasible limit theorem)]
In order to apply Theorem~\ref{mainthm} to real statistical problems
such as the construction of confidence intervals, we need an estimator
for the asymptotic covariance matrix given by $(\ref{avarW})$ and $(\ref
{avarZ})$. This will be achieved by combining the technique used in the
non-synchronously observed diffusion setting (e.g.,~a kernel approach
of \cite{HY2011} or a histogram-type method of \cite{Bibinger2012})
with the one used in the jump diffusion setting (e.g.,~a thresholding
and locally averaging method of
A{\"{\i}}t-Sahalia, Jacod and Li
\cite{AJL2011}). Or we can presumably
use an estimator of \citet{AX2014} for the equidistant sampling setting
without modification because the distribution of the variable
$n^{1/4}\overline{Y}_i$ is, roughly speaking, approximated by the
$d$-dimensional normal variable with mean 0 and covariance matrix
$\theta\psi_2\Sigma_{T_i}G_{T_i}+\frac{\psi_1}{\theta}\widetilde
{\Upsilon}_{T_i}$ in the absence of jumps conditionally on $\mathcal
{F}^{(0)}_{T_i}$, where $\widetilde{\Upsilon}$ is the same one as in
Remark~\ref{rmkW} (this is theoretically manifested by Lemma~\ref
{JPlem1637} in a sense).
\end{rmk}


\section{Examples of the observation times}\label{examples}

In this section, we give some illustrative examples of the observation
times that satisfy the condition~\textup{[A1]}.
We shall start to discuss univariate examples (i.e.,~we assume that
$d=1$), which are not encompassed with the restricted discretization schemes.

%
\begin{example}\label{exhit}
As an illustrative example of endogenous observation times, we consider
a simple model generated by hitting times of the underlying Brownian
motion $W$. This type of model is commonly used in the literature; see
\cite{Fu2010b,LMRZZ2014,RHW2014} among others. Here, we especially
focus on a simpler version of the specification from \cite{RHW2014}.
Specifically, $t^1_i$'s are defined as follows:\vspace*{-2pt}
\[
t^1_0=0,\qquad t^1_{i+1}=
\inf\bigl\{t>t^1_{i}\dvt W_t-W_{t^1_i}+
\sqrt{n}\mathfrak{a}_{t^1_i}\bigl(t-{t^1_i}
\bigr)=\mathfrak{b}_{t^1_i}/\sqrt{n} \bigr\},
\]
where $\mathfrak{a}$ and $\mathfrak{b}$ are cadlag $(\mathcal
{F}^{(0)}_t)$-adapted processes such that $\mathfrak{a}_t\mathfrak
{b}_t>0$ and $\mathfrak{a}_{t-}\mathfrak{b}_{t-}>0$ for every~$t$. In\vspace*{1pt}
this case,~\textup{[A1]} is satisfied with setting $T_p=\tau^1_p=t^1_p$
for every $p$, as long as $G:=\mathfrak{b}/\mathfrak{a}$ satisfies~[A1](iii-c). In fact, noting that, conditionally on $\mathcal
{F}^{(0)}_{T_p}$, $n(T_{p+1}-T_p)$ follows the inverse Gaussian
distribution with mean $G_{T_p}$ and variance $G_{T_p}^2/\mathfrak
{a}_{T_p}$, $(\ref{A4})$ holds true for any $t>0$ and $\xi\in(0,1)$.
Moreover, we have $E[n(T_{p+1}-T_p)\mid{\mathcal
{G}^{(0)}_{T_p}}]=G_{T_p}$ for every $p$ because $W$ is independent of
{$\mu$}. Hence,~\textup{[A1]}(i)-(iii) are satisfied with $\mathcal
{N}^n=\varnothing$. Finally,~\textup{[A1]}(iv) is automatically satisfied.
\end{example}


%
\begin{example}\label{exhitJ}
We can also accommodate observation times generated by hitting times of
a Brownian motion plus finitely many jumps to our situation. For
example, let us consider the observation times defined as follows:\vspace*{-2pt}
\[
t^1_0=0,\qquad t^1_{i+1}=\inf
\bigl\{t>t^1_{i}\dvt W_t-W_{t^1_i}+
\sqrt{n}\mathfrak{a}_{t^1_i}\bigl(t-{t^1_i}
\bigr)+\delta'\star\mu_t-\delta'\star
\mu_{t^1_i}=\mathfrak{b}_{t^1_i}/\sqrt{n} \bigr\},\vspace*{-1pt}
\]
where $\mathfrak{a}$ and $\mathfrak{b}$ are the same one as in Example
\ref{exhit} and $\delta'$ is an $(\mathcal{F}^{(0)}_t)$-optional
real-valued function on $\Omega^{(0)}\times\mathbb{R}_+\times E$ such
that $1_{\{\delta'\neq0\}}\star\mu_t<\infty$ for all $t$. Therefore,
the process $\delta'\star\mu$ has finitely many jumps. 
Then, it can easily been seen that~\textup{[A1]} is satisfied in this
case under the same situation as that of Example~\ref{exhit}, except
for setting $\mathcal{N}^n=\{p\in\mathbb{Z}_+\dvt\delta'\star\mu
_{T_{p+1}}-\delta'\star\mu_{T_{p}}>0\}$.
\end{example}


%
\begin{example}
Let us consider the observation times discussed in Example 3 of \citet
{BibV2013}. Namely, $t^1_i=i/n$ if $i$ is even and $t^1_i=(i+\alpha)/n$
if $i$ is odd, where $\alpha\in(0,1)$ is a constant. \cite{BibV2013}
showed that this observation times produce an additional randomness in
the asymptotic distribution of the realized covariance estimator even
though they are deterministic. In fact, in this case the variable $\eta
(t)$ in $(\ref{BVresult})$ takes the values $(1+\alpha)$ and $(1-\alpha
)$ with probabilities $(1+\alpha)/2$ and $(1-\alpha)/2$, respectively.
On the other hand, setting $T_p=(p+1)/n$ and $\tau^1_p=t^1_p$,~[A1] is satisfied. Hence, in our case this example has the same
impact as that of the regular observation times on the asymptotic distribution.
\end{example}

Next, we turn to the multivariate and non-synchronous examples. As the
data synchronization method, we focus on the refresh sampling method.


%
\begin{example}
We shall discuss the Poisson sampling, which is one of the most popular
models in this area; see, for example,~\cite
{Bibinger2012,BibV2013,HY2005,Zhang2011}. Let\vspace*{1pt} $(t^k_i)$ be a sequence
of Poisson arrival times with the intensity $np_k$ for each $k$ and
suppose that $(t^1_i),\dots,(t^d_i)$ are mutually independent and
independent of $X$ and $\epsilon$. Then,~\textup{[A1]} is satisfied with
\[
G_s\equiv\sum_{k=1}^d\sum
_{1\leq l_1<\cdots<l_k\leq d}\frac{(-1)^{k-1}}{p_{l_1}+\cdots+
p_{l_k}},\qquad\chi^{kl}_s
\equiv
\cases{
1, &\quad if $k=l$,
\vspace*{3pt}\cr
0, &\quad otherwise.}
\]
%
\end{example}


%
\begin{example}\label{LoMac}
Here, we give an example of observation times which are possibly
endogenous and satisfy~\textup{[A1]} with the explicit $G$ and $\chi$.
More precisely, we give a continuous time analog of the Lo--MacKinlay
model of \cite{LM1990}.

Let $(\tau_i)_{i=0}^\infty$ be a sampling scheme and suppose that $\sup
_{i\geq0}(\tau_i\wedge t-\tau_{i-1}\wedge t)=\mathrm{o}_p(n^{-\xi})$ as $n\to
\infty$ for any $t>0$ and $\xi\in(0,1)$. For each $k=1,\dots,d$, let
$(M_k(n,i))_{i=0}^\infty$ be a sequence of \mbox{$\mathbb{Z}_+$-}valued
variables defined on $\mathcal{B}^{(0)}$ and independent of $X$ and
$(\tau_i)$ such that $\mathcal{M}_k:=(M_k(n,i+1)-M_k(n,i))_{i=0}^\infty
$ is independent and geometrically distributed with the common success
probability $p_k\in(0,1)$. Moreover, suppose that, for each $i$,
$M_k(n,i)$ is an $(\mathcal{F}^{(0)}_{\tau_j})_{j=0}^\infty$-stopping
time so that $t^k_i:=\tau_{M_k(n,i)}$ is an $(\mathcal
{F}^{(0)}_t)$-stopping time and that $M_k(n,i+1)-M_k(n,i)$ is\vspace*{1pt}
independent of $\mathcal{F}^{(0)}_{t^k_i}$. Finally, assume that
$\mathcal{M}_1,\dots,\mathcal{M}_d$ are mutually independent and that
\textup{[A1]}(i)--(iv) are satisfied with replacing $(T_i)$ by $(\tau
_i)$. Then it can easily be shown that~\textup{[A1]} holds true with
\begin{eqnarray*}
&& G_s\equiv\sum_{k=1}^d\sum
_{1\leq l_1<\cdots<l_k\leq d}\frac{(-1)^{k-1}G^0_s}{1-(1-p_{l_1})\cdots
(1-p_{l_k})}
\end{eqnarray*}
and
\[
\chi^{kl}_s\equiv\cases{ 1, &\quad if $k=l$,
\vspace*{3pt}\cr
p_kp_l/(p_k+p_l-p_kp_l),
&\quad otherwise.}
\]
%
Here, $G^0$ denotes the asymptotic conditional expected duration
process corresponding to $(\tau_i)$. By taking an endogenous sampling
scheme as the underlying sampling scheme $(\tau_i)$, we can obtain
endogenous observation times.
\end{example}


\section{Simulation study}\label{simulation}

In this section, we assess the finite sample accuracy of the central
limit theory developed in this paper and confirm our theoretical
findings via Monte Carlo experiments.

We simulate over the unit interval $[0,1]$, and basically follow the
design of \citet{AX2014}. To simulate the latent semimartingale $X$,
the following bivariate Heston model with jumps is considered:
\begin{eqnarray}
\mathrm{d}X^k_t=\sigma_{k,t}
\,\mathrm{d}W^k_t+\mathrm{d}Z^k_t,
\qquad\mathrm{d}\sigma^2_{k,t}=\kappa_k\bigl(
\bar{\sigma}^2_k-\sigma^2_{k,t}
\bigr)\,\mathrm{d}t+s_k\sigma_{k,t}\,\mathrm{d}B^k_t+
\mathrm{d}J^k_t-\lambda^V_k
\tau^V_k\,\mathrm{d}t,\nonumber
\\
\eqntext{k=1,2.}
\end{eqnarray}
Here, $W^1,W^2,B^1,B^2$ are correlated standard Brownian motions such that
\begin{eqnarray*}
\mathrm{d}\bigl[W^1,W^2\bigr]_t&=&
\rho_B\,\mathrm{d}t,\qquad\mathrm{d}\bigl[W^k,B^k
\bigr]_t=\rho_k\,\mathrm{d}t,
\\
\mathrm{d}
\bigl[W^1,B^2\bigr]_t&=&\mathrm{d}
\bigl[W^2,B^1\bigr]_t=\mathrm{d}
\bigl[B^1,B^2\bigr]_t=0.
\end{eqnarray*}
$J^k$ is a compound Poisson process with jump size uniformly
distributed on $[0,2\tau^V_k]$ and jump intensity $\lambda^V_k$. $J^1$
and $J^2$ are assumed to be mutually independent. $Z^k$ is a pure jump
L\'evy process specified as follows. First, $Z^2$ is linearly
correlated with $Z^1$ as $Z^2=\rho_JZ^1+\sqrt{1-\rho_J^2}Z^0$, where
$Z^0$ is another L\'evy process independent of $Z^1$. For each $m=0,1$,
$Z^m$ is a CGMY process with L\'evy density given by
\begin{eqnarray*}
f_m(x)=c_m\frac{\ee^{-\gamma_{m-}\llvert x\rrvert }}{\llvert x\rrvert
^{1+\beta_m}}1_{\{x<0\}}+c_m
\frac{\ee^{-\gamma_{m+}x}}{x^{1+\beta_m}}1_{\{x>0\}}.
\end{eqnarray*}

%
\begin{table}[t]
\tabcolsep=0pt
\caption{The parameters of the stochastic volatility processes}\label{svparam}
\begin{tabular*}{\tablewidth}{@{\extracolsep{\fill}}@{}llllllll@{}}
\hline
$k$&$\kappa_k$&$s_k$&$\bar{\sigma}_k$&$\rho_k$&$\lambda^V_k$&$\tau
^V_k$&$\rho_B$\\
\hline
1&5&0.3&0.25&$-$0.6&\phantom{0}5&0.05&0.5\\
2&4&0.4&0.3&$-$0.75&10&0.01&---\\ \hline
\end{tabular*}
\end{table}
%

The parameter values of the stochastic volatility processes used in the
simulation are reported in Table~\ref{svparam}.
The initial value for the volatility processes $\sigma^2_{k,t}$ is set
at $\bar{\sigma}^2_k$ for each $k=1,2$, which ensures that $E[\sigma
^2_{k,t}]=\bar{\sigma}^2_k$ for all $t\in[0,1]$.
The specification of the parameters in the CGMY processes is as
follows. We set $\gamma_{m+}=3,\gamma_{m-}=5,\beta_m=0.5$ for every
$m=0,1$. $c_1$ is selected such that the quadratic variation
contributed by jumps in $X^1$ amounts to 15\% in expectation, that
is,~$E([Z^1,Z^1]_1)/E([X^1,X^1]_1)=0.15$. Then $c_0$ is selected such
that $E([Z^2,Z^2]_1)/E([X^2,X^2]_1)=0.15$. Finally, the correlation
parameter $\rho_J$ between the jump processes are set at 0.2. Note that
$Z^1$ and $Z^2$ can be exactly simulated because we only consider the
situation where they are of finite variation; see, for example,~\cite
{KM2011} for details.

To generate observation times, we consider Lo--MacKinlay type sampling
schemes illustrated in Example~\ref{LoMac}. Two kinds of sequence $(\tau
_i)_{i=0}^\infty$ of latent observation times are considered: one is
the equidistant sampling scheme $\tau_i=i/n$ and the other is the
endogenous sampling scheme defined by
%
\begin{equation}
\label{simulatehitting} 
\hspace*{-8pt}\tau_0=0,\qquad\tau_{i+1}=\inf\bigl\{t>
\tau_i\dvt W^1_t-W^1_{\tau_i}-2
\sqrt{n}(t-\tau_i)=-2/\sqrt{n}\bigr\},\qquad i=0,1,\dots,\quad
\end{equation}
where we set $n={}$23\,400. Note that in the latter case the sequence
$(\tau_{i+1}-\tau_i)_{i=0}^\infty$ is independent and identically
distributed with the inverse Gaussian distribution with mean $1/n$ and
variance $4/n^2$, thus we can exactly simulate $\tau_i$'s (and
construct the exactly discretized path $\{W_{\tau_i}\}$ from $\{\tau_i\}
$). Furthermore, in both cases the corresponding conditional expected
duration processes $G^0$ are identical with 1.
The parameters $p_1$ and $p_2$ from Example~\ref{LoMac}, which denote
the probabilities of observations occurring, are assumed to be
identical each other and varied thorough $1/3,1/5,1/10$ and $1/30$.



In constructing noisy prices $Y$, we first generate a discretized path
$X_{\tau_0},X_{\tau_1},\dots$ of $X$ using a standard Euler scheme.
After that, we add simulated microstructure noise $Y_{\tau_i}=X_{\tau
_i}+\epsilon_{\tau_i}$ by generating centered Gaussian i.i.d.~variables
$\epsilon^k_{\tau_0},\epsilon^k_{\tau_1},\dots$ with standard deviation
0.005. $\epsilon^1$ and $\epsilon^2$ are assumed to be mutually independent.
Simulation results are based on 10\,000 Monte Carlo iterations for each scenario.

%
\begin{table}[t]
\tabcolsep=0pt
\caption{Simulation results of the standardized estimates}\label{tableStd}
\begin{tabular*}{\tablewidth}{@{\extracolsep{\fill}}@{}lllllllll@{}}
\hline
& \multicolumn{4}{c}{$\tau_i=i/n$} & \multicolumn{4}{c}{$\tau_i$'s are defined by $(\protect\ref{simulatehitting})$} \\[-6pt]
& \multicolumn{4}{c}{\hrulefill} & \multicolumn{4}{c}{\hrulefill} \\
&& & Coverage & Coverage & & & Coverage & Coverage \\
& Mean & SD & (95\%) & (99\%) & Mean & SD & (95\%) & (99\%) \\
\hline
$\theta=1/3$\\
$p_1=p_2=1/3$ & $-$0.00 & 1.01 & 0.949 & 0.987 & $-$0.00 & 1.02 &0.948 & 0.987 \\
$p_1=p_2=1/5$ & $-$0.01 & 1.02 & 0.946 & 0.987 & $-$0.01 & 1.04 &0.943 & 0.986 \\
$p_1=p_2=1/10$ & $-$0.01 & 1.05 & 0.939 & 0.983 & $-$0.01 & 1.06 &0.937 & 0.984 \\
$p_1=p_2=1/30$ & $-$0.03 & 1.09 & 0.928 & 0.979 & $-$0.03 & 1.10 &0.929 & 0.980
\\[3pt]
$\theta=1$\\
$p_1=p_2=1/3$ & $-$0.01 & 1.01 & 0.948 & 0.987 & $-$0.01 & 1.01 &0.949 & 0.989 \\
$p_1=p_2=1/5$ & $-$0.01 & 1.01 & 0.948 & 0.987 & $-$0.01 & 1.01 &0.951 & 0.987 \\
$p_1=p_2=1/10$ & $-$0.02 & 1.02 & 0.947 & 0.986 & $-$0.02 & 1.02 &0.948 & 0.987 \\
$p_1=p_2=1/30$ & $-$0.03 & 1.03 & 0.946 & 0.985 & $-$0.03 & 1.03 &0.943 & 0.985 \\ \hline
\end{tabular*}
\tabnotetext{dd}{\textit{Note}. We report the sample mean, standard deviation (SD) as
well as the 95\% and 99\% coverages of the standardized statistics (\protect\ref{infeasible}) included in the simulation study.}
\end{table}

Following \cite{CKP2010}, the MRC estimator is implemented using the
weight function $g(x)=x\wedge(1-x)$ and the refresh time sampling
method (the finite sample corrections explained in \cite{CKP2010} are
also included).
We consider the window size $k_n$ of the form $k_n=\lceil\theta\sqrt
{N^n_1}\rceil$, and $\theta$ is selected among $1/3$ and $1$. The
former value of $\theta$ corresponds to the one used in \citet
{JLMPV2009}, while the latter one does to the one used in \citet{CKP2010}.
We assess the accuracy of the standard normal approximation of the
infeasible standardized statistic\looseness=-1
%
\begin{equation}
\label{infeasible} n^{1/4}\frac{\MRC[Y]^{n,12}_1-[X^1,X^2]_1}{\sqrt
{\mathbf{AVAR}}},
\end{equation}
where $\mathbf{AVAR}$ is the theoretical asymptotic variance given in
Theorem~\ref{mainthm}. Table~\ref{tableStd} reports the sample mean and
standard deviation as well as 95\% and 99\% coverages of $(\ref
{infeasible})$. As the table reveals, the central limit theorem for
$(\ref{infeasible})$ fairly works.
As was expected from the theory developed in the above, we find no
significant difference of the results between the exogenous and the
endogenous sampling cases.
At relatively low frequencies like $p_1=p_2=1/10$ or $1/30$, the
results for $\theta=1$ show the better performance than those for
$\theta=1/3$. This would be because $k_n$ is not sufficiently large in
such a situation, in order to work the averaging effect of the
pre-averaging procedure explained in Remark~\ref{rmkavg}.



\section{Proof of Theorem \texorpdfstring{\protect\ref{mainthm}}{3.1}}\label{proofs}

\subsection{Preliminaries}

\subsubsection{Localization}

Before starting the proof, we strengthen our assumptions~[A1]--\textup{[A3]} by localization procedures. First, a standard
localization procedure, described in detail in Lemma 4.4.9 of \cite
{JP2012}, for instance, allows us to replace the conditions~[A2] and~\textup{[A3]} by the following strengthened versions, respectively:

%
\begin{longlist}[{[SA2]}]
\item[{[SA2]}] We have~\textup{[A2]}, and the processes $X_t$, $b_t$, $\sigma_t$,
$\widetilde{b}_t$ and $\widetilde{\sigma}_t$ are bounded. Also, $b_t$
is $(\mathcal{H}^\wedge_t)$-progressively measurable and $\sigma_t$ is
$(\mathcal{H}^\wedge_t)$-adapted. Moreover, there are a constant
$\Lambda$ and a non-negative bounded function $\gamma$ on $E$ such that
$\int\gamma(z)^2\lambda(\mathrm{d}z)<\infty$ and $\llVert \delta(\omega
^{(0)},t,z)\rrVert \vee\llVert \widetilde{\delta}(\omega
^{(0)},t,z)\rrVert \leq\gamma(z)$ and
\begin{eqnarray*}
E \bigl[\llVert b_{t_1}-b_{t_2}\rrVert^2\mid
\mathcal{F}_{t_1\wedge t_2} \bigr] &\leq&\Lambda E \bigl[\llvert
t_1-t_2\rrvert^\varpi\mid
\mathcal{F}_{t_1\wedge t_2} \bigr],
\\
E \bigl[\bigl\llVert\delta(t_1,z)-\delta(t_2,z)\bigr
\rrVert^2\mid\mathcal{F}_{t_1\wedge t_2} \bigr] &\leq&\Lambda
\gamma(z)^2E \bigl[\llvert t_1-t_2\rrvert
^\varpi\mid\mathcal{F}_{t_1\wedge t_2} \bigr]
\end{eqnarray*}
for any bounded $(\mathcal{F}^{(0)}_t)$-stopping times $t_1$ and $t_2$.

\item[{[SA3]}] There are a constant $\Gamma>4$ and a constant $\Lambda'$ such that the
process $\int\llVert z\rrVert ^\Gamma Q_t(\mathrm{d}z)$ is bounded and
\[
E \bigl[\llVert \Upsilon_{t_1}-\Upsilon_{t_2}\rrVert ^2\mid\mathcal
{F}_{t_1\wedge t_2} \bigr]
\leq\Lambda' E \bigl[\llvert t_1-t_2\rrvert ^\varpi\mid\mathcal
{F}_{t_1\wedge t_2} \bigr]
\]
for any bounded $(\mathcal{F}^{(0)}_t)$-stopping times $t_1$ and $t_2$.
Moreover, $\Upsilon_t$ is cadlag and $(\mathcal{H}^\wedge_t)$-adapted.
\end{longlist}

Next, we introduce a strengthened version of~\textup{[A1]}. In the
following, we fix a constant $\xi\in(0,1)$ such that
%
\begin{equation}
\label{estxi} \xi>\tfrac{7}{8}\vee\tfrac{1}{2} \bigl(\kappa+
\tfrac{3}{2} \bigr)\vee(1-\varpi),
\end{equation}
and we set $\bar{r}_n=n^{-\xi}$.
%
\begin{longlist}[{[SA1]}]

\item[{[SA1]}] We have~\textup{[A1]}, and for every $n$ it holds that
%
\begin{equation}
\label{SA4} \sup_{p\geq0}(T_p-T_{p-1})
\leq\bar{r}_n.
\end{equation}
\end{longlist}

The following lemma allows us to replace~\textup{[A1]} by~\textup{[SA1]}
via another localization argument. The proof is similar to that of
Lemma 6.3 from \cite{Koike2014time}, so we omit it.

\begin{lem}\label{HJYlem41}
Assume~\textup{[A1]}. One can find sampling schemes $(\widetilde
{T}_{p})$ and $(\widetilde{\tau}^k_{p})$ $(k=1,\dots,d)$ satisfying the
following conditions:
%
\begin{longlist}[(ii)]
\item[(i)] $(\widetilde{T}_{p})$ and $(\widetilde{\tau}^k_{p})$ satisfy~\textup{[SA1]} with the same limiting processes $G$ and $\chi$ as those of
the original sampling schemes.

\item[(ii)] For any $t>0$ there is a subset $\Omega^{(0)}_{n,t}$ of $\Omega
^{(0)}$ such that $\lim_nP^{(0)}(\Omega^{(0)}_{n,t})=1$. Moreover, on
$\Omega^{(0)}_{n,t}$ we have $T_p\wedge t=\widetilde{T}_p\wedge t$ and
$\tau^k_p\wedge t=\widetilde{\tau}^k_{p}\wedge t$ for all $k,p$.
\end{longlist}
\end{lem}

\subsubsection{Outline of the proof}


Here, we give a brief description of the scheme of the proof. 
First, for the proof it is convenient to realize the processes $\mathcal
{W}$ and $\mathcal{Z}$ on an extension of $\mathcal{B}^{(0)}$ as in
Remarks~\ref{rmkW}--\ref{rmkZ} (so we will use the notation introduced
in these remarks in the following). For notational simplicity, we use
the same letters $P$ and $E$ for the probability and the expectation
with respect to this extension.
Next we introduce some notation. We denote by $\mathcal{R}_m$ the set
of all indices $r$ such that $S_r=S(m',j)$ for some $j\geq1$ and some
$m'\leq m$. Also, we set
\[
\cases{\displaystyle b(m)_t=b_t-\int_{A_m\cap\{z\dvt\llvert \delta(t,z)\rrvert
\leq1\}}
\delta(t,z)\lambda(\mathrm{d}z),\qquad B(m)_t=\int_0^tb(m)_s
\,\mathrm{d}s,\qquad
\vspace*{3pt}\cr
\displaystyle M_t=\int_0^t
\sigma_s\,\mathrm{d}W_s,
\vspace*{3pt}\cr
\displaystyle C(m)_t=X_0+B(m)_t+M_t,
\qquad J(m)_t=\delta1_{A_m}\star\mu_t,\qquad
X(m)_t=C(m)_t+J(m)_t,
\vspace*{3pt}\cr
\displaystyle Z(m)_t=X_t-X(m)_t=\delta1_{A_m^c}
\star(\mu-\nu)_t.}
\]
These processes are well-defined under~\textup{[{SA2}]}. Furthermore, set
$I_p=[T_{p-1},T_p)$ for every $p\in\mathbb{Z}_+$.
On the other hand, for any process $V$ and any (random) interval
$I=[S,T)$, we define the random variable $V(I)$ by $V(I)=V_T-V_S$. We
also set $I(t)=I\cap[0,t)=[S\wedge t,T\wedge t)$ for any $t\in\mathbb
{R}_+$ and $\llvert I\rrvert =T-S$. For any real-valued function $u$
on $[0,1]$, we set $u^n_p=u(p/k_n)$ for $p=0,1,\dots,k_n$. For any
$d$-dimensional processes $U$, $V$, any $k,l\in\{1,\dots,d\}$ and any
$u,v\in\{g,g'\}$, we define the process $\Xi^{(k,l)}_{u,v}(U,V)^n$ by
\begin{eqnarray*}
\Xi^{(k,l)}_{u,v}(U,V)^n_t=
\frac{1}{\psi_2k_n}\sum_{i=1}^{N^n_t-k_n+1}
\overline{U}(u)^k_i\overline{V}(v)^l_i,
\qquad t\in\mathbb{R}_+,
\end{eqnarray*}
where $\overline{U}(u)^k_i=\sum_{p=0}^{k_n-1}u^n_pU^k(I_{i+p})$ and
$\overline{V}(v)^l_i$ is defined analogously. Moreover, we define the
$d$-dimensional process $\mathfrak{E}$ by
\begin{eqnarray*}
\mathfrak{E}^k_t=-\frac{1}{k_n}\sum
_{p=1}^\infty\epsilon^k_{\tau^k_p}1_{\{\tau^k_p\leq t\}},
\qquad t\in\mathbb{R}_+, k=1,\dots,d.
\end{eqnarray*}
It can easily been seen that $\mathfrak{E}$ is a purely discontinuous
locally square-integrable martingale on $\mathcal{B}$ under~[SA3]. Finally, for any $d$-dimensional process $V$ we define the
$\mathbb{R}^d\otimes\mathbb{R}^d$-valued process $\bolds{\Xi}[V]^n$ by
\[
\bolds{\Xi}[V]^{n,kl}=\Xi^{(k,l)}_{g,g}(V,V)^n+
\Xi^{(k,l)}_{g,g'}(V,\mathfrak{E})^n +
\Xi^{(l,k)}_{g,g'}(V,\mathfrak{E})^n +
\Xi^{(k,l)}_{g',g'}(\mathfrak{E},\mathfrak{E})^n,\qquad
k,l=1,\dots,d.
\]


Now we turn to the outline of the proof. In the first step, we show
that the errors from the interpolations to the synchronized sampling
times are asymptotically negligible:

\begin{prop}\label{synchronization}
Assume~\textup{[SA1]}--\textup{[{SA3}]} and~\textup{[A4]}\textup{(ii)}. Then
$n^{1/4} (\MRC[Y]^{n}-\bolds{\Xi}[X]^{n}  +\frac{\psi_1}{\psi
_2k_n^2}[Y,Y]^{n} )\xrightarrow{ucp}0$.
\end{prop}

The proof of this proposition is an easy extension of that of
Proposition 6.1 from \cite{Koike2014time}, so we omit it.


In the next step, we decompose the quantity $\bolds{\Xi}[X]^n$ as
$\bolds{\Xi}[X]^n=\bolds{\Xi}[X(m)]^n_t+ (\bolds{\Xi}[X]^n-\bolds{\Xi
}[X(m)]^n )$ for each $m$, and show that the first term enjoys a
central limit theorem for any fixed~$m$ and the second term is
negligible as $m\to\infty$. More precisely, we prove the following propositions.

\begin{prop}\label{mainCLT}
Suppose that~\textup{[SA1]}--\textup{[{SA3}]} and~\textup{[A4]} are
satisfied. Then
\begin{eqnarray*}
n^{1/4} \biggl(\bolds{\Xi}\bigl[X(m)\bigr]^n_t-
\bigl[X(m),X(m)\bigr]_t-\frac{\psi_1}{\psi_2k_n^2}[Y,Y]^{n}_t
\biggr) \to^{d_s}\mathcal{W}_t+\mathcal{Z}(m)_t
\end{eqnarray*}
as $n\to\infty$ for any $t>0$ and any $m\geq1$, where $\mathcal
{Z}(m)_t=\sum_{r\in\mathcal{R}_m\dvt S_r\leq t}(\mathfrak{Z}_r+\mathfrak
{Z}_r^*)$.

When further $X$ is continuous, the processes $n^{1/4}(\bolds{\Xi
}[X]^n-[X,X]-\frac{\psi_1}{\psi_2k_n^2}[Y,Y]^{n})$ converge stably in
law to the process $\mathcal{W}$ for the Skorokhod topology.
\end{prop}


%
\begin{prop}\label{lemsmalljumps}
Suppose that~\textup{[SA1]}--\textup{[{SA3}]} and~\textup{[A4]}\textup{(ii)}
are satisfied. Then
%
\begin{equation}
\label{eqsmalljumps1} \mathcal{Z}(m)_t\to^{d_s}
\mathcal{Z}_t
\end{equation}
as $m\to\infty$ and
%
\begin{equation}
\label{eqsmalljumps2} \limsup_{m\to\infty}\limsup_{n\to\infty}P
\bigl(n^{1/4}\bigl\llVert\bolds{\Xi}[X]^n_t-
\bolds{\Xi}\bigl[X(m)\bigr]^n_t\bigr\rrVert>\eta
\bigr)=0
\end{equation}
for any $t,\eta>0$.
\end{prop}

Combining Propositions~\ref{synchronization}--\ref{lemsmalljumps} with
Proposition 2.2.4 of \cite{JP2012}, we obtain Theorem~\ref{mainthm}.


\subsection{Proof of Proposition \texorpdfstring{\protect\ref{mainCLT}}{6.2}}

Throughout the discussions, for (deterministic) sequences $(x_n)$ and
$(y_n)$, $x_n\lesssim y_n$ means that there is a (non-random) constant
$K\in[0,\infty)$ such that $x_n\leq Ky_n$ for large $n$. We also denote
by $E_0$ the conditional expectation given $\mathcal{F}^{(0)}$, that
is, $E_0[\cdot]=E[\cdot\mid\mathcal{F}^{(0)}]$.

The proof of Proposition~\ref{mainCLT} is divided into the following steps:
%
\begin{longlist}[(iii)]
\item[(i)] Approximating the estimation error due to the diffusion part by a
more tractable one, 

\item[(ii)] Proving a central limit theorem for the approximation constructed
in (i), 

\item[(iii)] Approximating the estimation error due to the jump part by a more
tractable one (Section~\ref{secappdiscont}),

\item[(iv)] Proving a local stable convergence result corresponding to Lemma
16.3.7 of \cite{JP2012} (Section~\ref{seclocalCLT}),

\item[(v)] Proving a joint limit theorem for the pair of the above
approximations and completing the proof of the proposition (Section~\ref
{secjointCLT}).
\end{longlist}

The first two steps have already been carried out in \cite
{Koike2014time}. For the later use, we summarize the result in the
following. We begin by introducing some notation. For any
$d$-dimensional processes $U,V$, any $k,l\in\{1,\dots,d\}$ and any
real-valued functions $u,v$ on $[0,1]$, we define the processes $\mathbb
{M}^{(k,l)}_{u,v}(U,V)^n$ and $\mathbb{L}^{(k,l)}_{u,v}(U,V)^n$ by
\begin{eqnarray*}
\mathbb{M}^{(k,l)}_{u,v}(U,V)^n_t=\sum
_{q=2}^{N^n_t+1}C^n_{u,v}(U)^k_qV^l(I_q),
\qquad\mathbb{L}^{(k,l)}_{u,v}(U,V)^n_t=
\mathbb{M}^{(k,l)}_{u,v}(U,V)^n_t+
\mathbb{M}^{(l,k)}_{v,u}(V,U)^n_t,
\end{eqnarray*}
where
\begin{eqnarray*}
C^n_{u,v}(U)^k_q=\sum
_{p=(q-k_n)\vee1}^{q-1}c^n_{u,v}(p,q)U^k(I_p),
\qquad c^n_{u,v}(p,q)=\frac{1}{\psi_2k_n}\sum
_{i=(p\vee q-k_n+1)\vee1}^{p\wedge q}u^n_{p-i}v^n_{q-i}.
\end{eqnarray*}
Moreover, define the $\mathbb{R}^d\otimes\mathbb{R}^d$-valued process
$\mathbf{L}[M]^n$ by
\begin{eqnarray*}
\mathbf{L}[M]^{n,kl}=\mathbb{L}^{(k,l)}_{g,g}(M,M)^n+
\mathbb{L}^{(k,l)}_{g,g'}(M,\mathfrak{E})^n +
\mathbb{L}^{(l,k)}_{g,g'}(M,\mathfrak{E})^n +
\mathbb{L}^{(k,l)}_{g',g'}(\mathfrak{E},\mathfrak{E})^n.
\end{eqnarray*}
Then, we have the following results, which are proved as Proposition
6.2 and equation (6.7) from~\cite{Koike2014time}:

\begin{prop}\label{appcont}
Suppose that~\textup{[SA1]}--\textup{[{SA3}]} and~\textup{[A4]}\textup{(ii)}
are satisfied. Then
\begin{eqnarray*}
n^{1/4} \biggl(\bolds{\Xi}\bigl[C(m)\bigr]^n-
\mathbf{L}[M]^n-[M,M]-\frac{\psi_1}{\psi_2k_n^2}[Y,Y]^{n} \biggr)
\xrightarrow{ucp}0
\end{eqnarray*}
as $n\to\infty$ for any $m\geq1$.
\end{prop}

%
\begin{prop}\label{contCLT}
Suppose that~\textup{[SA1]}--\textup{[{SA3}]} and~\textup{[A4]}\textup{(ii)}
are satisfied. Then the processes $n^{1/4}\mathbf{L}[M]^n$ converge
stably in law to $\mathcal{W}$ for the Skorokhod topology.
\end{prop}

From the next section, we start the proofs of the remaining steps.


\subsubsection{Approximation of the estimation error due to the jump
part}\label{secappdiscont}

In this subsection we fix $t>0$ and $m\in\mathbb{N}$, and denote by
$\Omega_n(t,m)$ the set on which $k_n-1\leq N^n_{S_r-}\leq N^n_t-k_n$
for all $r\in\mathcal{R}_m$ such that $S_r\leq t$. On this set, we have
%
\begin{eqnarray}\label{CJformula}
n^{1/4}\Xi^{(k,l)}_{g,g}\bigl(C(m),J(m)
\bigr)^n_t &=&\frac{n^{1/4}}{\psi_2k_n}\sum
_{i=0}^{N^n_t-k_n+1}\sum_{p,q=0}^{k_n-1}g^n_p
g^n_qC(m)^k(I_{i+p})J(m)^l(I_{i+q})
\nonumber\\[-8pt]\\[-8pt]\nonumber
&=&\sum_{r\in\mathcal{R}_m\dvt S_r\leq t} \bigl\{
\eta_+(n,r)^k+\eta_-(n,r)^k \bigr\}\Delta
X^l_{S_r},
\end{eqnarray}
where
\[
\cases{ \displaystyle\eta_+(n,r)=n^{1/4}\sum_{p=N^{n}_{S_r-}+1}^{N^{n}_{S_r-}+k_n}c^n_{g,g}
\bigl(p,N^{n}_{S_r-}+1\bigr)C(m) (I_p),
\vspace*{3pt}\cr
\displaystyle\eta_-(n,r)=n^{1/4}\sum_{p=(N^{n}_{S_r-}-k_n+2)_+}^{N^{n}_{S_r-}}c^n_{g,g}
\bigl(p,N^{n}_{S_r-}+1\bigr)C(m) (I_p).}
\]
Similarly, on $\Omega_n(t,m)$ we have
%
\begin{equation}
\label{EJformula} n^{1/4}\Xi^{(k,l)}_{g',g}\bigl(
\mathfrak{E},J(m)\bigr)^n_t =\sum
_{r\in\mathcal{R}_m\dvt S_r\leq t} \bigl\{\eta'_+(n,r)^k+
\eta'_-(n,r)^k \bigr\}\Delta X^l_{S_r},
\end{equation}
where
\[
\cases{ \displaystyle\eta'_+(n,r)^k=-\frac{n^{1/4}}{k_n}\sum
_{p=N^{n}_{S_r-}+1}^{N^{n}_{S_r-}+k_n}c^n_{g',g}
\bigl(p,N^{n}_{S_r-}+1\bigr)\epsilon^k_{\tau^k_p},
\vspace*{3pt}\cr
\displaystyle\eta'_-(n,r)^k=-\frac{n^{1/4}}{k_n}\sum
_{p=(N^{n}_{S_r-}-k_n+2)_+}^{N^{n}_{S_r-}}c^n_{g',g}
\bigl(p,N^{n}_{S_r-}+1\bigr)\epsilon^k_{\tau^k_p}.}
\]
The aim of this subsection is to approximate $\eta_{\pm}(n,r)$ and $\eta
'_{\pm}(n,r)$ by more tractable quantities. Here,\vspace*{1pt} the major difficulty
coming from the irregularity of the observation times is the fact that
$N^n_{S_r-}-k_n+1$ might not be {a $(\mathcal
{G}^{(0)}_{T_p})_{p=0}^\infty$}-stopping time. Therefore,\vspace*{-1pt} we first
``approximate'' $N^n_{S_r-}-k_n+1$ by {a $(\mathcal
{G}^{(0)}_{T_p})_{p=0}^\infty$}-stopping time.

More precisely, set $\underline{S}_r=(S_r-\frac{k_n}{n}\log n)_+$ and
$S^\dagger_r=(S_r-\frac{k_n}{n}G^n_{\underline{S}_r}\wedge\log n)_+$.
Then, $S^\dagger_r$ is {a $(\mathcal{G}^{(0)}_{t})$}-stopping time
(Lemma~\ref{lemstoptime}), and thus $N^n_{S^\dagger_r}+1$ is {a
$(\mathcal{G}^{(0)}_{T_p})_{p=0}^\infty$}-stopping time, and this
variable gives an approximation of $N^n_{S_r-}-k_n+1$ (Lemma~\ref{appstoptime}).

Now we can define our tractable approximations of $\eta_{\pm}(n,r)$ and
$\eta'_{\pm}(n,r)$ as follows. For any non-negative random variable
$\rho$ and any real-valued function $\phi$ on $[0,1]$, we define the
$d'$-dimensional variable $L(\phi,\rho)_n=(L(\phi,\rho)^j_n)_{1\leq
j\leq d'}$ and the $d$-dimensional variable $L'(\phi,\rho)_n=(L^{\prime
}(\phi,\rho)^k_n)_{1\leq k\leq d}$ by
\begin{eqnarray*}
L(\phi,\rho)^j_n=n^{1/4}\sum
_{w=1}^{k_n-1}\phi^n_wW^j(I_{i(\rho)^n+w}),
\qquad L'(\phi,\rho)^k_n=\frac{n^{1/4}}{k_n}
\sum_{w=1}^{k_n-1}\phi^n_w
\epsilon^k_{\tau^k_{i(\rho)^n+w}},
\end{eqnarray*}
where we set $i(\rho)^n=N^n_\rho+1$ (recall that $\phi^n_w=\phi
(w/k_n)$). We also define the function $\widetilde{\phi}$ on $[0,1]$ by
$\widetilde{\phi}(x)=\phi(1-x)$. Then we set
\[
\cases{ z^n_{r-}=\psi_2^{-1}L
\bigl(\widetilde{\phi_{g,g}},S^\dagger_r
\bigr)_n, &\quad$z^n_{r+}=
\psi_2^{-1}L(\phi_{g,g},S_r)_n$,
\vspace*{3pt}\cr
z'^n_{r-}=-\psi_2^{-1}L'
\bigl(\widetilde{\phi_{g',g}},S^\dagger_r
\bigr)_n, &\quad$z'^n_{r+}=-
\psi_2^{-1}L'(\phi_{g',g},S_r)_n$.}
\]
The aim of this subsection is to prove the following proposition.

\begin{prop}\label{JPlem1218}
Suppose that~\textup{[SA1]}--\textup{[{SA3}]} and~\textup{[A4]}\textup{(ii)}
are satisfied. Then
%
\begin{eqnarray}
\eta_-(n,r)&=&\sigma_{S^\dagger_r}z^{n}_{r-}+\mathrm{o}_p(1),
\qquad\eta_+(n,r)=\sigma_{S_r}z^{n}_{r+}+\mathrm{o}_p(1),\label{JPlem1218eq1}
\\
\eta'_-(n,r)&=&z'^n_{r-}+\mathrm{o}_p(1),
\qquad\eta'_+(n,r)=z'^n_{r+}+\mathrm{o}_p(1).\label{JPlem1218eq2}
\end{eqnarray}
\end{prop}

Now we start to justify that the variable $N^n_{S^\dagger_r}+1$ is an
appropriate approximation of $N^n_{S_r-}-k_n+1$. In the remainder of
this subsection we fix an index $r\in\mathcal{R}_m$ such that
$S_r<\infty$.

\begin{lem}\label{lemstoptime}
Under~\textup{[SA1]}, $S^\dagger_r$ is {a $(\mathcal
{G}^{(0)}_t)$}-stopping time.
\end{lem}

\begin{pf}
For any $t\geq0$, we have
$\{S^\dagger_r\leq t\}
= \{S_r\leq t+(k_n/n)G^n_{\underline{S}_r}\wedge\log n \}
\cap\{\underline{S}_r\leq t \}$.
Therefore, noting that $S_r$ is {$\mathcal{G}^{(0)}_0$}-measurable, we
obtain $\{S^\dagger_r\leq t\}\in{\mathcal{G}^{(0)}_t}$.
\end{pf}

%
\begin{lem}\label{left}
Under~\textup{[SA1]}, $\sup_{0<h<h_0}\llvert
G_{S_r-}-G_{(S_r-h)_+}\rrvert =\mathrm{O}_p(\sqrt{h_0})$ as $h_0\downarrow0$.
\end{lem}

\begin{pf}
Define the processes $G(m)$, $G'(m)$ and $G''(m)$ by
$G(m)_t=\int_0^t\widehat{\sigma}_s\,\mathrm{d}W_s+(\widehat{\delta
}1_{A_m^c\cap\{\llvert \widetilde{\delta}\rrvert \leq1\}})\star(\mu
-\nu)_t$,
$G'(m)_t=(\widehat{\delta}1_{A_m\cap\{\llvert \widetilde{\delta
}\rrvert \leq1\}}+\widehat{\delta}1_{\{\llvert \widetilde{\delta
}\rrvert >1\}})\star\mu_t$,
and $G''(m)=G-G(m)-G'(m)$. Since $G'(m)$ is piecewise constant, it is
evident that $\sup_{0<h<h_0}\llvert
G'(m)_{S_r-}-G'(m)_{(S_r-h)_+}\rrvert =\mathrm{O}_p(\sqrt{h_0})$. Moreover,
since $G''(m)$ is absolutely continuous with a locally bounded
derivative, it also holds that $\sup_{0<h<h_0}\llvert
G''(m)_{S_r-}-G''(m)_{(S_r-h)_+}\rrvert =\mathrm{O}_p(\sqrt{h_0})$. On the
other hand, let $(\mathcal{G}^{A_m}_t)$ be the smallest filtration
containing $(\mathcal{F}^{(0)}_t)$ {such that $\mathcal{G}^{A_m}_0$
contains the $\sigma$-field generated by the restriction of the measure
$\mu$ to $\mathbb{R}_+\times A_m$}. Then, by Proposition 2.1.10 of \cite
{JP2012} $G(m)$ is a locally square integrable martingale with
{respect} to $(\mathcal{G}^{A_m}_t)$ and its predictable quadratic
variation is given by $\langle G(m)\rangle=\int_0^t\widehat{\sigma
}_s\widehat{\sigma}_s^*\,\mathrm{d}s+(\widetilde{\delta}^21_{A_m^c\cap\{
\llvert \widetilde{\delta}\rrvert \leq1\}})\star\nu$, and
$G(m)_{S_r-}=G(m)_{S_r}$. Since $S_r$ is {$\mathcal
{G}^{A_m}_0$}-measurable, $(S_r-h)_+$ is a $(\mathcal
{G}^{A_m}_t)$-stopping time for every $h\geq0$. Therefore, the Lenglart
inequality implies that
\begin{eqnarray*}
&& P \Bigl(\sup_{0\leq h\leq h_0}\bigl\llvert h_0^{-1/2}
\bigl\{G(m)_{S_r}-G(m)_{(S_r-h)_+}\bigr\}\bigr\rrvert
^2>K \Bigr)
\\
&&\quad \leq\frac{K'}{K}+P \bigl(h_0^{-1}
\bigl\llvert\bigl\langle G(m)\bigr\rangle_{S_r}-\bigl\langle G(m)\bigr
\rangle_{(S_r-h_0)_+}\bigr\rrvert>K' \bigr)
\end{eqnarray*}
for any $K,K'>0$, and thus a standard localization argument yields $\sup
_{0\leq h\leq h_0}\llvert G(m)_{S_r}-G(m)_{(S_r-h)_+}\rrvert
=\mathrm{O}_p(\sqrt{h_0})$. This completes the proof of the lemma.
\end{pf}



%
\begin{lem}\label{appstoptime}
Under~\textup{[SA1]}, $N^n_{S_r-}-N^n_{S^\dagger_r}=k_n+\mathrm{o}_p(n^{\sklfrac
{1}{2}-\alpha'})$ for any $\alpha'\in(0,(\xi-\kappa-\frac{1}{2})\wedge
(\xi-\frac{3}{4})\wedge\varpi)$.
\end{lem}

\begin{pf}
Since $\alpha'<\xi-\kappa-\frac{1}{2}$, by~\textup{[SA1]}, we have
\[
N^n_{S_r-}-N^n_{S^\dagger_r}=\sum
_{p=1}^{N^n_{S_r-}+1}\frac{E [n\llvert I_p\rrvert \mid{\mathcal
{G}^{(0)}_{T_{p-1}}} ]}{G^n_{T_{p-1}}}1_{\{T_{p-1}>S^\dagger_r\}}+\mathrm{o}_p
\bigl(n^{1/2-\alpha'}\bigr).
\]
{In particular, from this expression and~\textup{[SA1]}, we deduce
$N^n_{S_r-}-N^n_{S^\dagger_r}=\mathrm{O}_p(\sqrt{n}\log n)$.} Therefore,~\textup{[SA1]} yields
\begin{eqnarray*}
\sum_{p=1}^{N^n_{S_r-}+1}E \biggl[\biggl\llvert
n^{\alpha'-\sfrac{1}{2}}\frac{n\llvert I_p\rrvert }{G^n_{T_{p-1}}}1_{\{
T^{p-1}>S^\dagger_r\}}\biggr\rrvert^2
\Big| {\mathcal{G}^{(0)}_{T_{p-1}}} \biggr] =n^{1+2\alpha'}
\bar{r}_n^2\sum_{p=1}^{N^n_{S_r-}+1}
\frac{1}{G^n_{T_{p-1}}}1_{\{T^{p-1}>S^\dagger_r\}}=\mathrm{o}_p(1),
\end{eqnarray*}
hence {Lemma 2.3 of \cite{Fu2010b}} implies that
\begin{eqnarray*}
N^n_{S_r-}-N^n_{S^\dagger_r}=n\sum
_{p=1}^{N^n_{S_r-}+1}\frac{\llvert I_p\rrvert }{G^n_{T_{p-1}}}1_{\{
T^{p-1}>S^\dagger_r\}}+\mathrm{o}_p
\bigl(n^{\sklfrac{1}{2}-\alpha'} \bigr).
\end{eqnarray*}
Now~\textup{[SA1]}, Lemma~\ref{left} and the fact that $\alpha'<\varpi
\wedge\frac{1}{4}$ yield
\begin{eqnarray*}
n\sum_{p=1}^{N^n_{S_r-}+1} \biggl(
\frac{\llvert I_p\rrvert }{G^n_{T_{p-1}}}-\frac{\llvert I_p\rrvert
}{G^n_{\underline{S}_{r}}} \biggr)1_{\{T^{p-1}>S^\dagger_r\}} &=&n\sum
_{p=1}^{N^n_{S_r-}+1} \biggl(\frac{\llvert I_p\rrvert }{G_{T_{p-1}}}-
\frac{\llvert I_p\rrvert }{G_{\underline{S}_{r}}} \biggr)1_{\{
T^{p-1}>S^\dagger_r\}}+\mathrm{o}_p \bigl(n^{\sklfrac{1}{2}-\alpha'}
\bigr)
\\
&=& \mathrm{o}_p \bigl(n^{\sklfrac{1}{2}-\alpha'} \bigr),
\end{eqnarray*}
thus we have
\begin{eqnarray*}
N^n_{S_r-}-N^n_{S^\dagger_r}=n\sum
_{p=1}^{N^n_{S_r-}+1}\frac{\llvert I_p\rrvert }{G^n_{\underline
{S}_r}}1_{\{T^{p-1}>S^\dagger_r\}}+\mathrm{o}_p
\bigl(n^{\sklfrac{1}{2}-\alpha'} \bigr) =k_n\frac{G^n_{\underline
{S}_r}\wedge\log n}{G^n_{\underline{S}_r}}+\mathrm{o}_p
\bigl(n^{\sklfrac{1}{2}-\alpha'} \bigr).
\end{eqnarray*}
Since $\lim_nP(G^n_{\underline{S}_r}>\log n)=0$, we obtain the desired result.
\end{pf}

Now we proceed to the main body of the proof of Proposition~\ref
{JPlem1218}. Denote by $\Omega_{n}(m)$ the set on which $\llvert
S_{r_1}-S_{r_2}\rrvert >(k_n/n)\log n$ for any $r_1,r_2\in\mathcal
{R}_m$ such that $r_1\neq r_2$ and $S_{r_1},S_{r_2}<\infty$. Since
$S_{r_1}\neq S_{r_2}$ if $r_1\neq r_2$ and $S_{r_1}, S_{r_2}<\infty$,
we have $P(\Omega_{n}(m))\to1$ as $n\to\infty$.

\begin{lem}\label{sigmalc}
Under~\textup{[{SA2}]}, $E[\sup_{\underline{S}_r\leq s<S_r}\llVert \sigma
_s-\sigma_{\underline{S}_r}\rrVert ^2;\Omega_{n}(m)]\lesssim
(k_n/n)\log n$.
\end{lem}

\begin{pf}
Since no jump of the Poisson process $1_{A_m}\star\mu$ occurs in
$[\underline{S}_r,S_r)$ on the set $\Omega_{n}(m)$, we have $\sigma
_s=\sigma(m)_s$ for every $s\in[\underline{S}_r,S_r)$ on this set, where
\[
\cases{\displaystyle \sigma(m)_s=\sigma_0+\int_0^s
\widetilde{b}(m)_u\,\mathrm{d}u+\int_0^s
\widetilde{\sigma}_u\,\mathrm{d}W_u +(\widetilde{
\delta}1_{A_m^c})\star(\mu-\nu)_s,
\vspace*{3pt}\cr
\displaystyle\widetilde{b}(m)_u=
\widetilde{b}_u-\int_{A_m\cap\{z\dvt\llvert \widetilde{\delta
}(u,z)\rrvert \leq1\}}\widetilde{
\delta}(u,z)\lambda(\mathrm{d}z).}
\]
On the other hand, by Proposition 2.1.10 of \cite{JP2012} we have
$E [\sup_{\underline{S}_r\leq s<S_r}\llVert \sigma(m)_s-\sigma
(m)_{\underline{S}_r}\rrVert ^2 ]
\lesssim\frac{k_n}{n}\log n$,
which implies the desired result.
\end{pf}

\begin{pf*}{Proof of Proposition~\ref{JPlem1218}}
Throughout the proof, we fix a constant $\alpha'$ such that $1-\xi
<\alpha'<(\xi-\kappa-\frac{1}{2})\wedge(\xi-\frac{3}{4})\wedge\varpi$.
Such an $\alpha'$ exists due to $(\ref{estxi})$.

First we prove the first equation of $(\ref{JPlem1218eq1})$. Set $\Omega
_n=\{N^n_{S_{r-}}-k_n+1\geq0\}$. By the Lipschitz continuity of $g$,
\begin{eqnarray*}
\eta_-(n,r)^k&=&\frac{n^{1/4}}{\psi_2}\sum_{p=N^{n}_{S_r-}-k_n+2}^{N^{n}_{S_r-}}(
\phi_{g,g})^n_{N^{n}_{S_r-}+1-p}M^k(I_p)+\mathrm{o}_p(1)
\\
&=&\frac{n^{1/4}}{\psi_2}\sum_{w=1}^{k_n-1}(
\widetilde{\phi_{g,g}})^n_wM^k(I_{N^n_{S_r-}+1-k_n+w})+\mathrm{o}_p(1)
\end{eqnarray*}
on $\Omega_n$. On the other hand, noting that we have
\[
\sup\bigl\{\llVert M_s-M_r\rrVert\dvt\llvert
s-r\rrvert\leq h,s,r\in[0,t] \bigr\}=\mathrm{O}_p \bigl(\sqrt{h}\llvert\log
h\rrvert\bigr)
\]
as $h\downarrow0$ for any $t>0$ due to a representation of a continuous
local martingale with Brownian motion and L\'{e}vy's theorem on the
uniform modulus of continuity of Brownian motion, summation by parts,
$(\ref{SA4})$ and Lemma~\ref{appstoptime} imply that
\begin{eqnarray*}
&& n^{1/4}\sum_{w=1}^{k_n-1}(
\widetilde{\phi_{g,g}})^n_w \bigl\{M^k(I_{N^n_{S_r-}+1-k_n+w})-M^k(I_{i(S^\dagger_r)^n+w}) \bigr\}
\\
&&\quad = n^{1/4}\sum_{w=1}^{k_n-2} \bigl\{(
\widetilde{\phi_{g,g}})^n_w-(\widetilde{
\phi_{g,g}})^n_{w+1} \bigr\} \bigl
(M^k_{T_{N^n_{S_r-}+1-k_n+w}}-M^k_{T_{i(S^\dagger_r)^n+w}}
\bigr)
\\
&&\qquad{} +n^{1/4}(\widetilde{\phi_{g,g}})^n_{k_n-1}
\bigl(M^k_{T_{N^n_{S_r-}}}-M^k_{T_{i(S^\dagger_r)^n+k_n-1}} \bigr)
-n^{1/4}(\widetilde{\phi_{g,g}})^n_1
\bigl(M^k_{T_{N^n_{S_r-}+1-k_n}}-M^k_{T_{i(S^\dagger_r)^n}} \bigr)
\\
&&\quad =\mathrm{o}_p \bigl(n^{1/4}\sqrt{n^{1/2-\alpha'-\xi}\log n} \bigr)
=\mathrm{o}_p \bigl(n^{(1-\xi-\alpha')/2}\sqrt{\log n} \bigr)=\mathrm{o}_p(1)
\end{eqnarray*}
on $\Omega_n$. Since $\lim_nP(\Omega_n)=1$, we conclude that
\begin{eqnarray*}
\eta_-(n,r)^k =\frac{n^{1/4}}{\psi_2}\sum_{w=1}^{k_n-1}(
\widetilde{\phi_{g,g}})^n_wM^k(I_{i(S^\dagger_r)+w})+\mathrm{o}_p(1).
\end{eqnarray*}
Next, noting that $W$ is a $d'$-dimensional {$(\mathcal
{G}^{(0)}_t)$}-Brownian motion (recall that {$(\mathcal{G}^{(0)}_t)$}
is the smallest filtration containing $(\mathcal{F}^{(0)}_t)$ such that
{$\mathcal{G}^{(0)}_0$ contains the $\sigma$-field generated by $\mu
$}), we have
\begin{eqnarray*}
&& n^{1/4}\sum_{w=1}^{k_n-1}(
\widetilde{\phi_{g,g}})^n_wM^k(I_{i(S^\dagger_r)+w})-
\sum_{j=1}^{d'}\sigma^{kj}_{\underline{S}_r}z^{n,j}_{r-}
\\
&&\quad =n^{1/4}\sum_{j=1}^{d'}\sum
_{w=1}^{k_n-1}(\widetilde{\phi_{g,g}})^n_w
\int_{T_{i(S^\dagger_r)^n+w-1}}^{T_{i(S^\dagger_r)^n+w}} \bigl(\sigma^{kj}_s-
\sigma^{kj}_{\underline{S}_r} \bigr)\,\mathrm{d}W^j_s,
\end{eqnarray*}
hence the Lenglart inequality implies that it is enough to show that
\begin{eqnarray*}
\Delta_n:=\sqrt{n}\sum_{j=1}^{d'}
\sum_{w=1}^{k_n-1}\bigl\llvert(\widetilde{
\phi_{g,g}})^n_w\bigr\rrvert^2E
\biggl[\int_{T_{i(S^\dagger_r)^n+w-1}}^{T_{i(S^\dagger_r)^n+w}} \bigl
(\sigma^{kj}_s-
\sigma^{kj}_{\underline{S}_r} \bigr)^2\,\mathrm{d}s\Big| {
\mathcal{G}^{(0)}_{T_{i(S^\dagger_r)^n+w-1}}} \biggr]\to^p0.
\end{eqnarray*}
Set
\begin{eqnarray*}
\Delta_n' =\sqrt{n}\sum_{j=1}^{d'}
\sum_{w=1}^{k_n-1}\bigl\llvert(\widetilde{
\phi_{g,g}})^n_w\bigr\rrvert^2E
\biggl[\int_{T_{i(S^\dagger_r)^n+w-1}}^{T_{i(S^\dagger_r)^n+w}} \bigl
(\sigma^{kj}_s-
\sigma^{kj}_{\underline{S}_r} \bigr)^2\,\mathrm{d}s\Big|  {
\mathcal{G}^{(0)}_{T_{i(S^\dagger_r)^n+w-1}}} \biggr]1_{ \{
T_{i(S^\dagger_r)^n+w-1}\leq S_r \}}.
\end{eqnarray*}
Then, since $\Omega_{n}(m)\in{\mathcal{G}^{(0)}_0}$, it holds that
\begin{eqnarray*}
E \bigl[\Delta'_n;\Omega_{n}(m)
\bigr]
\lesssim\sqrt{n}\sum_{j=1}^{d'}E
\biggl[\int_{T_{i(S^\dagger_r)^n}}^{T_{i(S_r)^n}} \bigl(\sigma^{kj}_s-
\sigma^{kj}_{S^\dagger_r} \bigr)^2\,\mathrm{d}s;
\Omega_{n}(m) \biggr].
\end{eqnarray*}
Now, Lemma~\ref{sigmalc}, the boundedness of $\sigma$ and $(\ref{SA4})$
imply that
\begin{eqnarray*}
E \biggl[\int_{T_{i(S^\dagger_r)^n}}^{T_{i(S_r)^n}} \bigl(
\sigma^{kj}_s-\sigma^{kj}_{S^\dagger_r}
\bigr)^2\,\mathrm{d}s;\Omega_{n}(m) \biggr] 
&\lesssim&\frac{k_n}{n}(\log n)E \Bigl[\sup_{\underline{S}_r\leq s<S_r}
\bigl(
\sigma^{kj}_s-\sigma^{kj}_{S^\dagger_r}
\bigr)^2;\Omega_{n}(m) \Bigr] +\bar{r}_n
\lesssim\bar{r}_n,
\end{eqnarray*}
hence we obtain $E [\Delta'_n;\Omega_{n}(m) ]\lesssim\sqrt{n}\bar
{r}_n=\mathrm{o}(1)$. Therefore, the equation $\lim_nP(\Omega_{n}(m))=1$ and the
Chebyshev inequality yield $\Delta'_n=\mathrm{o}_p(1)$. On the other hand, the
boundedness of $\sigma$, $(\ref{SA4})$ and Lemma~\ref{appstoptime}
imply that
$
\llvert \Delta_n-\Delta'_n\rrvert
=\mathrm{o}_p(n^{1-\xi-\alpha'})=\mathrm{o}_p(1)$.
Consequently, we obtain $\Delta_n=\mathrm{o}_p(1)$ and the first equation of
$(\ref{JPlem1218eq1})$ has been proved. On the other hand, noting that
$N^n_{S_r}-N^n_{S_r-}\leq1$ and $S_r$ is an {$(\mathcal
{F}^{(0)}_t)$}-stopping time, the second equation of $(\ref
{JPlem1218eq1})$ can be shown in a similar (and simpler) manner.

Next, we prove the first equation of $(\ref{JPlem1218eq2})$. By the
(piecewise) Lipschitz continuity of $g$ and $g'$, we have on $\Omega_n$
\begin{eqnarray*}
\eta'_-(n,r)^k=-\frac{n^{1/4}}{\psi_2k_n}\sum
_{p=(N^{n}_{S_r-}-k_n+2)_+}^{N^{n}_{S_r-}}(\phi_{g',g})^n_{N^n_{S_r-}+1-p}
\epsilon^k_{\tau^k_p}+\mathrm{o}_p(1).
\end{eqnarray*}
Moreover, by Lemma~\ref{appstoptime},~\textup{[{SA3}]} and the Lipschitz
continuity of $\phi_{g',g}$ we have
\begin{eqnarray*}
&& E_0 \Biggl[\Biggl\llvert\frac{n^{1/4}}{k_n} \Biggl\{\sum
_{p=(N^{n}_{S_r-}-k_n+2)_+}^{N^{n}_{S_r-}}(\phi_{g',g})^n_{N^n_{S_r-}+1-p}
\epsilon^k_{\tau^k_p}-\sum_{p=i(S^\dagger_r)^n+1}^{i(S^\dagger_r)^n+k_n-1}(
\phi_{g',g})^n_{i(S^\dagger_r)^n+k_n-p}\epsilon^k_{\tau^k_p}
\Biggr\}\Biggr\rrvert^2 \Biggr]
\\
&&\quad  =\mathrm{O}_p
\bigl(n^{-\alpha'} \bigr)
\end{eqnarray*}
on $\Omega_n$. Since $\lim_nP(\Omega_n)=1$, we conclude that
\begin{eqnarray*}
\eta'_-(n,r)^k=-\frac{n^{1/4}}{\psi_2k_n}\sum
_{p=i(S^\dagger_r)^n+1}^{i(S^\dagger_r)^n+k_n-1}(\phi
_{g',g})^n_{i(S^\dagger_r)^n+k_n-p}
\epsilon^k_{\tau^k_p}+\mathrm{o}_p(1) =z'^{n,k}_{r-}+\mathrm{o}_p(1).
\end{eqnarray*}
Similarly, we can prove the second equation of $(\ref{JPlem1218eq2})$.
\end{pf*}


\subsubsection{An auxiliary local stable convergence result}\label{seclocalCLT}

In this subsection, we prove an auxiliary local stable convergence
result corresponding to Lemma~16.3.7 of \cite{JP2012}. The proof is
close to that of the aforementioned lemma, but there is a difference
due to the additional randomness coming from the sampling times.
Furthermore, we can also simplify some parts of the proof because it is
sufficient for our purpose to prove a simpler consequence than that of
the aforementioned lemma. For these reasons, we give a complete proof.

The following lemma is a direct consequence of the Skorokhod
representation theorem, so we omit the proof:

\begin{lem}\label{lemSk}
Let $(f_n)$ be a sequence of real-valued functions on $\mathbb{R}^D$
such that there exists a constant $K$ satisfying $\llvert
f_n(x)\rrvert \leq K$ and $\llvert f_n(x)-f_n(y)\rrvert \leq
K\llVert x-y \rrVert $ for all $x,y\in\mathbb{R}^D$ and every~$n$. If
a sequence $(x_n)$ of $\mathbb{R}^D$-valued random variables converges
in law to a variable $x$, then $E[f_n(x_n)]-E[f_n(x)]\to0$.
\end{lem}

The following lemma is the main result of this subsection. We denote by
$\mathfrak{N}_D$ the $D$-dimensional standard normal distribution.

\begin{lem}\label{JPlem1637}
Assume that $[\mathrm{SA}1]$, $[\mathrm{SA}3]$ and $[\mathrm
{A}4](\mathrm{ii})$ are satisfied. Suppose that for each $n$ there is
{a $(\mathcal{G}^{(0)}_t)$}-stopping time $\rho_n$. Suppose also that
there is a finite-valued variable $\rho$ such that \mbox{$\rho_n\to\rho$} as
$n\to\infty$ and one of the following two condition is satisfied:
%
\begin{equation}
\label{JPeq16318} \lleft.
\begin{array} {l@{\quad}l} (1) &\mbox{$\rho>0$,
$P(T_{i(\rho_n)^n+k_n-\lfloor n^{\beta}\rfloor}<\rho)\to1$ as $n\to
\infty$ for some $\beta\in(0,\xi-1/2)$,}
\\
& \mbox{in which case we set $G_{(\rho)}=G_{\rho-}$, $\widetilde{
\upsilon}_{(\rho)}=\widetilde{\upsilon}_{\rho-}$ and {$
\mathcal{G}^{(0)}_{(\rho)}=\mathcal{G}^{(0)}_{\rho-}$},}
\\
(2)&\mbox{$\rho_n\geq\rho$ for all $n$, in which case we set
$G_{(\rho)}=G_{\rho}$, $\widetilde{\upsilon}_{(\rho)}=
\widetilde{\upsilon}_{\rho}$ and {$\mathcal{G}^{(0)}_{(\rho)}=
\mathcal{G}^{(0)}_{\rho}$}}
\end{array}
\rright\}.\quad
\end{equation}
Let $\phi_1$ and $\phi_2$ be continuous real-valued functions $\phi_1$
and $\phi_2$ on $[0,1]$.
Then, for any $\mathcal{F}$-measurable bounded variable $U$ and any
bounded Lipschitz function $f$ on $\mathbb{R}^{d'+d}$ we have
%
\begin{eqnarray}\label{eqlocalCLT}
&& E\bigl[Uf\bigl(L_n,L'_n
\bigr)\mid\mathcal{G}_{\rho_n}\bigr]
\nonumber\\[-8pt]\\[-8pt]\nonumber
&&\quad \to^p E \biggl[U\int f
\bigl(\llVert\phi_1\rrVert\sqrt{\theta G_{(\rho)}}x,\llVert
\phi_2\rrVert\sqrt{\theta^{-1}}\widetilde{
\upsilon}_{(\rho)}y \bigr)\mathfrak{N}_{d'}(\mathrm{d}x)
\mathfrak{N}_d(\mathrm{d}y)\Big| {\mathcal{G}^{(0)}_{(\rho)}}
\biggr],
\end{eqnarray}
where $L_n=L(\phi_1,\rho_n)_n$, $L'_n=L'(\phi_2,\rho_n)_n$ and
$\llVert \phi_j\rrVert ^2=\int_0^1\phi_j(x)^2\,\mathrm{d}x$ for $j=1,2$.
\end{lem}

\begin{pf}
 \textit{Step}~1. For $k,l=1,\dots,d$ we set
$D^{kl}_n=\frac{1}{k_n}\sum_{w=1}^{k_n-1}\llvert (\phi
_{2})^n_w\rrvert ^21_{\{\tau^k_{i(\rho_n)^n+w}=\tau^l_{i(\rho_n)^n+w}\}}$.
We begin by proving $D^{kl}_n\to^p\llVert \phi_{2}\rrVert ^2\chi
^{kl}_{(\rho)}$, where\vspace*{1pt} we set $\chi^{kl}_{(\rho)}=\chi^{kl}_{\rho-}$ in
case~(1) and $\chi^{kl}_{(\rho)}=\chi^{kl}_{\rho}$ in case~(2). Since
$i(\rho_n)$ is {a $(\mathcal{G}^{(0)}_{T_p})_{p=0}^\infty$}-stopping
time, [SA1] and {Lemma 2.3 of \cite{Fu2010b}} yield
$D^{kl}_n=\frac{1}{k_n}\sum_{w=1}^{k_n-1}\llvert (\phi
_{2})^n_w\rrvert ^2\chi^{kl}_{T_{i(\rho_n)^n+w-1}}+\mathrm{o}_p(1)$.
Since $\chi^{kl}$ is cadlag, $(\ref{JPeq16318})$ implies that
$D^{kl}_n=\llVert \phi_{2}\rrVert ^2\chi^{kl}_{(\rho)}+\mathrm{o}_p(1)$.\vspace*{2pt}

 \textit{Step} 2. 
From\vspace*{1pt} step~1, by considering an appropriate subsequence if necessary,
without loss of generality\vspace*{1pt} we may assume that there is a subset $\Omega
_0$ of $\Omega^{(0)}$ such that $P^{(0)}(\Omega_0)=1$ and
$D^{kl}_n(\omega^{(0)})\to\llVert \phi_{2}\rrVert ^2\chi^{kl}_{(\rho
)}(\omega^{(0)})$ for all $\omega^{(0)}\in\Omega_0$.

 \textit{Step} 3. Fix\vspace*{1pt} $\omega^{(0)}\in\Omega_0$, and consider
the probability space $(\Omega^{(1)},\mathcal{F}^{(1)},Q_0)$, where
$Q_0(\cdot)=Q(\omega^{(0)},\cdot)$. Our aim in this step is to show
that under $Q_0$
%
\begin{equation}
\label{noiseCLT} L'_n\to^d\llVert
\phi_2\rrVert\sqrt{\theta^{-1}}\widetilde{
\upsilon}_{(\rho)}\bigl(\omega^{(0)}\bigr)\zeta',
\end{equation}
where $\zeta'$ is a standard $d$-dimensional normal variable
independent of $\mathcal{F}$.

For each $w=1,\dots, k_n-1$ we define the $d$-dimensional variable
$y^n_w=(y^{n,k}_w)_{1\leq k\leq d}$ by
\[
y^{n,k}_w=\frac{n^{1/4}}{k_n}(\phi_{2})^n_w\epsilon^k_{\tau^k_{i(\rho
_n)^n(\omega^{(0)})+w}}.
\]
Then $y^n_1,\dots,y^n_{k_n-1}$ are independent under $Q_0$ and we have
$L'_n(\omega^{(0)},\cdot)=\sum_{w=1}^{k_n-1}y^{n}_w$. Moreover, by
[SA3] we have
\begin{eqnarray*}
E_{Q_0}\bigl(y^n_w\bigr)&=&0,\qquad
E_{Q_0}\bigl(\bigl\llVert y^n_w\bigr\rrVert
^4\bigr)\lesssim k_n^{-2},\qquad\sum
_{w=1}^{k_n-1}E_{Q_0}\bigl(\bigl\llVert
y^n_w\bigr\rrVert^4\bigr)\to0,
\\
\sum_{w=1}^{k_n-1}E_{Q_0}
\bigl(y^{n,k}_w y^{n,l}_w\bigr)&=&
\frac{n^{1/2}}{k_n^2}\sum_{w=1}^{k_n-1}\bigl
\llvert(\phi_{2})^n_w\bigr\rrvert
^2\Upsilon^{kl}\bigl(\omega^{(0)}
\bigr)_{\tau^k_{i(\rho_n)^n(\omega^{(0)})+w}}1_{\{\tau^k_{i(\rho
_n)^n(\omega^{(0)})+w}=\tau^l_{i(\rho_n)^n(\omega^{(0)})+w}\}}.
\end{eqnarray*}
Since $\Upsilon$ is cadlag, $(\ref{JPeq16318})$ and the fact that
$n^{1/2}/k_n\to\theta^{-1}$ yield
$\sum_{w=1}^{k_n-1}E_{Q_0}(y^{n,k}_w y^{n,l}_w)
=\theta^{-1}\Upsilon^{kl}_{(\rho)}(\omega^{(0)})D^{kl}_n(\omega^{(0)})+\mathrm{o}_p(1)$,
where we\vspace*{1pt} set $\Upsilon^{kl}_{(\rho)}=\Upsilon^{kl}_{\rho-}$ in case~(1)
and $\Upsilon^{kl}_{(\rho)}=\Upsilon^{kl}_{\rho}$ in case~(2). Since
$\omega^{(0)}\in\Omega^{(0)}$, this implies that
\begin{eqnarray*}
\sum_{w=1}^{k_n-1}E_{Q_0}
\bigl(y^{n,k}_w y^{n,l}_w\bigr)
\to^p\llVert\phi_{2}\rrVert^2
\theta^{-1}\Upsilon^{kl}_{(\rho)}\bigl(
\omega^{(0)}\bigr)\chi^{kl}_{(\rho)}\bigl(
\omega^{(0)}\bigr).
\end{eqnarray*}
Now a standard central limit theorem on row-wise independent triangular
arrays of infinitesimal variables (e.g.,~Theorem 2.2.14 of \cite
{JP2012}) yields $(\ref{noiseCLT})$.

 \textit{Step} 4. In this step, we shall show the following
convergence for $L_n$:
%
\begin{equation}
\label{wienerCLT} L_n\to^d\sqrt{\Phi_{22}\theta
G_{(\rho)}}\zeta,
\end{equation}
where $\zeta$ is a standard $d'$-dimensional normal variable
independent of $\mathcal{F}$. Unlike step~3, here the limiting variable
is mixed normal, so we cannot rely on the standard central limit
theorem used in step~3. Instead, we use the classic mixed normal limit
theorem of \citet{Hall1977}.

Fix $u\in\mathbb{R}^{d'}$ arbitrarily and set $y(u)^n_w=n^{1/4}(\phi
_{1})^n_wu^*W(I_{i(\rho_n)^n+w})$ for each $w=1,\dots, k_n-1$.
Then $y(u)^n_w$ is $\mathcal{G}_{T_{i(\rho_n)^n+w}}$-measurable and
$u^*L_n=\sum_{w=1}^{k_n-1}y(u)^n_w$. Therefore, noting that $G$ and
$G_-$ do not vanish, it suffices to verify the following four
conditions according to \cite{Hall1977} and the Cram\'er--Wold method:
%
\begin{eqnarray}
E \Bigl[\max_{1\leq w\leq k_n-1}\bigl\llvert y(u)^n_w
\bigr\rrvert^2 \Bigr]&\to&0,\label{Halleq1}
\\
\sum_{w=1}^{k_n-1}\bigl\llvert
y(u)^n_w\bigr\rrvert^2-\llVert u\rrVert
^2\llVert\phi_{1}\rrVert^2\theta
G_{\rho_n}&\to^p&0,\label{Halleq2}
\\
 \llVert u\rrVert^2\llVert\phi_{1}\rrVert
^2\theta G_{\rho_n}&\to^p&\llVert u\rrVert
^2\llVert\phi_{1}\rrVert^2\theta
G_{(\rho)},\label{Halleq3}
\\
 \sum_{w=1}^{k_n-1}\bigl\llvert E
\bigl[y(u)^n_w\mid\mathcal{G}_{T_{i(\rho_n)^n+w-1}}\bigr]
\bigr\rrvert&\to^p&0.\label{Halleq4}
\end{eqnarray}
Equation (\ref{Halleq1}) follows from $(\ref{SA4})$ and L\'evy's theorem on
the uniform modulus of continuity of Brownian motion. Next, [SA1] and
{Lemma 2.3 of \cite{Fu2010b}} imply that
\begin{eqnarray*}
\sum_{w=1}^{k_n-1}\bigl\llvert
y(u)^n_w\bigr\rrvert^2&=&\llVert u\rrVert
^2\sqrt{n}\sum_{w=1}^{k_n-1}(
\phi_{1})^n_w\llvert I_{i(\rho_n)^n+w}\rrvert
+\mathrm{o}_p(1) =\frac{\llVert u\rrVert ^2}{\sqrt{n}}\sum_{w=1}^{k_n-1}(
\phi_{1})^n_wG_{T_{i(\rho_n)^n+w-1}}+\mathrm{o}_p(1),
\end{eqnarray*}
hence we obtain $(\ref{Halleq2})$ because $G$ is cadlag. Finally, the
fact that $G$ is cadlag and $(\ref{JPeq16318})$ yield~$(\ref
{Halleq3})$, while we have $E[y(u)^n_w\mid\mathcal{G}_{T_{i(\rho
_n)^n+w-1}}]=0$ because\vspace*{1pt} $W$ is a $d'$-dimensional $(\mathcal
{F}_t)$-Brownian motion independent of $\mathcal{G}$, hence $(\ref
{Halleq4})$ holds true.

 \textit{Step} 5. We denote by $\Psi_n(U)$ and $\Psi(U)$ the
left-hand and right-hand sides of $(\ref{eqlocalCLT})$, respectively.
In this step, we show that it is enough to prove
%
\begin{equation}
\label{eqlocalCLT1} \Psi_n(1)\to^p\Psi(1).
\end{equation}
In fact, assume this, and take an arbitrary bounded variable $U$. We
consider the cadlag version of the bounded martingale $U_t=E(U\mid
{\mathcal{G}^{(0)}_t})$.

First, suppose that we are in case~(1). Set $\underline
{k}_n=k_n-\lfloor n^\beta\rfloor$ and define the $d'$-dimensional
variable $\underline{L}_n=(\underline{L}^j_n)_{1\leq j\leq d'}$ and the
$d$-dimensional variable $\underline{L}'_n=(\underline{L}^{\prime
}_n)_{1\leq k\leq d}$ by
\begin{eqnarray*}
\underline{L}^j_n=n^{1/4}\sum
_{w=1}^{\underline{k}_n}(\phi_1)^n_wW^j(I_{i(\rho)^n+w}),
\qquad\underline{L}'^k_n=\frac{n^{1/4}}{k_n}
\sum_{w=1}^{\underline{k}_n}(\phi_2)^n_w
\epsilon^k_{\tau^k_{i(\rho)^n+w}}.
\end{eqnarray*}
Then, since $E[\llVert L_n-\underline{L}_n\rrVert ^2]\lesssim\sqrt
{n}n^\beta\bar{r}_n$ and $E[\llVert L'_n-\underline{L}'_n\rrVert
^2]\lesssim\sqrt{n}k_n^{-2}n^\beta$, by the boundedness of $U$ and the
Lipschitz continuity of $f$ it holds that $\Psi_n(U)-\underline{\Psi
}_n(U)\to^p0$ and $\Psi_n(1)-\underline{\Psi}_n(1)\to^p0$, where
$\underline{\Psi}_n(U)=E[U f(\underline{L}_n,\underline{L}'_n)\mid
\mathcal{G}_{\rho_n}]$. In particular, to prove $(\ref{eqlocalCLT1})$
it is enough to show that $\underline{\Psi}_n(U)\to^p\Psi(U)$. Now,
since both $G_{\rho-}$ and $\widetilde{\upsilon}_{\rho-}$ are {$\mathcal
{G}^{(0)}_{\rho-}$}-measurable, we have $\Psi(U)=U_{\rho-}\Psi(1)$
because $U_{\rho-}=E[U\mid{\mathcal{G}^{(0)}_{\rho-}}]$. Also,
$f(\underline{L}_n,\underline{L}'_n)$ in restriction to the set $\Omega
_n=\{\rho>T_{i(\rho_n)^n+\underline{k}_n}\}$ is $\mathcal{G}_{\rho
-}$-measurable, so $\underline{\Psi}_n(U)=\underline{\Psi}_n(U_{\rho
-})$ on $\Omega_n$. We also obviously have $\underline{\Psi}_n(U_{\rho
_n})=\underline{\Psi}_n(1)U_{\rho_n}\to^p\Psi(1)U_{\rho-}$ by $(\ref
{eqlocalCLT1})$, $\Psi_n(1)-\underline{\Psi}_n(1)\to^p0$ and $U_{\rho
_n}\to U_{\rho-}$, while $P(\Omega_n)\to1$ by assumption. Now, since
$E[\llvert \underline{\Psi}_n(U_{\rho_n})-\underline{\Psi}_n(U_{\rho
-})\rrvert ]\leq\llVert f\rrVert _\infty E[\llvert U_{\rho
_n}-U_{\rho-}\rrvert ]\to0$ by the boundedness of $f$ and $U$, $U_{\rho
_n}\to U_{\rho-}$ on $\Omega_n$ and the fact that $P(\Omega_n)\to1$, we
obtain the desired result.

Next, suppose that we are in case~(2). Then $\Psi(U)=U_\rho\Psi(1)$
because $\Psi(1)$ is {$\mathcal{G}^{(0)}_\rho$}-measurable, and also
$\Psi_n(U_\rho)=U_\rho\Psi_n(1)$ because $\rho_n\geq\rho$. Moreover,
setting $\rho'_n=\rho_n+k_n\bar{r}_n$, $\Psi_n(1)$ is $\mathcal{G}_{\rho
'_n}$-measurable due to $(\ref{SA4})$, so $\Psi_n(U)=U_{\rho'_n}\Psi
_n(1)$. Since $\rho'_n\to\rho$ and $\rho'_n>\rho$, we have $U_{\rho
'_n}\to U_\rho$, and the same arguments as above shows that $\Psi
_n(U_{\rho'_n})-\Psi_n(U_{\rho})\to^p0$, thus the desired result is obtained.

 \textit{Step} 6. Now we finish the proof by proving the
convergence $(\ref{eqlocalCLT1})$. First, for each $\omega^{(0)}\in
\Omega_0$ define the function $h^n_{\omega^{(0)}}$ on $\mathbb{R}^d$ by
$h^n_{\omega^{(0)}}(y)=f(L_n(\omega^{(0)}),y)$. Then, noting that $f$
is bounded and Lipschitz continuous, Lemma~\ref{lemSk} and $(\ref
{noiseCLT})$ imply that
\begin{eqnarray*}
\int h^n_{\omega^{(0)}}\bigl(L'_n\bigl(
\omega^{(1)}\bigr)\bigr)Q\bigl(\omega^{(0)},\mathrm{d}
\omega^{(1)}\bigr) -\int h^n_{\omega^{(0)}}\bigl(\llVert
\phi_2\rrVert\sqrt{\theta^{-1}}\widetilde{
\upsilon}_{(\rho)}\bigl(\omega^{(0)}\bigr)y\bigr)
\mathfrak{N}_d(\mathrm{d}y) \to0.
\end{eqnarray*}
Since $f$ is bounded and $P(\Omega_0)=1$, this convergence and the
bounded convergence theorem yield
\[
\Psi_n(1)- E \biggl[\int f \bigl(L_n,\llVert
\phi_2\rrVert\sqrt{\theta^{-1}}\widetilde{
\upsilon}_{(\rho)}y \bigr)\mathfrak{N}_d(\mathrm{d}y)\Big|
\mathcal{G}_{\rho_n} \biggr] \to^p0.
\]
Next, since $f$ is Lipschitz continuous and $\widetilde{\upsilon}$ is
cadlag and bounded, by $(\ref{JPeq16318})$ we obtain
%
\begin{equation}
\label{resultlocalCLT1} \Psi_n(1)- E \biggl[\int f \bigl(L_n,
\llVert\phi_2\rrVert\sqrt{\theta^{-1}}\widetilde{
\upsilon}_{\rho_n}y \bigr)\mathfrak{N}_d(\mathrm{d}y)\Big|
\mathcal{G}_{\rho_n} \biggr] \to^p0.
\end{equation}
Now, noting that $W$ is a standard $d'$-dimensional Brownian motion
with respect to {$(\mathcal{G}_t)$}, by the strong Markov property of a
Brownian motion $(W_{\rho_n+t}-W_{\rho_n})_{t\geq0}$ is independent of
$\mathcal{G}_{\rho_n}$, hence we have
%
\begin{equation}
\label{resultlocalCLT2} \hspace*{-15pt}
E \biggl[\int f \bigl(L_n,\llVert\phi_2
\rrVert\sqrt{\theta^{-1}}\widetilde{\upsilon}_{\rho_n}y \bigr)
\mathfrak{N}_d(\mathrm{d}y)\Big|\mathcal{G}_{\rho_n} \biggr] =
\int f \bigl(x,\llVert\phi_2\rrVert\sqrt{\theta^{-1}}
\widetilde{\upsilon}_{\rho_n}y \bigr)\mathbb{P}^n(\mathrm{d}x)
\mathfrak{N}_d(\mathrm{d}y),
\end{equation}
where $\mathbb{P}^n$ is the law of $L_n$ under $P^{(0)}$. Then, again
using the Lipschitz continuity of $f$ and the cadlag property of
$\widetilde{\upsilon}$ as well as $(\ref{JPeq16318})$, we obtain
\begin{eqnarray*}
\int f \bigl(x,\llVert\phi_2\rrVert\sqrt{\theta^{-1}}
\widetilde{\upsilon}_{\rho_n}y \bigr)\mathbb{P}^n(\mathrm{d}x)
\mathfrak{N}_d(\mathrm{d}y) -\int f \bigl(x,\llVert
\phi_2\rrVert\sqrt{\theta^{-1}}\widetilde{
\upsilon}_{(\rho)}y \bigr)\mathbb{P}^n(\mathrm{d}x)
\mathfrak{N}_d(\mathrm{d}y) \to^p0,
\end{eqnarray*}
hence $(\ref{wienerCLT})$ yields
%
\begin{eqnarray}\label{resultlocalCLT3}
&& \int f \bigl(x,\llVert\phi_2\rrVert\sqrt{
\theta^{-1}}\widetilde{\upsilon}_{\rho_n}y \bigr)
\mathbb{P}^n(\mathrm{d}x)\mathfrak{N}_d(\mathrm{d}y)
\nonumber\\[-8pt]\\[-8pt]\nonumber
&&\quad \to^p\int f \bigl(\llVert\phi_1\rrVert\sqrt{\theta
G_{(\rho)}}x,\llVert\phi_2\rrVert\sqrt{
\theta^{-1}}\widetilde{\upsilon}_{(\rho)}y \bigr)
\mathfrak{N}_{d'}(\mathrm{d}x)\mathfrak{N}_d(\mathrm{d}y).
\end{eqnarray}
Equations (\ref{resultlocalCLT1})--(\ref{resultlocalCLT3}) imply that $(\ref
{eqlocalCLT1})$ holds true, and thus we complete the proof.
\end{pf}


\subsubsection{A joint convergence result and the proof of Proposition
\texorpdfstring{\protect\ref{mainCLT}}{6.2}}\label{secjointCLT}

In this subsection, we prove a joint convergence result for the pair
$(n^{1/4}\mathbf{L}[M]^n,(z^n_{r-},z'^n_{r-},\break z^n_{r+}, z'^n_{r+})_{r\geq
1})$ and complete the proof of Proposition~\ref{mainCLT}.

For the proof, we use some elementary results on the Skorokhod
topology. For any $k\in\mathbb{N}$, denote by $\mathbb{D}^k$
(resp.,~$\mathbb{D}^{k\times k}$) the space of $\mathbb{R}^k$-valued
(resp.,~$\mathbb{R}^k\otimes\mathbb{R}^k$-valued) cadlag functions on
$\mathbb{R}_+$ equipped with the Skorokhod topology. For any $x\in
\mathbb{D}^k$ and any $t\in\mathbb{R}_+$, we define the function $x^t$
by $x^t(s)=x(s\wedge t)$ for $s\in\mathbb{R}_+$. We evidently have
$x^t\in\mathbb{D}^k$. On the other hand, for any $S\geq0$ we define the
function $\pi_S$ from $\mathbb{D}^{d\times d}$ into itself by $\pi
_S(x)(t)=x(t-S)1_{\{t\geq S\}}$.

The following two lemmas can be shown using basic properties of the
Skorokhod topology, so we omit the proofs.

%
\begin{lem}\label{contstop}
The map $\mathbb{R}_+\times\mathbb{D}^k\ni(t,x)\mapsto x^t\in\mathbb
{D}^k$ is continuous at every point $(t,x)\in\mathbb{R}_+\times\mathbb
{D}^k$ such that $x$ is continuous at $t$, where the space $\mathbb
{R}_+\times\mathbb{D}^k$ is equipped with the product topology.
\end{lem}

%
\begin{lem}\label{shiftcont}
$\pi_S$ is a continuous function of $\mathbb{D}^{d\times d}$ into itself.
\end{lem}

Now we are ready to prove the following joint convergence result:

\begin{prop}\label{jointCLT}
Suppose that~\textup{[SA1]}--\textup{[{SA3}]} and~\textup{[A4]} are
satisfied. Then
\[
\bigl(n^{1/4}\mathbf{L}[M]^n,\bigl(z^n_{r-},z'^n_{r-},z^n_{r+},z'^n_{r+}
\bigr)_{r\geq1}\bigr) \to^{d_s}\bigl(\mathcal{W},
\bigl(z_{r-},z'_{r-},z_{r+},z'_{r+}
\bigr)_{r\geq1}\bigr)
\]
as $n\to\infty$ for the product topology on the space $\mathbb
{D}^{d\times d}\times(\mathbb{R}^{2(d+d')})^{\mathbb{N}}$, where
\[
\cases{ z_{r-}=\psi_2^{-1}\sqrt{
\Phi_{22}\theta G_{S_r-}}\Psi_{r-}, &\quad
$z_{r+}=\psi_2^{-1}\sqrt{\Phi_{22}
\theta G_{S_r}}\Psi_{r+}$,
\vspace*{3pt}\cr
z'_{r-}=
\psi_2^{-1}\sqrt{\Phi_{12}\theta^{-1}}
\widetilde{\upsilon}_{S_r-}\Psi'_{r-}, &\quad
$z'_{r+}=\psi_2^{-1}\sqrt{
\Phi_{12}\theta^{-1}}\widetilde{\upsilon}_{S_r}
\Psi'_{r+}$.}
\]
\end{prop}

\begin{pf}
\textit{Step} 1. It suffices to prove
%
\begin{equation}
\label{jointCLTaim} \bigl(n^{1/4}\mathbf{L}[M]^n,
\bigl(z^n_{r-},z'^n_{r-},z^n_{r+},z'^n_{r+}
\bigr)_{r\in\mathcal{R}}\bigr) \to^{d_s}\bigl(\mathcal{W},
\bigl(z_{r-},z'_{r-},z_{r+},z'_{r+}
\bigr)_{r\in\mathcal{R}}\bigr)
\end{equation}
in $\mathbb{D}^{d\times d}\times\mathbb{R}^{2(d'+d)\#\mathcal{R}}$ for
any finite subset $\mathcal{R}$ of $\mathbb{N}$, and we prove this by
induction on the number $\#\mathcal{R}$ of the elements in the set
$\mathcal{R}$. First, $(\ref{jointCLTaim})$ holds true when $\#\mathcal
{R}=0$ due to Proposition~\ref{contCLT}. Next, let $J\in\mathbb{N}$ and
assume that $(\ref{jointCLTaim})$ holds true when $\#\mathcal{R}=J-1$.
Then, we need to prove~$(\ref{jointCLTaim})$ for the case that $\#
\mathcal{R}=J$. We write $\mathcal{R}=\{r_1,\dots,r_J\}$ with
$S_{r_1}<\cdots<S_{r_J}$.


 \textit{Step} 2. Before stating the detailed proof, we briefly
explain the intuition behind the proof. The basic idea is the same as
in the proof of Theorem 4.3.1 from \cite{JP2012}. Namely, for each
$\beta>0$ we set $S^{\beta-}=(S_{r_J}-\beta)_+$ and $S^{\beta+}=S^{\beta
+}$, and divide $n^{1/4}\mathbf{L}[M]^n$ into the summands containing
the data observed in the interval $[S^{\beta-},S^{\beta+}]$ and the
remaining ones. Then we prove the negligibility of the former part (as
$\beta\to0$) and the joint limit theorem of the latter part and
$(z^n_{r-},z'^n_{r-},z^n_{r+},z'^n_{r+})_{r\in\mathcal{R}}$. More
formally, we set
$\widehat{\mathbf{L}}(\beta)^n=n^{1/4} (\mathbf{L}[M]^n )^{S^{\beta
+}}-n^{1/4} (\mathbf{L}[M]^n )^{S^{\beta-}}$
and
$\widehat{\mathcal{W}}(\beta)=\mathcal{W}^{S^{\beta+}}-\mathcal
{W}^{S^{\beta-}}$,
and show that
%
\begin{eqnarray}\label{jointaim1}
\limsup_{\beta\to0}\limsup_{n\to\infty}P
\Bigl(\sup_{0\leq t\leq T}\bigl\llVert\widehat{\mathbf{L}}(
\beta)^n_t\bigr\rrVert>\eta\Bigr)&=&0,
\nonumber\\[-8pt]\\[-8pt]\nonumber
\limsup
_{\beta\to0}P \Bigl(\sup_{0\leq t\leq T}\bigl\llVert
\widehat{\mathcal{W}}(\beta)_t\bigr\rrVert>\eta\Bigr)&=&0
\end{eqnarray}
for any $T,\eta>0$ and that
%
\begin{eqnarray}\label{jointaim11}
&& \bigl(n^{1/4}\mathbf{L}[M]^n-\widehat{
\mathbf{L}}(\beta)^n,\bigl(z^n_{r-},z'^n_{r-},z^n_{r+},z'^n_{r+}
\bigr)_{r\in\mathcal{R}}\bigr)
\nonumber\\[-8pt]\\[-8pt]\nonumber
&&\quad \to^{d_s}\bigl(\mathcal{W}-\widehat{
\mathcal{W}}(\beta),\bigl(z_{r-},z'_{r-},z_{r+},z'_{r+}
\bigr)_{r\in\mathcal{R}}\bigr)
\end{eqnarray}
in $\mathbb{D}^{d\times d}\times\mathbb{R}^{2(d'+d)J}$ as $n\to\infty$
for any fixed $\beta>0$. Then, Proposition 2.2.4 of \cite{JP2012}
yields~(\ref{jointCLTaim}).

However, to prove $(\ref{jointaim11})$ we need a different approach
from the one of \cite{JP2012} because we cannot argue conditionally on
the increments of $W$ consisting of the observations in $[S^{\beta
-},S^{\beta+}]$ as \cite{JP2012} do, which is due to the time
endogeneity. For this reason we further decompose $n^{1/4}\mathbf
{L}[M]^n-\widehat{\mathbf{L}}(\beta)^n$ as
$n^{1/4}\mathbf{L}[M]^n-\widehat{\mathbf{L}}(\beta)^n=\check{\mathbf
{L}}^n+\widetilde{\mathbf{L}}(\beta)^n$,
where $\check{\mathbf{L}}^n=(n^{1/4}\mathbf{L}[M]^n)^{S^{\beta-}}$.
Roughly speaking, $\check{\mathbf{L}}^n$ consists of the data observed
before $S^{\beta-}$, while $\widetilde{\mathbf{L}}(\beta)^n$ consists
of those observed after $S^{\beta+}$. By Proposition VI-1.23 of \cite
{JS} and the continuous mapping theorem $(\ref{jointaim11})$ follows
once we show that
%
\begin{equation}
\label{jointaim12} \bigl(\check{\mathbf{L}}^n,\widetilde{\mathbf{L}}(
\beta)^n,\bigl(z^n_{r-},z'^n_{r-},z^n_{r+},z'^n_{r+}
\bigr)_{r\in\mathcal{R}}\bigr) \to^{d_s}\bigl(\mathcal{W}^{S_{r_J}},
\widetilde{\mathcal{W}}(\beta),\bigl(z_{r-},z'_{r-},z_{r+},z'_{r+}
\bigr)_{r\in\mathcal{R}}\bigr)
\end{equation}
in $\mathbb{D}^{d\times d}\times\mathbb{D}^{d\times d}\times\mathbb
{R}^{2(d'+d)J}$ as $n\to\infty$. The strategy of the proof of $(\ref
{jointaim12})$ is, roughly speaking, as follows. We first prove a
stable limit theorem for $\widetilde{\mathbf{L}}(\beta)^n$
conditionally on $\mathcal{F}_{S^{\beta+}}$ (this will be done in step~5;
the assumption~\textup{[A4]}(i) is necessary for this part). Then we
obtain a joint stable limit theorem for $\widetilde{\mathbf{L}}(\beta
)^n$ and $(z^n_{r_J-},z'^n_{r_J-},z^n_{r_J+},z'^n_{r_J+})$
conditionally on $\mathcal{G}_{S^\dagger_{r_J}}$ by virtue of Lemma~\ref
{JPlem1637} (step~7). Finally, from the assumption of the induction we
will obtain the desired result (step~8).


 \textit{Step} 3. We begin with proving $(\ref{jointaim1})$.
The second equation immediately follows from the continuity of the
process $\mathcal{W}$. On the other hand, for any $\beta>0$ we have
$(S^{\beta-},S^{\beta+},n^{1/4}\mathbf{L}[M]^n)\to^{d_s}(S^{\beta
-},S^{\beta+},\mathcal{W})$ as $n\to\infty$ in $\mathbb{R}_+\times
\mathbb{R}_+\times\mathbb{D}^{d\times d}$ by Proposition~\ref{contCLT},
hence Lemma~\ref{contstop} and the continuous mapping theorem imply
that $(n^{1/4} (\mathbf{L}[M]^n )^{S^{\beta+}},n^{1/4} (\mathbf{L}[M]^n
)^{S^{\beta-}})\to^{d_s}(\mathcal{W}^{S^{\beta+}},\mathcal{W}^{S^{\beta
-}})$ as $n\to\infty$ in $\mathbb{D}^{d\times d}\times\mathbb
{D}^{d\times d}$. Therefore, Propositions VI-1.23 and VI-2.4 of \cite
{JS} as well as the continuous mapping theorem yield $\sup_{0\leq t\leq
T}\llVert \widehat{\mathbf{L}}(\beta)^n_t\rrVert \to^{d_s}\sup_{0\leq
t\leq T}\llVert \mathcal{W}(\beta)_t\rrVert $ as $n\to\infty$. In
particular, we have
\[
\limsup_{n\to\infty}P \Bigl(\sup_{0\leq t\leq T}\bigl
\llVert\widehat{\mathbf{L}}(\beta)^n_t\bigr\rrVert>\eta
\Bigr)\leq P \Bigl(\sup_{0\leq t\leq T}\bigl\llVert\mathcal{W}(
\beta)_t\bigr\rrVert\geq\eta\Bigr),
\]
hence we also obtain the first equation $(\ref{jointaim1})$.


 \textit{Step} 4. Now we start the proof of $(\ref
{jointaim12})$. First, due to the property of the product topology it
suffices to prove the following convergence:
%
\begin{equation}
\label{jointaim2} E \Biggl[\zeta f_1\bigl(\check{\mathbf{L}}^n
\bigr)f_2\bigl(\widetilde{\mathbf{L}}(\beta)^n\bigr)\prod
_{j=1}^JY^n_{j-}Y^n_{j+}
\Biggr] \to E \Biggl[\zeta f_1\bigl(\mathcal{W}^{S_{r_J}}
\bigr)f_2\bigl(\widetilde{\mathcal{W}}(\beta)\bigr)\prod
_{j=1}^JY_{j-}Y_{j+} \Biggr]
\end{equation}
as $n\to\infty$, where $\zeta$ is any bounded {$\mathcal
{F}^{(0)}$}-measurable variable, $f_1$ and $f_2$ are bounded Lipschitz
functions on $\mathbb{D}^{d\times d}$, and $Y^n_{j\pm}=F_{j\pm
}(z^n_{r_j\pm},z'^n_{r_j\pm})$ and $Y_{j\pm}=F_{j\pm}(z_{r_j\pm
},z'_{r_j\pm})$ with $F_{j-}$ and $F_{j+}$ being bounded Lipschitz
functions on $\mathbb{R}^{d'+d}$ for every $j=1,\dots,J$.


 \textit{Step} 5. We begin with proving
%
\begin{equation}
\label{jointaim3} E \Biggl[\zeta f_1\bigl(\check{\mathbf{L}}^n
\bigr)f_2\bigl(\widetilde{\mathbf{L}}(\beta)^n\bigr)\prod
_{j=1}^JY^n_{j-}Y^n_{j+}
\Biggr] -E \Biggl[\zeta f_1\bigl(\check{\mathbf{L}}^n
\bigr)f_2\bigl(\widetilde{\mathcal{W}}(\beta)\bigr)\prod
_{j=1}^JY^n_{j-}Y^n_{j+}
\Biggr] \to0
\end{equation}
as $n\to\infty$. First, we introduce some notation. For any
$d$-dimensional processes $U,V$, any $u,v\in\{g,g'\}$ and any
$k,l=1,\dots,d$, we define the process $\widetilde{\mathbb
{L}}^{(k,l)}_{u,v}(U,V)^n$ in the same way as that of $\mathbb
{L}^{(k,l)}_{u,v}(U,V)^n$ with replacing $(T_p)_{p\geq0}$ by
$(\widetilde{T}_p)_{p\geq0}:=(T_{i(S^{\beta+})^n+1+p}-S^{\beta+})_{p\geq
0}$. Also, for any process $V$ we define the process $\mathring{V}$ by
$\mathring{V}_t=V_{S^{\beta+}+t}-V_{S^{\beta+}}$, and define the
$\mathbb{R}^d\otimes\mathbb{R}^d$-valued process $\widetilde{\mathbf
{L}}^n$ by
\begin{eqnarray*}
\widetilde{\mathbf{L}}^{n,kl}=n^{1/4} \bigl\{\widetilde{
\mathbb{L}}^{(k,l)}_{g,g}(\mathring{M},\mathring{M})^n+
\widetilde{\mathbb{L}}^{(k,l)}_{g,g'}(\mathring{M},\mathring{
\mathfrak{E}})^n +\widetilde{\mathbb{L}}^{(l,k)}_{g,g'}(
\mathring{M},\mathring{\mathfrak{E}})^n +\widetilde{
\mathbb{L}}^{(k,l)}_{g',g'}(\mathring{\mathfrak{E}},\mathring{
\mathfrak{E}})^n \bigr\}.
\end{eqnarray*}
Then it can easily be seen that $\widetilde{\mathbf{L}}(\beta)^n-\pi
_{S^{\beta+}}(\widetilde{\mathbf{L}}^n)\xrightarrow{ucp}0$ as $n\to
\infty$. Therefore, by the Lipschitz continuity of $f_2$ as well as the
boundedness of $\zeta$, $f_1$ and $Y^n_{j\pm}$ we have
%
\begin{equation}
\label{jointshift} E \Biggl[\zeta f_1\bigl(\check{\mathbf{L}}^n
\bigr)f_2\bigl(\widetilde{\mathbf{L}}(\beta)^n\bigr)\prod
_{j=1}^JY^n_{j-}Y^n_{j+}
\Biggr] -E \Biggl[\zeta f_1\bigl(\check{\mathbf{L}}^n
\bigr)f_2\bigl(\pi_{S^{\beta+}}\bigl(\widetilde{
\mathbf{L}}^n\bigr)\bigr)\prod_{j=1}^JY^n_{j-}Y^n_{j+}
\Biggr] \to0
\end{equation}
as $n\to\infty$.

Now we consider a regular conditional probability $p^{(0)}(\omega
^{(0)},\cdot)$ of $P^{(0)}$ given $\mathcal{F}^{(0)}_{S^{\beta+}}$.
Such one exists because of the assumption~\textup{[A4]}(i). We also
consider a filtration $(\mathring{\mathcal{F}}^{(0)}_t)_{t\geq0}$ of
$\mathcal{F}^{(0)}$ defined by $\mathring{\mathcal{F}}^{(0)}_t=\mathcal
{F}^{(0)}_{S^{\beta+}+t}$, and for each $\omega_0\in\Omega^{(0)}$ we
introduce a stochastic basis $\mathcal{B}^{(0)}_{\omega_0}:=(\Omega
^{(0)},\mathcal{F}^{(0)},(\mathring{\mathcal
{F}}^{(0)}_{t}),p^{(0)}(\omega_0,\cdot))$. For each $t\in\mathbb{R}_+$
we also introduce a transition probability $\mathring{Q}_t(\omega
^{(0)},\mathrm{d}u)$ from $(\Omega^{(0)},\mathring{\mathcal
{F}}^{(0)}_{t})$ into $\mathbb{R}^d$ by setting $\mathring{Q}_t(\omega
^{(0)},A)=Q_{S^{\beta+}(\omega^{(0)})+t}(\omega^{(0)},A)$ for each
Borel set $A$ of $\mathbb{R}^d$. {Note that the process $(\mathring
{Q}_t(\cdot,A))_{t\geq0}$ is $(\mathring{\mathcal
{F}}_t^{(0)})$-progressively measurable because of~\textup{[A4]}(ii) and
Theorem IV-57 of \cite{Meyer1966}}. Now, by replacing $\mathcal
{B}^{(0)}$, $Q_t(\omega^{(0)},\mathrm{d}u)$ and $(\mathcal{T}^n_i)$
with $\mathcal{B}^{(0)}_{\omega_0}$, $\mathring{Q}_t(\omega
^{(0)},\mathrm{d}u)$ and the increasing reordering of $\{\widetilde{\tau
}^k_p:=\tau^k_{i(S^{\beta+})^n+1+p}-S^{\beta+}\dvt k=1,\dots,d$ and
$p\geq0\}$, respectively, we introduce the new stochastic basis
$\mathcal{B}_{\omega_0}=(\Omega,\mathcal{F},(\mathring{\mathcal
{F}}_t),P_{\omega_0})$ instead of~$\mathcal{B}$.

By the strong Markov property of a Brownian motion $\mathring{W}$ is a
$d'$-dimensional standard Brownian motion on $\mathcal{B}^{(0)}_{\omega
_0}$. Moreover, defining the random measure $\mathring{\mu}$ by
$\mathring{\mu}((0,t]\times A)=\mu((S^{\beta+},S^{\beta+}+t]\times A)$,
the strong Markov property of a Poisson random measure implies that
$\mathring{\mu}$ is a Poisson random measure on $\mathcal
{B}^{(0)}_{\omega_0}$ with compensator $\nu$. We also have
\begin{eqnarray*}
\mathring{M}_t&=&\int_0^t
\sigma_{S^{\beta+}+s}\,\mathrm{d}\mathring{W}_s,
\\
\sigma_{S^{\beta+}+t}&=&\sigma_{S^{\beta+}}+\int_0^t
\widetilde{b}_{S^{\beta+}+s}\,\mathrm{d}s+\int_0^t
\widetilde{\sigma}_{S^{\beta+}+s}\,\mathrm{d}\mathring{W}_s +(
\mathring{\widetilde{\delta}}1_{\{\llvert \mathring{\widetilde{\delta
}}\rrvert \leq1\}})\star(\mathring{\mu}-
\nu)_t+(\mathring{\widetilde{\delta}}1_{\{\llvert \mathring{\widetilde
{\delta}}\rrvert >1\}})\star\mathring{
\mu}_t,
\\
G_{S^{\beta+}+t}&=&G_{S^{\beta+}}+\int_0^t
\widehat{b}_{S^{\beta+}+s}\,\mathrm{d}s+\int_0^t
\widehat{\sigma}_{S^{\beta+}+s}\,\mathrm{d}\mathring{W}_s +(\mathring{
\widehat{\delta}}1_{\{\llvert \mathring{\widehat{\delta}}\rrvert \leq
1\}})\star(\mathring{\mu}-\nu)_t+(
\mathring{\widehat{\delta}}1_{\{\llvert \mathring{\widehat{\delta
}}\rrvert >1\}})\star\mathring{\mu}_t,
\end{eqnarray*}
where for a function $\eta$ on $\Omega^{(0)}\times\mathbb{R}_+\times E$
the function $\mathring{\eta}$ on $\Omega^{(0)}\times\mathbb{R}_+\times
E$ is defined by $\mathring{\eta}(\omega^{(0)},t,z)=\eta(\omega
^{(0)},S^{\beta+}(\omega^{(0)})+t,z)$. Therefore, noting that for any
$\mathcal{F}$-measurable variable $x$ and any sub $\sigma$-filed
$\mathcal{H}$ of $\mathcal{F}$ we have $E_{P_{\omega_0}}[x\mid\mathcal
{H}]=E[x\mid\mathcal{H}]$ as long as $\mathcal{F}^{(0)}_{S^{\beta
+}}\subset\mathcal{H}$, it can easily been shown that the conditions
\textup{[SA1]}--\textup{[{SA3}]} are satisfied with replacing {$\mathcal
{B}$, $X$, $(T_p)$, $(\tau^k_p)$, $G_t$ and $\chi_t$ by $\mathcal
{B}_{\omega_0}$, $\mathring{M}$, $(\widetilde{T}_p)$, $(\widetilde{\tau
}^k_{p})$, $G_{S^{\beta+}+t}$ and $\chi_{S^{\beta+}+t}$}, respectively.

Consequently, Proposition~\ref{contCLT} and Lemma~\ref{shiftcont} as
well as the continuous mapping theorem yield
\begin{eqnarray*}
E_{P_{\omega_0}}\bigl[\zeta f_2\bigl(\pi_{S_{r_J}(\omega_0)+\beta}\bigl(
\widetilde{\mathbf{L}}^n\bigr)\bigr)\bigr] \to E_{p^{(0)}(\omega_0,\cdot
)}\bigl[
\zeta f_2\bigl(\pi_{S_{r_J}(\omega_0)+\beta}(\mathring{\mathcal
{W}})\bigr)
\bigr].
\end{eqnarray*}
Therefore, noting that $E_{P_{\omega_0}}[\cdot]=E[\cdot\mid\mathcal
{F}_{S^{\beta+}}](\omega_0)$ and $E_{p^{(0)}(\omega_0,\cdot)}[\zeta
f_2(\pi_{S_{r_J}(\omega_0)+\beta}(\mathring{\mathcal{W}}))]= E[\zeta
f_2(\pi_{S_{r_J}(\omega_0)+\beta}(\mathring{\mathcal{W}}))\mid\mathcal
{F}_{S^{\beta+}}](\omega_0)$ for almost all $\omega_0$ (with respect to
$P^{(0)}$) and that $f_1(\check{\mathbf{L}}^n)\prod
_{j=1}^JY^n_{j-}Y^n_{j+}$ is bounded and $\mathcal{F}_{S^{\beta
+}}$-measurable (for sufficiently large $n$; note that $S^{\beta-}$ is
an $(\mathcal{F}_{S^{\beta+}+t})$-stopping time), the bounded
convergence theorem implies that
\[
E \Biggl[\zeta f_1\bigl(\check{\mathbf{L}}^n
\bigr)f_2\bigl(\pi_{S^{\beta+}}\bigl(\widetilde{
\mathbf{L}}^n\bigr)\bigr)\prod_{j=1}^JY^n_{j-}Y^n_{j+}
\Biggr] -E \Biggl[\zeta f_1\bigl(\check{\mathbf{L}}^n
\bigr)f_2\bigl(\pi_{S^{\beta+}}(\mathring{\mathcal{W}})\bigr)\prod
_{j=1}^JY^n_{j-}Y^n_{j+}
\Biggr] \to0.
\]
Since $\pi_{S^{\beta+}}(\mathring{\mathcal{W}})=\widetilde{\mathcal
{W}}(\beta)$, by $(\ref{jointshift})$ and the above convergence we
obtain $(\ref{jointaim3})$.

 \textit{Step} 6. In this step, we prove
%
\begin{equation}
\label{jointaim4} E \bigl[\zeta'Y^n_{J-}Y^n_{J+}
\mid\mathcal{G}_{S^\dagger_{r_J}} \bigr]\to^pE \bigl[
\zeta'Y_{J-}Y_{J+}\mid{\mathcal{G}^{(0)}_{S_{r_J-}}}
\bigr],
\end{equation}
where $\zeta'=E[\zeta f_2(\widetilde{\mathcal{W}}(\beta))\mid{\mathcal
{F}^{(0)}}]$. Fix a constant $\alpha'$ such that $1-\xi<\alpha'<(\xi
-\kappa-\frac{1}{2})\wedge(\xi-\frac{3}{4})\wedge\varpi$. Then, set
$\underline{k}_n=k_n-\lfloor n^{1/2-\alpha'}\rfloor$ and define the
$d'$-dimensional variable $\underline{z}^{n}_{r_J-}=(\underline
{z}^{n,j}_{r_J-})_{1\leq j\leq d'}$ and the $d$-dimensional variable
$\underline{z}'^{n}_{r_J-}=(\underline{z}'^{n,k}_{r_J-})_{1\leq k\leq
d}$ by
\begin{eqnarray*}
\underline{z}^{n,j}_{r_J-}=n^{1/4}\sum
_{w=1}^{\underline{k}_n}(\widetilde{\phi_{g,g}})^n_wW^j(I_{i(S^\dagger_r)^n+w}),
\qquad\underline{z}'^{n,k}_{r_J-}=-
\frac{n^{1/4}}{k_n}\sum_{w=1}^{\underline{k}_n}(
\widetilde{\phi_{g',g}})^n_w
\epsilon^k_{\tau^k_{i(S^\dagger_r)^n+w}},
\end{eqnarray*}
and put $\underline{Y}^n_{J-}=F_{J-}(\underline{z}^{n}_{r_J-},\underline
{z}'^{n}_{r_J-})$. Since $E[\llVert z^{n}_{r_J-}-\underline
{z}^{n}_{r_J-}\rrVert ^2]\lesssim n^{1-\alpha'-\xi}$ and $E[\llVert
z'^{n}_{r_J-}-\underline{z}'^{n}_{r_J-}\rrVert ^2]\lesssim n^{-\alpha
'}$ by~\textup{[SA1]}--\textup{[{SA3}]} and the optional sampling theorem,
we have $Y^n_{J-}-\underline{Y}^n_{J-}\to^p0$ due to the Lipschitz
continuity of $F_{J-}$. Therefore, by virtue of the boundedness of
$\zeta'$ and $F_{J+}$, for the proof of $(\ref{jointaim4})$ it is
enough to prove
%
\begin{equation}
\label{jointCLTaim3} E \bigl[\zeta'\underline{Y}^n_{J-}Y^n_{J+}
\bigl\llvert\mathcal{G}_{S^\dagger_{r_J}} \bigr]\to^pE \bigl[
\zeta'Y_{J-}Y_{J+}\bigr\rrvert{
\mathcal{G}^{(0)}_{S_{r_J-}}} \bigr].
\end{equation}
Now, Lemma~\ref{JPlem1637} yields $E [\zeta'Y^n_{J+}\llvert \mathcal
{G}_{S_{r_J}} ]\to^pE [\zeta'Y_{J+}\rrvert {\mathcal
{G}^{(0)}_{S_{r_J}}} ]$. Moreover, setting $\Omega_n= \{T_{i(S^\dagger
_{r_J})+\bar{k}_n}<S_{r_J}\}$, we have
\[
E [\zeta'\underline{Y}^n_{J-}Y^n_{J+}1_{\Omega_n}\mid\mathcal
{G}_{S^\dagger_{r_J}} ]=E [\underline{Y}^n_{J-}E [\zeta'Y^n_{J+}\mid
\mathcal{G}_{S_{r_J}} ]1_{\Omega_n}\mid\mathcal{G}_{S^\dagger_{r_J}} ].
\]
Since $\lim_nP(\Omega_n)=1$ by Lemma~\ref{appstoptime}, the boundedness
of $\zeta'$ and $F_{J\pm}$ and the bounded convergence theorem imply that
%
\begin{equation}
\label{limitYplus} E \bigl[\zeta'\underline{Y}^n_{J-}Y^n_{J+}
\mid\mathcal{G}_{S^\dagger_{r_J}} \bigr]-E \bigl[\zeta''
\underline{Y}^n_{J-}\mid\mathcal{G}_{S^\dagger_{r_J}} \bigr]
\to^p0,
\end{equation}
where $\zeta''=E [\zeta'Y_{r_J+}\mid{\mathcal{G}^{(0)}_{S_{r_J}}} ]$.
On the other hand, Lemma~\ref{JPlem1637} again yields $E [\zeta
''\underline{Y}^n_{J-}\mid\mathcal{G}_{S^\dagger_{r_J}} ]\to^pE [\zeta
''Y_{J-}\mid{\mathcal{G}^{(0)}_{S_{r_J}-}} ]$. Since
$Y^n_{J-}-\underline{Y}^n_{J-}\to^p0$ and $\zeta''$ and $F_{J-}$ are
bounded, the bounded convergence theorem again implies that
%
\begin{equation}
\label{limitYminus} E \bigl[\zeta''\underline{Y}^n_{J-}
\mid\mathcal{G}_{S^\dagger_{r_J}} \bigr]\to^pE \bigl[
\zeta''Y_{J-}\mid{\mathcal{G}^{(0)}_{S_{r_J-}}}
\bigr].
\end{equation}
Equations (\ref{limitYplus}) and (\ref{limitYminus}) yield (\ref{jointCLTaim3}).


 \textit{Step} 7. Set $\Omega'_n=\{S_{r_{J-1}}+k_n\bar
{r}_n<S^\dagger_{r_J}\}$ if $J>1$ and $\Omega'_n=\Omega$ otherwise.
Then we have
\[
E \Biggl[\zeta f_1\bigl(\check{\mathbf{L}}^n
\bigr)f_2\bigl(\widetilde{\mathcal{W}}(\beta)\bigr)\prod
_{j=1}^JY^n_{j-}Y^n_{j+};
\Omega'_n \Biggr] =E \Biggl[f_1\bigl(\check{
\mathbf{L}}^n\bigr)\prod_{j=1}^{J-1}Y^n_{j-}Y^n_{j+}E
\bigl[\zeta'Y^n_{J-}Y^n_{J+}
\mid\mathcal{G}_{S^\dagger_{r_J}} \bigr];\Omega'_n
\Biggr].
\]
Therefore, by $(\ref{jointaim4})$ and the boundedness of $\zeta'$,
$f_1$ and $F_{j\pm}$ we obtain
\[
E \Biggl[\zeta f_1\bigl(\check{\mathbf{L}}^n
\bigr)f_2\bigl(\widetilde{\mathcal{W}}(\beta)\bigr)\prod
_{j=1}^JY^n_{j-}Y^n_{j+};
\Omega'_n \Biggr] -E \Biggl[\check{\zeta}f_1
\bigl(\check{\mathbf{L}}^n\bigr)\prod_{j=1}^{J-1}Y^n_{j-}Y^n_{j+};
\Omega'_n \Biggr]\to0,
\]
where $\check{\zeta}=E [\zeta'Y_{J-}Y_{J+}\mid{\mathcal
{G}^{(0)}_{S_{r_J-}}} ]$. Since $\lim_nP(\Omega_n')=1$, we conclude that
%
\begin{equation}
\label{jointaim5} E \Biggl[\zeta f_1\bigl(\check{\mathbf{L}}^n
\bigr)f_2\bigl(\widetilde{\mathcal{W}}(\beta)\bigr)\prod
_{j=1}^JY^n_{j-}Y^n_{j+}
\Biggr] -E \Biggl[\check{\zeta}f_1\bigl(\check{\mathbf{L}}^n
\bigr)\prod_{j=1}^{J-1}Y^n_{j-}Y^n_{j+}
\Biggr]\to0.
\end{equation}


 \textit{Step} 8. Now we are ready to prove $(\ref
{jointaim2})$. From $(\ref{jointaim3})$ and $(\ref{jointaim5})$ it
remains to prove
%
\begin{equation}
\label{jointaim6} E \Biggl[\check{\zeta}f_1\bigl(\check{
\mathbf{L}}^n\bigr)\prod_{j=1}^{J-1}Y^n_{j-}Y^n_{j+}
\Biggr] \to E \Biggl[\zeta f_1\bigl(\mathcal{W}^{S_{r_J}}
\bigr)f_2\bigl(\widetilde{\mathcal{W}}(\beta)\bigr)\prod
_{j=1}^JY_{j-}Y_{j+}
\Biggr].
\end{equation}
By the assumption of the induction, we have
\[
\bigl(n^{1/4}\mathbf{L}[M]^n,\bigl(z^n_{r_j-},z'^n_{r_j-},z^n_{r_j+},z'^n_{r_j+}
\bigr)_{j=1}^{J-1}\bigr) \to^{d_s}\bigl(\mathcal{W},
\bigl(z_{r_j-},z'_{r_j-},z_{r_j+},z'_{r_j+}
\bigr)_{j=1}^{J-1}\bigr)
\]
in $\mathbb{D}^{d\times d}\times\mathbb{R}^{2(d'+d)(J-1)}$ as $n\to
\infty$. Therefore, by Lemma~\ref{contstop} and the continuous mapping
theorem we obtain
\[
\bigl(\check{\mathbf{L}}^n,\bigl(z^n_{r_j-},z'^n_{r_j-},z^n_{r_j+},z'^n_{r_j+}
\bigr)_{j=1}^{J-1}\bigr) \to^{d_s}\bigl(
\mathcal{W}^{S_{r_J}},\bigl(z_{r_j-},z'_{r_j-},z_{r_j+},z'_{r_j+}
\bigr)_{j=1}^{J-1}\bigr)
\]
in $\mathbb{D}^{d\times d}\times\mathbb{R}^{2(d'+d)(J-1)}$. This
implies that
$E [\check{\zeta}f_1(\check{\mathbf{L}}^n)\prod
_{j=1}^{J-1}Y^n_{j-}Y^n_{j+} ]
\to E [\check{\zeta}f_1(\mathcal{W}^{S_{r_J}})\*\prod
_{j=1}^{J-1}Y_{j-}Y_{j+} ]$.
Now, since $\zeta f_2(\widetilde{\mathcal{W}}(\beta))$ is independent
of {$Y_{J\pm}$ by construction}, we have $\check{\zeta}=E [\zeta
f_2(\widetilde{\mathcal{W}}(\beta))Y_{J-}Y_{J+}\mid{\mathcal
{G}^{(0)}_{S_{r_J-}}} ]$. Moreover, since $\zeta f_2(\widetilde{\mathcal
{W}}(\beta))Y_{J-}Y_{J+}$ is independent of {$Y_{1\pm},\dots,Y_{(J-1)\pm
}$} by construction, we conclude that
\[
E \Biggl[\check{\zeta}f_1\bigl(\mathcal{W}^{S_{r_J}}\bigr)
\prod_{j=1}^{J-1}Y_{j-}Y_{j+}
\Biggr] =E \Biggl[\zeta f_1\bigl(\mathcal{W}^{S_{r_J}}
\bigr)f_2\bigl(\widetilde{\mathcal{W}}(\beta)\bigr)\prod
_{j=1}^JY_{j-}Y_{j+}
\Biggr],
\]
hence we obtain $(\ref{jointaim6})$.
\end{pf}

\begin{pf*}{Proof of Proposition~\ref{mainCLT}}
First, by Propositions~\ref{appcont}--\ref{jointCLT} as well as
properties of stable convergence, we have
\begin{eqnarray*}
&& \biggl(n^{1/4} \biggl(\bolds{\Xi}\bigl[C(m)\bigr]^n-[M,M]-
\frac{\psi_1}{\psi_2k_n^2}[Y,Y]^{n} \biggr),\bigl(\eta_-(n,r),
\eta'_-(n,r),\eta_+(n,r),\eta'_+(n,r)
\bigr)_{r\geq1} \biggr)
\\
&&\quad \to^{d_s} \bigl(\mathcal{W},\bigl(\sigma_{S_r-}z^{n}_{r-},z'^{n}_{r-},
\sigma_{S_r}z^{n}_{r+},z'^{n}_{r+}
\bigr)_{r\geq1} \bigr).
\end{eqnarray*}
The second claim immediately follows from this convergence.
On the other hand, since equations $(\ref{CJformula})$--$(\ref
{EJformula})$ holds on the set $\Omega_n(t,m)$ (recall that $\Omega
_n(t,m)$ is defined at the beginning of Section~\ref{secappdiscont}),
we obtain
%
\begin{eqnarray}\label{jumplimiteq1}
&& n^{1/4} \biggl(\bolds{\Xi}\bigl[X(m)
\bigr]^n_t- \bigl(\Xi^{(k,l)}_{g,g}
\bigl(J(m),J(m)\bigr)^n_t \bigr)_{1\leq k,l,\leq d}-[M,M]_t-
\frac{\psi_1}{\psi_2k_n^2}[Y,Y]^{n}_t \biggr)
\nonumber\\[-8pt]\\[-8pt]\nonumber
&&\quad  \to^{d_s}
\mathcal{W}_t+\mathcal{Z}(m)_t\quad
\end{eqnarray}
by the continuous mapping theorem and the fact that $(\ref{HJYlem22})$
yields $\lim_{n\to\infty}P(\Omega_n(t,m))=1$.

Next, let $\Omega'_{n}(m)$ be the set on which $\llvert
S_{r_1}-S_{r_2}\rrvert >k_n\bar{r}_n$ for any $r_1,r_2\in\mathcal
{R}_m$ such that $r_1\neq r_2$ and $S_{r_1},S_{r_2}<\infty$. 
Then we have
\begin{eqnarray*}
\Xi^{(k,l)}_{g,g}\bigl(J(m),J(m)\bigr)^n_t
=\frac{1}{\psi_2}\sum_{r\in\mathcal{R}_m\dvt S_r\leq t} \Biggl(
\frac{1}{k_n}\sum_{i=1}^{k_n-1}g^n_ig^n_i
\Biggr)\Delta X^k_{S_r}\Delta X^l_{S_r},
\end{eqnarray*}
on the set $\Omega'_{n}(m)\cap\Omega_n(t,m)$. Since $\frac{1}{k_n}\sum
_{i=1}^{k_n-1}g^n_ig^n_i=\psi_2+\mathrm{O}(k_n^{-1})$ by the Lipschitz
continuity of $g$ and $\lim_nP(\Omega'_{n}(m)\cap\Omega_n(t,m))=1$, we obtain
%
\begin{equation}
\label{jumplimiteq2} n^{1/4} \biggl\{\Xi^{(k,l)}_{g,g}
\bigl(J(m),J(m)\bigr)^n_t-\sum
_{r\in\mathcal{R}_m\dvt S_r\leq t}\Delta X^k_{S_r}\Delta
X^l_{S_r} \biggr\}\to^p0.
\end{equation}

Finally, since
\[
\bigl[X(m)^k,X(m)^l\bigr]_t=\bigl[M^k,M^l\bigr]_t+\sum_{r\in\mathcal{R}_m\dvt S_r\leq
t}\Delta X^k_{S_r}\Delta X^l_{S_r},
\]
$(\ref{jumplimiteq1})$ and $(\ref{jumplimiteq2})$ imply the first claim
of the proposition.
\end{pf*}


\subsection{Proof of Proposition \texorpdfstring{\protect\ref{lemsmalljumps}}{6.3}}

We decompose the target quantity as
%
\begin{eqnarray}\label{decomposesmall1}
\hspace*{-30pt}&& \bolds{\Xi}[X]^{n,kl}_t-\bolds{\Xi}
\bigl[X(m)\bigr]^{n,kl}_t
\nonumber\\[-8pt]
\hspace*{-30pt} \\ [-8pt]\nonumber
\hspace*{-30pt}&&\quad = \bigl\{\Xi^{(k,l)}_{g,g}(X,X)^n-
\Xi^{(k,l)}_{g,g}\bigl(X(m),X(m)\bigr)^n \bigr\}+
\Xi^{(k,l)}_{g,g'}\bigl(Z(m),\mathfrak{E}\bigr)^n +
\Xi^{(l,k)}_{g,g'}\bigl(Z(m),\mathfrak{E}\bigr)^n.
\end{eqnarray}
We start by proving the negligibility of the second and the third terms
in the right-hand side of the above equation, which can be shown by an
easy calculation.

\begin{lem}\label{Esmalljumps}
Under the assumptions of Proposition~\ref{lemsmalljumps}, it holds that
\begin{eqnarray*}
\limsup_{m\to\infty}\limsup_{n\to\infty}P
\bigl(n^{1/4}\bigl\llvert\Xi^{(k,l)}_{g,g'}\bigl(Z(m),
\mathfrak{E}\bigr)^n_t\bigr\rrvert>\eta\bigr)=0
\end{eqnarray*}
for any $t,\eta>0$.
\end{lem}

\begin{pf}
First, since $E_0 [\overline{\mathfrak{E}}(g')^l_i\overline{\mathfrak
{E}}(g')^l_j ]=0$ if $\llvert i-j\rrvert \geq k_n$ and $\llvert E_0
[\overline{\mathfrak{E}}(g')^l_i\overline{\mathfrak{E}}(g')^l_j
]\rrvert \lesssim k_n^{-1}$ by~\textup{[{SA3}]} and the definition of
$\mathfrak{E}$, we have
%
\begin{eqnarray*}
E \bigl[\bigl\llvert n^{1/4}\Xi^{(k,l)}_{g,g'}
\bigl(Z(m),\mathfrak{E}\bigr)^n_t\bigr\rrvert
^2 \bigr] \lesssim\frac{\sqrt{n}}{k_n^2}E \Biggl[\sum
_{i=1}^{N^n_t-k_n+1}\bigl\llvert\overline{Z(m)}(g)^k_i
\bigr\rrvert^2 \Biggr].
\end{eqnarray*}
Next, the definition of $Z(m)$ and the optimal sampling theorem yield
\begin{eqnarray*}
E \Biggl[\frac{\sqrt{n}}{k_n^2}\sum_{i=1}^{N^n_t-k_n+1}
\bigl\llvert\overline{Z(m)}(g)^k_i\bigr\rrvert
^2 \Biggr] &\leq& E \Biggl[\frac{\sqrt{n}}{k_n^2}\sum
_{i=1}^\infty\Biggl\llvert\sum
_{p=0}^{k_n-1}g^n_p
Z(m)^k\bigl(I_{i+p}(t)\bigr)\Biggr\rrvert^2
\Biggr]
\\
&\leq&\frac{\sqrt{n}}{k_n^2}\llVert g\rrVert_\infty\overline{
\gamma}_m E \Biggl[\sum_{i=1}^\infty
\sum_{p=0}^{k_n-1}\bigl\llvert
I_{i+p}(t)\bigr\rrvert\Biggr] \leq\frac{\sqrt{n}}{k_n}t\llVert g\rrVert
_\infty\overline{\gamma}_m,
\end{eqnarray*}
where $\overline{\gamma}_m=\int_{A_m^c}\gamma(z)^2\lambda(\mathrm
{d}z)$. Since $\sqrt{n}/k_n=\mathrm{O}(1)$ as $n\to\infty$ and $\lim_m\overline
{\gamma}_m=0$ by the dominated convergence theorem, we conclude that
\begin{eqnarray*}
\limsup_{m\to\infty}\limsup_{n\to\infty}E \Biggl[
\frac{\sqrt{n}}{k_n^2}\sum_{i=1}^{N^n_t-k_n+1}\bigl
\llvert\overline{Z(m)}(g)^k_i\bigr\rrvert
^2 \Biggr]=0.
\end{eqnarray*}
Therefore, the Chebyshev inequality implies the desired result.
\end{pf}

Next, we prove the negligibility of the term $\Xi
^{(k,l)}_{g,g}(X,X)^n_t-\Xi^{(k,l)}_{g,g}(X(m),X(m))^n_t$. We further
decompose it as
%
\begin{eqnarray}\label{decomposesmall2}
&& \Xi^{(k,l)}_{g,g}(X,X)^n_t-
\Xi^{(k,l)}_{g,g}\bigl(X(m),X(m)\bigr)^n_t
\nonumber\\[-8pt]\\[-8pt]\nonumber
&&\quad
=\Xi^{(k,l)}_{g,g}\bigl(Z(m),X\bigr)^n_t+
\Xi^{(k,l)}_{g,g}\bigl(X,Z(m)\bigr)^n_t-
\Xi^{(k,l)}_{g,g}\bigl(Z(m),Z(m)\bigr)^n_t.
\end{eqnarray}
Therefore, using the decomposition
%
\begin{equation}
\label{decomposeX} X_t=X_0+B'_t+M_t+Z_t,
\end{equation}
where $B'_t=\int_0^tb'_s\,\mathrm{d}s$, $b'_s=b_s+\int_{\{\llVert \delta
(s,z)\rrVert >1\}}\delta(s,z)\lambda(\mathrm{d}z)$ and $Z_t=\delta\star
(\mu-\nu)_t$,
it is enough to prove the negligibility of $\Xi
^{(k,l)}_{g,g}(Z(m),V)^n_t$ for $V\in\{B',M,Z,Z(m)\}$.
In the following we fix $V\in\{B',M,Z,Z(m)\}$.

\begin{lem}\label{edgecont}
Assume~\textup{[{SA2}]}. Then,
%
\begin{longlist}[(a)]
\item[(a)]$\sup_{0\leq h\leq h_0}\llVert V_t-V_{(t-h)_+}\rrVert =\mathrm{O}_p(\sqrt
{h_0})$ as $h_0\downarrow0$,

\item[(b)]$\sup_{0\leq h\leq h_0}|
[Z(m)^k,V^l]_t-[Z(m)^k,V^l]_{(t-h)_+}| =\mathrm{O}_p(h_0)$ as $h_0\downarrow0$ for all $k,l$.
\end{longlist}
\end{lem}

\begin{pf}
The claim is evident if $V=B'$, so we assume that $V\neq B'$. Then the
Doob inequality and~\textup{[{SA2}]} yield $E[\sup_{0\leq h\leq
h_0}\llVert V_t-V_{(t-h)_+}\rrVert ^2]\lesssim h_0$, which implies
(a). On the other hand, the Kunita--Watanabe and Schwarz inequalities
as well as~\textup{[{SA2}]} yield $E[\sup_{0\leq h\leq h_0}\llvert
[Z(m)^k,V^l]_t-[Z(m)^k,V^l]_{(t-h)_+}\rrvert ]\lesssim h_0$, which
implies (b).
\end{pf}

%
\begin{lem}\label{supC}
$\sup_{1\leq q\leq N^n_t+1}\llvert C^n_{g,g}(V)^k_{q}\rrvert =\mathrm{O}_p(1)$
as $n\to\infty$ for any $t>0$ and $k=1,\dots,d$.
\end{lem}

\begin{pf}
This can be shown in the same manner as the proof of Lemma 6.8 from
\cite{Koike2014time}.
\end{pf}

%
\begin{lem}\label{XirepZ}
Under the assumptions of Proposition~\ref{lemsmalljumps}, it holds that
\begin{eqnarray*}
n^{1/4} \bigl\{\Xi^{(k,l)}_{g,g}\bigl(Z(m),V
\bigr)^n_t-\mathbb{L}^{(k,l)}_{g,g}
\bigl(Z(m),V\bigr)^n_t-\bigl[Z(m)^k,V^l
\bigr]_t \bigr\}\to^p0
\end{eqnarray*}
as $n\to\infty$ for any $t>0$.
\end{lem}

\begin{pf}
Simple calculations and Lemma~\ref{edgecont}(a) yield
\[
\mathbf{B}:=\Xi^{(k,l)}_{g,g}\bigl(Z(m),V\bigr)^n_t-
\mathbb{L}^{(k,l)}_{g,g}\bigl(Z(m),V\bigr)^n_t=
\sum_{p=1}^{N^n_t+1}c^n_{g,g}(p,p)Z(m)^k(I_{p})_tV^l(I_{p})+\mathrm{o}_p
\bigl(n^{-1/4}\bigr).
\]
Therefore, using Lemma~\ref{edgecont}(b), we can prove $\mathbf
{B}=[Z(m)^k,V^l]_t+\mathrm{o}_p(n^{-1/4})$ analogously to the proof of equation
(6.24) from \cite{Koike2014time}.
\end{pf}

%
\begin{lem}\label{Lsmalljumps}
Under the assumptions of Proposition~\ref{lemsmalljumps}, it holds that
%
\begin{equation}
\label{Lsmalljumpsaim} \limsup_{m\to\infty}\limsup_{n\to\infty}P
\bigl(n^{1/4}\bigl\llvert\mathbb{L}^{(k,l)}_{g,g}
\bigl(Z(m),V\bigr)^n_t\bigr\rrvert>\eta\bigr)=0
\end{equation}
for any $t,\eta>0$.
\end{lem}

\begin{pf}
First, if $V=B'$, we can adopt an analogous argument to the proof of
Lemma 6.10 from~\cite{Koike2014time} and deduce $n^{1/4}\mathbb
{L}^{(k,l)}_{g,g}(Z(m),V)^n_t\to^p0$ as $n\to\infty$ for every $m$, so
$(\ref{Lsmalljumpsaim})$ holds true.

Next, we suppose that $V\neq B'$. It is enough to prove
%
\begin{eqnarray}
\limsup_{m\to\infty}\limsup_{n\to\infty}P
\bigl(n^{1/4}\bigl\llvert\mathbb{M}^{(k,l)}_{g,g}
\bigl(Z(m),V\bigr)^n_t\bigr\rrvert>\eta'
\bigr)=0,\label{Msmalljumpsaim1}
\\
\limsup_{m\to\infty}\limsup_{n\to\infty}P
\bigl(n^{1/4}\bigl\llvert\mathbb{M}^{(k,l)}_{g,g}
\bigl(V,Z(m)\bigr)^n_t\bigr\rrvert>\eta'
\bigr)=0\label{Msmalljumpsaim2}
\end{eqnarray}
for any $\eta'>0$. Since we can prove $(\ref{Msmalljumpsaim2})$ in a
similar manner to the proof of $(\ref{Msmalljumpsaim1})$, we only prove
$(\ref{Msmalljumpsaim1})$.

Since $V$ is an {$(\mathcal{F}^{(0)}_t)$}-martingale for any $n$ due to
\textup{[{SA2}]}, by the Lenglart inequality it suffices to show that
%
\begin{equation}
\label{eqU1} \limsup_{m\to\infty}\limsup_{n\to\infty}P
\bigl(U(n,m)_t>\eta' \bigr)=0
\end{equation}
for any $\eta'>0$, where $U(n,m)_t=\sqrt{n}\sum
_{q=2}^{N^n_t+1}E[\llvert C^n_{g,g}(Z(m))^k_q V^l(I_q)\rrvert ^2\mid
{\mathcal{F}^{(0)}_{T_{q-1}}}]$. To prove $(\ref{eqU1})$, for each
$j\geq1$ we set $\Lambda(j)^n_q=\{E[n\llvert I_q\rrvert \mid{\mathcal
{F}^{(0)}_{T_{q-1}}}]\leq j\}$ and decompose $U(n,m)$ as
\begin{eqnarray*}
U(n,m)_t &=&\sqrt{n}\sum_{q=2}^{N^n_t+1}E
\bigl[\bigl\llvert C^n_{g,g}\bigl(Z(m)\bigr)^k_q
V^l(I_q)\bigr\rrvert^2\mid{
\mathcal{F}^{(0)}_{T_{q-1}}} \bigr] (1_{\Lambda(j)^n_q}+1_{(\Lambda(j)^n_q)^c}
)
\\
&=&U(n,m,j)_t+U'(n,m,j)_t.
\end{eqnarray*}

First, we prove
%
\begin{equation}
\label{eqU2} \limsup_{m\to\infty}\limsup_{n\to\infty}P
\bigl(U(n,m,j)_t>\eta' \bigr)=0
\end{equation}
for any fixed $j$ and any $\eta'>0$. We have
\begin{eqnarray*}
E \bigl[U(n,m,j)_{t} \bigr] &=&\sqrt{n}E \Biggl[\sum
_{q=2}^{N^n_{t}+1}\bigl\llvert C^n_{g,g}
\bigl(Z(m)\bigr)^k_q\bigr\rrvert^2E \bigl[
\bigl\langle V^l\bigr\rangle(I_q)\mid{
\mathcal{F}^{(0)}_{T_{p-1}}} \bigr]1_{\Lambda(j)^n_q} \Biggr]
\\
&\lesssim&\frac{j}{\sqrt{n}}E \Biggl[\sum
_{q=2}^\infty\sum_{p=(q-k_n)\vee1}^{q-1}
\bigl\llvert c^n_{u,v}(p,q)\bigr\rrvert^2\bigl
\langle Z(m)^k\bigr\rangle\bigl(I_p(t)\bigr) \Biggr]
\\
&\lesssim&\overline{\gamma}_m\frac{1}{\sqrt{n}}E \Biggl[\sum
_{q=2}^\infty\sum_{p=(q-k_n)\vee1}^{q-1}
\bigl\llvert I_p(t)\bigr\rrvert\Biggr] \leq\overline{
\gamma}_mk_nn^{-1/2}t,
\end{eqnarray*}
hence it holds that $\limsup_m\limsup_nE [U(n,m,j)_{t} ]=0$. Therefore,
we obtain $(\ref{eqU2})$ by the Chebyshev inequality.

Next, we prove
%
\begin{equation}
\label{eqU3} \limsup_{j\to\infty}\limsup_{m\to\infty}
\limsup_{n\to\infty}P \bigl(U'(n,m,j)_t>
\eta' \bigr)=0
\end{equation}
for any $\eta'>0$. Since $E[n\llvert I_q\rrvert \mid{\mathcal
{F}^{(0)}_{T_{q-1}}}]=E[E[n\llvert I_q\rrvert \mid{\mathcal
{G}^{(0)}_{T_{q-1}}}]\mid{\mathcal{F}^{(0)}_{T_{q-1}}}]$, we have
$\Lambda(j)^n_q\supset\{E[n\llvert I_q\rrvert \mid{\mathcal
{G}^{(0)}_{T_{q-1}}}]\leq j\}$. Therefore, 
\begin{eqnarray*}
U'(n,m,j)_t 
&\leq&\sqrt{n}\sum
_{q=2}^{N^n_t+1}E \bigl[\bigl\llvert C^n_{g,g}
\bigl(Z(m)\bigr)^k_q V^l(I_q)
\bigr\rrvert^2\mid{\mathcal{F}^{(0)}_{T_{q-1}}}
\bigr] (1_{\{G^n_{T_{q-1}}>j\}}+1_{\{q-1\in\mathcal{N}^n\}} )
\\
&=:&U'(n,m,j)_t^{(1)}+U'(n,m,j)_t^{(2)}.
\end{eqnarray*}
Since $\{U'(n,m,j)^{(1)}>0\}\subset\{\sup_{0\leq s\leq t}G^n_{s}>j\}$,
we have
$\limsup_{j}\limsup_{m}\limsup_{n}P (U'(n,\break m, j)^{(1)}_t>0 )=0$.
On the other hand,~\textup{[{SA2}]} and $(\ref{SA4})$ imply that
\begin{eqnarray*}
U'(n,m,j)_t^{(2)} 
&\leq&\sqrt{n}
\bar{r}_n\sup_{1\leq q\leq N^n_t+1}\bigl\llvert
C^n_{g,g}\bigl(Z(m)\bigr)^k_q\bigr
\rrvert^2\#\bigl(\mathcal{N}^n\cap\{q\dvt T_q
\leq t\}\bigr),
\end{eqnarray*}
hence Lemma~\ref{supC}, $(\ref{estxi})$ and~\textup{[A1]}(i) yield
$\limsup_nP(U'(n,m,j)^{(2)}_t>\eta')=0$. Consequently, we obtain $(\ref{eqU3})$.

From $(\ref{eqU2})$ we have $\limsup_{m}\limsup_{n}P (U(n,m)_t>\eta'
)\leq\limsup_{m}\limsup_{n}P (U'(n,\break m,j)_t>\eta' )$ for any $j\geq1$ and
any $\eta'>0$. Hence, $(\ref{eqU3})$ yields $(\ref{eqU1})$, which
completes the proof.
\end{pf}

\begin{pf*}{Proof of Proposition~\ref{lemsmalljumps}}
$(\ref{eqsmalljumps2})$ immediately follows from equations $(\ref
{decomposeX})$ and $(\ref{decomposesmall1})$--$(\ref{decomposesmall2})$
as well as Lemmas~\ref{Esmalljumps} and~\ref{XirepZ}--\ref
{Lsmalljumps}. Equation~(\ref{eqsmalljumps1}) follows from the equation
$\widetilde{E} [\llvert \mathcal{Z}(m)^{kl}_t-\mathcal
{Z}^{kl}_t\rrvert ^2\mid{\mathcal{F}^{(0)}} ]
=\frac{1}{\psi_2^2}\sum_{r\notin\mathcal{R}_m\dvt S_r\leq t} (\mathfrak
{J}^{kk}_{S_r}+2\mathfrak{J}^{kl}_{S_r}+\mathfrak{J}^{ll}_{S_r} )$
and the fact that $\sum_{r\notin\mathcal{R}_m\dvt S_r\leq t}\llVert
\Delta X_{S_r}\rrVert ^2\to^p0$ as $m\to\infty$.
\end{pf*}


\section*{Acknowledgements}
I wish to thank the Editors and an anonymous referee for their careful
reading and constructive comments that substantially improved this paper.
I am also grateful to Teppei Ogihara who pointed out a problem on the
mathematical construction of the noise process in a draft of this paper.
This work was supported by Grant-in-Aid for JSPS Fellows.


%

\printhistory
\end{document}